\documentclass[final,onefignum,onetabnum]{siamonline250211}
\usepackage[scale=0.8]{geometry} 

\usepackage{lipsum}
\usepackage{amsfonts}
\usepackage{epstopdf}
\ifpdf
  \DeclareGraphicsExtensions{.eps,.pdf,.png,.jpg}
\else
  \DeclareGraphicsExtensions{.eps}
\fi

\usepackage{datetime}
\newdateformat{monthyeardate}{%
  \monthname[\THEMONTH] \THEDAY, \THEYEAR}

\usepackage{academicons}
\usepackage{xcolor}

\usepackage{amsmath}
\allowdisplaybreaks
\usepackage{amssymb}
\usepackage{commath}
\usepackage{mathtools}
\usepackage{bbm}

\usepackage{color}
\usepackage{graphicx}
\usepackage[small]{caption}
\usepackage{subcaption}

\usepackage{relsize}
\usepackage{adjustbox}
\usepackage{algorithm}
\usepackage[noend]{algpseudocode}
\usepackage{booktabs}
\usepackage{tikz}

\usepackage{verbatim}

\usepackage{enumitem}
\setlist[enumerate]{leftmargin=.5in}
\setlist[itemize]{leftmargin=.5in}

\newsiamthm{problem}{Problem}
\newsiamremark{remark}{Remark}
\newsiamremark{hypothesis}{Hypothesis} 
\crefname{hypothesis}{Hypothesis}{Hypotheses}
\newsiamremark{example}{Example}
\newsiamthm{claim}{Claim}
\newsiamthm{conjecture}{Conjecture}

\headers{Efficient sampling for SBL using prior normalization}{J.\ Glaubitz and Y.\ Marzouk}

\title{Efficient sampling for sparse Bayesian learning using hierarchical prior normalization%
\thanks{
Submitted to the editors \monthyeardate\today 
\funding{
JG and YM acknowledge support from the US Department of Defense, MURI program under grant number \#N00014-20-1-2595 and from the US Department of Energy, SciDAC program (ASCR/HEP) under grant number \#DE-SC0012704.
JG further acknowledges support by the Swedish Research Council (VR) Starting Grant \#2025-05370, the Zenith Career Development Grant \#26.07, and the National Academic Infrastructure for Supercomputing in Sweden (NAISS) grants \#2025/22-1599 and \#2024/22-1207.
}}
}

\author{
Jan Glaubitz\thanks{Department of Mathematics, Link\"oping University, Sweden (\email{jan.glaubitz@liu.se})}
\and 
Youssef M.\ Marzouk\thanks{Department of Aeronautics and Astronautics, Massachusetts Institute of Technology, USA (\email{ymarz@mit.edu})} 
}

\usepackage{amsopn}

\DeclareMathOperator{\diag}{diag}

\newcommand{\intd}{\, \mathrm{d}}

\newcommand{\R}{\mathbb{R}}

\usepackage{lineno}

\ifpdf
\hypersetup{
  pdftitle={Hierarchical prior normalization},
  pdfauthor={Jan Glaubitz}
}
\fi

\begin{document}

\maketitle

\begin{abstract}
	We introduce an approach for efficient Markov chain Monte Carlo (MCMC) sampling for challenging high-dimensional distributions in sparse Bayesian learning (SBL). 
The core innovation involves using \emph{hierarchical prior-normalizing transport maps (TMs)}, which are deterministic couplings that transform the sparsity-promoting SBL prior into a standard normal one. 
We analytically derive these prior-normalizing TMs by leveraging the product-like form of SBL priors and Knothe--Rosenblatt (KR) rearrangements.
These transform the complex target posterior into a simpler reference distribution equipped with a  \emph{standard normal prior} that can be sampled more efficiently. 
Specifically, one can leverage the standard normal prior by using more efficient, structure-exploiting samplers.  
Our numerical experiments on various inverse problems---including signal deblurring, inverting the non-linear inviscid Burgers equation, and recovering an impulse image---demonstrate significant performance improvements for standard MCMC techniques.
\end{abstract}



\begin{Code}
    \url{https://github.com/jglaubitz/paper-2025-SBL-priorNormalization}
\end{Code}

\begin{DOI}
	\url{https://doi.org/10.1137/25M1790427}
\end{DOI}

\section{Introduction} 
\label{sec:introduction} 


Many applications involve recovering an unknown high-dimensional parameter vector $\mathbf{x} \in \R^n$ from noisy, indirect, and limited observational data $\mathbf{y} \in \R^m$ by solving an inverse problem. 
For instance, for additive noise, one considers the data model   
\begin{equation}\label{eq:data_model} 
	\mathbf{y} = F(\mathbf{x}) + \mathbf{e}, 
\end{equation}
where $F: \R^n \to \R^m$ is a known forward operator and $\mathbf{e} \in \R^m$ denotes an unknown additive noise component. 
The Bayesian approach \cite{stuart2010inverse,calvetti2023bayesian} frames \cref{eq:data_model} as a statistical inference problem based on the posterior distribution $\pi^{y}$, which combines the likelihood function $f( \mathbf{x}; \mathbf{y} )$ implied by \cref{eq:data_model} with a prior density $\pi^0$ that encodes our structural beliefs about $\mathbf{x}$. 

One can often assume that $\mathbf{x}$ is sparse or has a sparse representation---for instance, in the edge domain when $\mathbf{x}$ contains the nodal values of a piecewise smooth function. 
There are various sparsity-promoting priors, including Laplace \cite{figueiredo2007majorization,babacan2009bayesian}, 
Cauchy \cite{markkanen2019cauchy,suuronen2022cauchy}, 
and horseshoe \cite{carvalho2009handling,uribe2023horseshoe,dong2023inducing} priors. 
Here, we consider the particularly potent class of hierarchical SBL priors $\pi^0( \mathbf{x}, \boldsymbol{\theta} ) = \pi^0( \mathbf{x} | \boldsymbol{\theta} ) \, \pi^0( \boldsymbol{\theta} )$. 
These combine a conditionally Gaussian prior $\pi^0( \mathbf{x} | \boldsymbol{\theta} )$ with a generalized gamma hyper-prior $\pi^0( \boldsymbol{\theta} )$. 
SBL priors have been demonstrated to provide efficient algorithms for finding sparse solutions to ill-posed inverse problems.  
See \cite{tipping2001sparse,calvetti2007gaussian} for early versions of SBL based on (inverse) gamma hyper-priors, \cite{calvetti2020sparse,calvetti2020sparsity} for their extension to generalized gamma hyper-priors, \cite{calvetti2019hierachical,sanz2024hierarchical} for their analysis, and \cite{glaubitz2023generalized,xiao2023sequential,glaubitz2024leveraging,lindbloom2024generalized,lindbloom2025priorconditioned} for generalized SBL (GSBL) models for promoting linear transforms, $R \mathbf{x}$, to be sparse. 
The SBL posterior density $\pi^y$ for $\mathbf{x}, \boldsymbol{\theta}| \mathbf{y}$ is given by Bayes' theorem as 
\begin{equation}\label{eq:posterior}
	\pi^y( \mathbf{x}, \boldsymbol{\theta} ) 
		= \frac{1}{Z} f( \mathbf{x}; \mathbf{y} ) \pi^0( \mathbf{x}, \boldsymbol{\theta} ),
\end{equation}
where $Z$ is a normalizing constant that is generally unknown.
Most works on SBL have concentrated on maximum a posteriori (MAP) estimation. 
However, MAP estimates often do not adequately describe the posterior \cref{eq:posterior}, do not allow quantifying uncertainty, and can be unstable w.r.t.\ perturbations in $\mathbf{y}$, particularly in the case of multimodal posterior densities. 

In this work, we aim to characterize the SBL posterior by evaluating its moments or by computing the probability of an event of interest, which can be cast as computing expectations under the posterior. 
The workhorse algorithms in this setting are sampling methods, with MCMC \cite{brooks2011handbook,liu2013monte,robert2013monte,sanz2024first} being among the most broadly useful. 
They provide a flexible approach for generating correlated samples, using only evaluations of an unnormalized probability density, which can nonetheless be used to compute the desired expectations. 
MCMC algorithms become inefficient, however, if the correlations between successive samples decay too slowly.
SBL priors, in particular, lead to posterior distributions that are challenging to sample for the following reasons: 
(1) The unknown parameter vector $\mathbf{x} \in \R^n$ is often high-dimensional, and the introduction of the hyper-parameter vector $\boldsymbol{\theta} \in \mathbb{R}^n$ within the SBL prior further increases the dimensionality, yielding high per-sample computational cost.
(2) SBL priors typically result in posteriors that are non-log-concave \cite{calvetti2020sparse,lindbloom2024generalized,glaubitz2024leveraging,si2024path} and that can have multiple modes separated by low-density regions, which makes it difficult for samplers to explore the posterior distribution efficiently. 
As we illustrate in \Cref{sec:SBL}, such complexities arise even in relatively simple scenarios, such as when the forward operator is linear and the noise is Gaussian additive.
(3) 
SBL posteriors often exhibit strong correlations between the $\mathbf{x}$ and $\boldsymbol{\theta}$, which can severely impair the sampler's mixing and slow down convergence.
To address the above challenges, \emph{we propose using hierarchical prior-normalizing TMs as an alternative to direct sampling from the SBL posterior}.

\subsection*{Our approach: Hierarchical prior normalization}

The convergence of MCMC algorithms is primarily influenced by the posterior landscape, which is impacted by the prior. 
Hence, we use prior-normalizing TMs \cite{cui2022prior}, transforming the challenging SBL prior $\pi^0(\mathbf{x},\boldsymbol{\theta})$ into a standard normal prior with density $\phi^0(\mathbf{u},\boldsymbol{\tau}) \propto \exp( -[ \| \mathbf{u} \|_2^2 + \| \boldsymbol{\tau} \|_2^2 ] / 2 )$. 
This involves finding a bijective map $S: (\mathbf{x},\boldsymbol{\theta}) \mapsto ( \mathbf{u}, \boldsymbol{\tau} )$ such that $( \mathbf{u}, \boldsymbol{\tau} ) = S(\mathbf{x},\boldsymbol{\theta})$ follows $\phi^0$ whenever $(\mathbf{x},\boldsymbol{\theta})$ follows $\pi^0$.
In other words: the standard normal density $\phi^0$ is the pushforward of the prior density $\pi^0$ under the map $S$, i.e., $\phi^0 = S_{\sharp} \pi^0$; see \cite{villani2009optimal,marzouk2016sampling,parno2018transport}.
Because $S$ is a bijection, we can push forward the challenging original posterior $\pi^y(\mathbf{x},\boldsymbol{\theta})$ to obtain a simpler reference posterior $\phi^y(\mathbf{u},\boldsymbol{\tau})$ equipped with a standard normal prior. 
The resulting \emph{prior-normalized posterior}, $\phi^y = S_{\sharp} \pi^y$, can be expressed as 
\begin{equation}\label{eq:posterior_reference}
	\phi^y( \mathbf{u}, \boldsymbol{\tau} ) 
		= \frac{1}{Z} g( \mathbf{u}, \boldsymbol{\tau}; \mathbf{y} ) \phi^0( \mathbf{u}, \boldsymbol{\tau} ), 
\end{equation} 
where $g( \mathbf{u}, \boldsymbol{\tau}; \mathbf{y} ) := f( S^{-1}(\mathbf{u}, \boldsymbol{\tau}); \mathbf{y} )$ is the pushforward likelihood. 
We can then implement MCMC algorithms with \cref{eq:posterior_reference} as the target distribution to obtain posterior samples in the reference coordinates and transform them via $S^{-1}$ to obtain posterior samples in the original coordinate.
In particular, the change of variables formula 
\begin{equation}
	\int h( \mathbf{x}, \boldsymbol{\theta} ) \, \pi^y( \mathbf{x}, \boldsymbol{\theta} ) \intd(\mathbf{x}, \boldsymbol{\theta}) 
		= \int (h \circ S^{-1})( \mathbf{u},\boldsymbol{\tau} ) \, \phi^y( \mathbf{u},\boldsymbol{\tau} ) \intd(\mathbf{u},\boldsymbol{\tau})
\end{equation}
allows for the direct computation of posterior expectations \cite{athreya2006measure}, where $h$ is a generic test function. 
Notably, the prior $\phi^0$ in the reference coordinate is a standard normal one, and we observe the prior-normalized posterior $\phi^y$ to be easier to explore with standard MCMC methods---which we demonstrate for various examples in \Cref{sec:tests}.  
Specifically, we can utilize the prior-normalized posterior having a standard normal prior by using more efficient, structure-exploiting samplers, such as the elliptical slice sampler \cite{murray2010elliptical}.
Moreover, we analytically derive the desired hierarchical prior-normalizing TMs by leveraging the product-like form of the SBL prior and considering KR \cite{rosenblatt1952remarks,knothe1957contributions,bogachev2005triangular} rearrangements, which are a particular class of TMs. 
Numerical demonstrations on various inverse problems, including signal deblurring, inverting the nonlinear inviscid Burgers equation, and recovering an impulse image, show that standard MCMC techniques sample the prior-normalized posterior more efficiently than the original one. 
Specifically, even in settings where the Gibbs sampler---widely used in hierarchical Bayesian models due to its relatively low per-sample cost \cite{tan2010efficient,markkanen2019cauchy,churchill2022sampling,uribe2022hybrid,uribe2023horseshoe}---is applicable, we find that it may fail to explore the original posterior adequately. 
In contrast, existing MCMC samplers exhibit substantially improved exploration efficiency when applied to the proposed prior-normalized posterior.

\subsection*{Related works}

Hierarchical posteriors are often explored using (block) Gibbs or Metropolis-within-Gibbs (MwG) algorithms that alternate between sampling from the conditional densities $\pi^y( \mathbf{x} | \boldsymbol{\theta} ) \propto f( \mathbf{x}; \mathbf{y} ) \, \pi^{0}( \mathbf{x} | \boldsymbol{\theta} )$ and $\pi^y( \boldsymbol{\theta} | \mathbf{x} ) \propto \pi^{0}( \mathbf{x} | \boldsymbol{\theta} ) \, \pi^{0}( \boldsymbol{\theta} )$. 
This strategy has been analyzed in \cite{agapiou2014analysis,ascolani2024scalability} and applied to sparsity-promoting hierarchical priors in  \cite{tan2010efficient,markkanen2019cauchy,churchill2022sampling,uribe2022hybrid,uribe2023horseshoe}. 
Notably, Gibbs and MwG are so-called centered algorithms \cite{papaspiliopoulos2003non,papaspiliopoulos2007general,chada2018parameterizations,dunlop2020hyperparameter}, which mix poorly when $\pi^0( \mathbf{x}, \boldsymbol{\theta} )$ exhibits strong correlations. 
One approach to mitigate this pathology is to use a non-centered parameterization \cite[Section 4]{agapiou2014analysis}, transforming the conditional prior $\pi^0( \mathbf{x} | \boldsymbol{\theta} )$ into a standard normal one. 
In contrast, the proposed prior-normalizing TMs normalize the joint prior $\pi^0( \mathbf{x}, \boldsymbol{\theta} )$ rather than just $\pi^0( \mathbf{x} | \boldsymbol{\theta} )$. 

Existing works that exploit prior normalization to accelerate MCMC for problems with heavy-tailed priors, such as SBL priors, include \cite{fleischer2007transformations,wang2017bayesian,chen2018robust,cui2022prior}. 
Specifically, \cite{fleischer2007transformations} discusses sampling a one-dimensional multimodal distribution, \cite{cui2022prior} explores prior normalization in the context of dimension reduction, and \cite{wang2017bayesian,chen2018robust} implement the preconditioned Crank-Nicholson (pCN) and randomize-then-optimize algorithms with prior normalization. 
Notably, none of these works considered prior normalization for SBL or other \emph{hierarchical} sparsity-promoting priors.

Finally, \cite{calvetti2024computationally} combined a (modified) pCN algorithm with partial prior normalization for SBL. 
Specifically, the authors decomposed the joint prior into the product of an exponential term and a power function, applying a transformation to normalize only the exponential component to a standard normal density. 
This approach represents a middle ground between traditional non-centered parameterizations and full prior normalization as proposed here.

\subsection*{Outline} 

We begin with some preliminaries on SBL in \Cref{sec:SBL}.
In \Cref{sec:transport}, we introduce and analytically derive the proposed hierarchical prior-normalizing TMs.
\Cref{sec:tests} presents numerical experiments comparing the performance of various MCMC samplers applied to both the original and prior-normalized posteriors.
Finally, we conclude in \Cref{sec:summary}. 
\section{Background on sparse Bayesian learning} 
\label{sec:SBL} 

\begin{figure}[tb]
	\centering 
	\begin{subfigure}[b]{0.45\textwidth}
		\includegraphics[width=\textwidth]{%
      		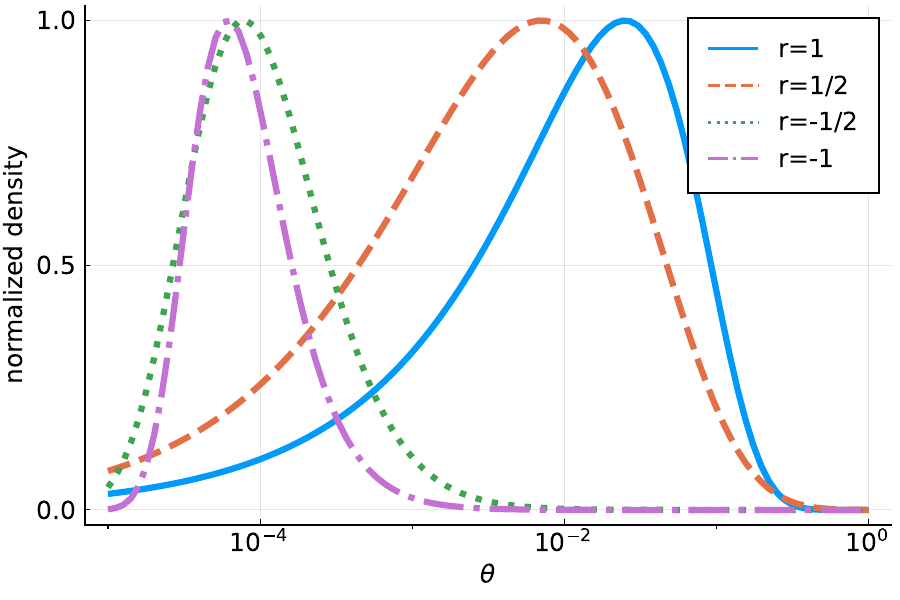} 
    		\caption{Generalized gamma hyper-prior}
    		\label{fig:PDF_hyper_prior}
  	\end{subfigure}%
	~
	\begin{subfigure}[b]{0.45\textwidth}
		\includegraphics[width=\textwidth]{%
      		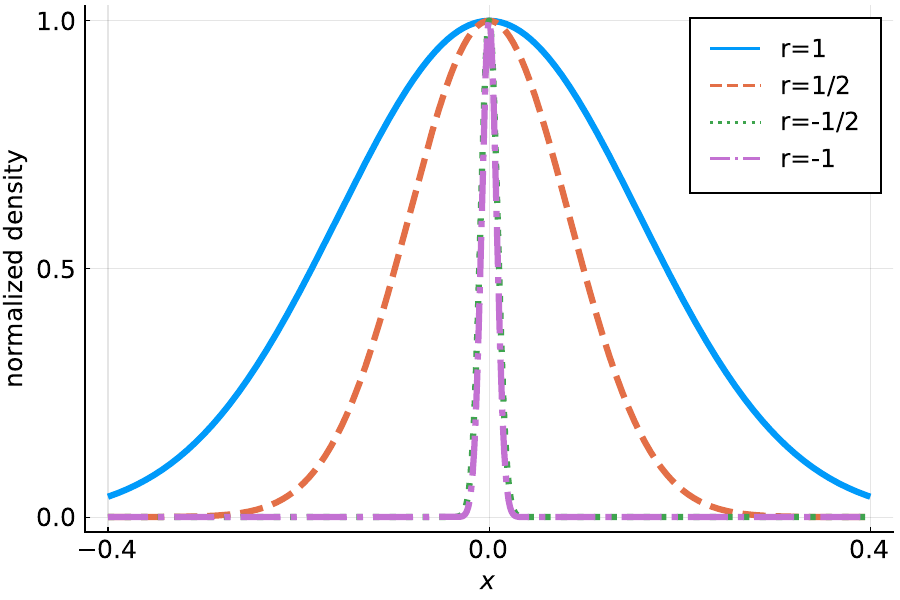} 
    		\caption{Conditionally Gaussian priors}
    		\label{fig:PDF_cond_prior}
  	\end{subfigure}%
  	\caption{
  	Normalized densities of the univariate generalized gamma hyper-prior $\pi( \theta )$ and the conditionally Gaussian prior $\pi( x | \theta )$. 
	The hyper-prior is illustrated for the parameters in \cref{tab:parameters} and the conditional prior for $\theta$ equal to the mode (maximum value) of $\pi( \theta )$, which is $\max\{ 0, \vartheta ( [ r\beta - 1 ]/r )^{1/r} \}$. 
  	}
  	\label{fig:prior_model}
\end{figure}

We review the SBL approach and illustrate the challenges associated with sampling for a one-dimensional toy problem. 
The SBL prior is 
\begin{equation}\label{eq:prior_model}
	\pi^0( \mathbf{x}, \boldsymbol{\theta} ) 
		= \pi^0( \mathbf{x} | \boldsymbol{\theta} ) \, 
		\pi^0( \boldsymbol{\theta} )
\end{equation}
with conditional Gaussian prior $\pi^0( \mathbf{x} | \boldsymbol{\theta} )$ and generalized gamma hyper-prior $\pi^0( \boldsymbol{\theta} )$.  
Specifically, SBL assumes that the components of $\mathbf{x}$ are independently normal-distributed with mean zero and variances $\boldsymbol{\theta}$, i.e., $x_i | \theta_i \sim \mathcal{N}(0,\theta_i)$.
The variance $\theta_i$ is also modeled as a random variable, 
which is generalized gamma-distributed, i.e., $\theta_i \sim \mathcal{GG}(r,\beta,\vartheta)$ with parameters $r \in \R \setminus \{0\}$, $\beta > 0$, and $\vartheta > 0$; see \cite{gomes2008parameter}.  
Hence, the conditional prior and hyper-prior are
\begin{equation}\label{eq:SBL_priors}
\begin{aligned} 
	\pi^0( \mathbf{x} | \boldsymbol{\theta} ) 
		\propto \left( \prod_{i=1}^n \theta_i^{-1/2} \right)
		\exp\left( 
			- \sum_{i=1}^n \frac{x_i^2}{2 \theta_i}
		\right), \
	\pi^0( \boldsymbol{\theta} ) 
		\propto \left( \prod_{i=1}^n \theta_i^{r \beta -1} \right)
		\exp\left( 
			- \sum_{i=1}^n \left[ \frac{\theta_i}{\vartheta} \right]^r
		\right).
\end{aligned}
\end{equation}
\cref{tab:parameters} lists a few common parameter combinations used in \cite{calvetti2024computationally}.
Furthermore, \cref{fig:prior_model} illustrates the density functions of the univariate generalized gamma hyper-prior and conditionally Gaussian prior for different parameters. 
We recall that the parameter $r$ smoothly interpolates between different distribution classes.
For example, the choice $r=1$ yields a gamma distribution, whereas $r=-1$ leads to an inverse gamma distribution.
More generally, $r$ governs the rate at which probability mass decays for large arguments and influences the log-concavity properties of the resulting posterior.
See \cite{calvetti2020sparse} for more details.

If we consider the data model \cref{eq:data_model} with $\mathbf{e} \sim \mathcal{N}(\mathbf{0}, \Gamma_{\rm noise})$, then the likelihood function is $f(\mathbf{x};\mathbf{y}) \propto \exp( - \| \Gamma_{\rm noise}^{-1/2} ( F(\mathbf{x}) - \mathbf{y} ) \|_2^2 / 2 )$ and the SBL posterior density \cref{eq:posterior} becomes 
\begin{equation}\label{eq:posterior_expl}
	\pi^y(\mathbf{x},\boldsymbol{\theta}) 
		\propto \exp\left( 
			-\frac{1}{2} \norm{ \Gamma_{\rm noise}^{-1/2} \left[ F(\mathbf{x}) - \mathbf{y} \right] }_2^2 
			- \sum_{i=1}^n \left[ 
				\frac{x_i^2}{2 \theta_i} 
				+ \left( \frac{\theta_i}{\vartheta} \right)^r 
				+ (r \beta - 3/2) \log \theta_i 
			\right]
		\right).
\end{equation}
We illustrate the geometry of the SBL posterior density \cref{eq:posterior_expl} for the following toy problem.

\begin{table}[tb]
\renewcommand{\arraystretch}{1.3}
	\centering 
  	\begin{tabular}{c | c c c c} 
		\toprule 
    		$r$ &$1$ & $1/2$ & $-1/2$ & $-1$ \\  
		$\beta$& $1.501$ & $3.0918$ & $2.0165$ & $1.0017$ \\
    		$\vartheta$ & $5 \cdot 10^{-2}$ & $5.9323 \cdot 10^{-3}$ & $1.2583 \cdot 10^{-3}$ & $1.2308 \cdot 10^{-4}$ \\  
		\bottomrule
	\end{tabular} 
	\caption{
	Typical parameter choices for the generalized gamma hyper-prior
	}
	\label{tab:parameters}
\end{table}

\begin{example}[Toy problem]\label{expl:toy_problem}
	Consider the scalar data model $y = x + e$ with $e \sim \mathcal{N}(0,\sigma^2)$. 
	In this case, \cref{eq:posterior_expl} reduces to $\pi^y( x, \theta ) \propto \theta^{r \beta - 3/2} \exp( -[x-y]^2/[2 \sigma^2] - x^2/[2 \theta] - [ \theta/\vartheta ]^r )$. 
\end{example}

\begin{figure}[tb]
	\centering 
	\begin{subfigure}[b]{0.245\textwidth}
		\includegraphics[width=\textwidth]{%
      		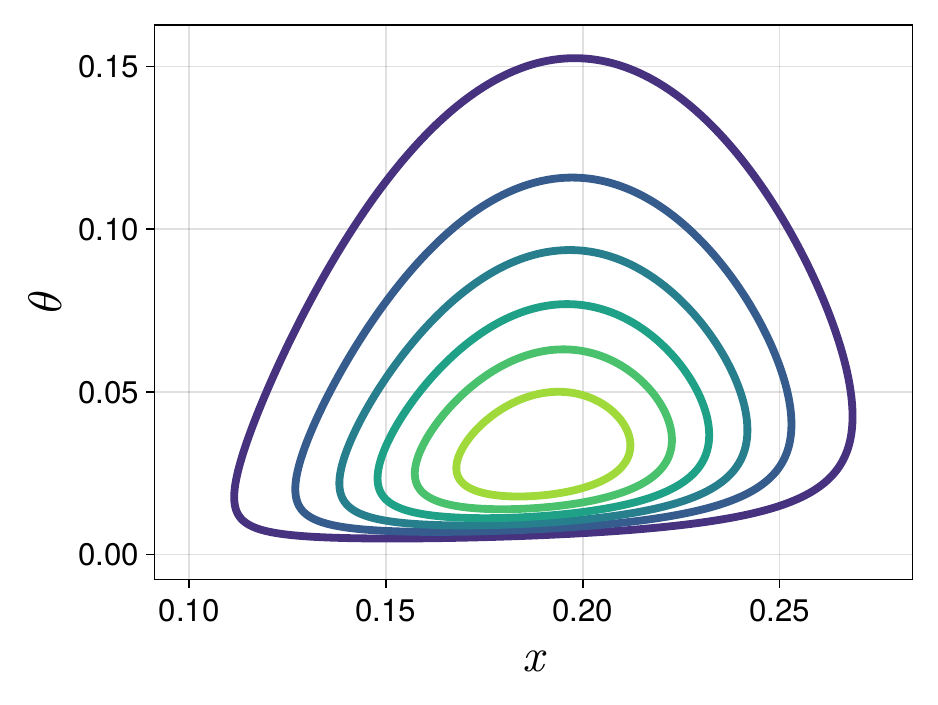} 
    		\caption{$r=1$}
    		\label{fig:contour_posterior_toy_rp1}
  	\end{subfigure}%
	\begin{subfigure}[b]{0.245\textwidth}
		\includegraphics[width=\textwidth]{%
      		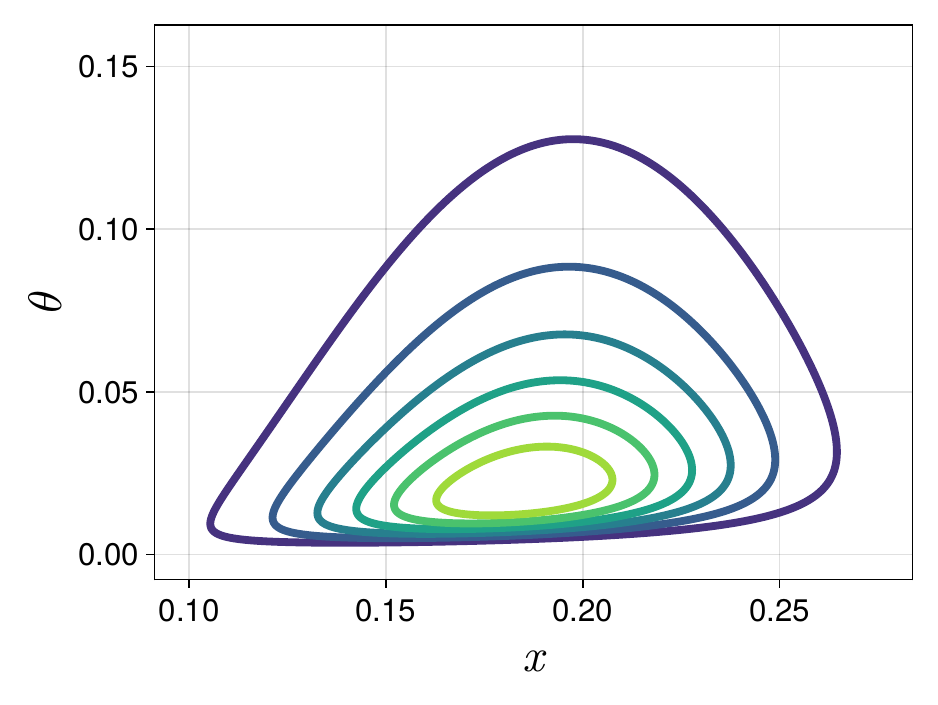} 
    		\caption{$r=1/2$}
    		\label{fig:contour_posterior_toy_rp05}
  	\end{subfigure}%
	\begin{subfigure}[b]{0.245\textwidth}
		\includegraphics[width=\textwidth]{%
      		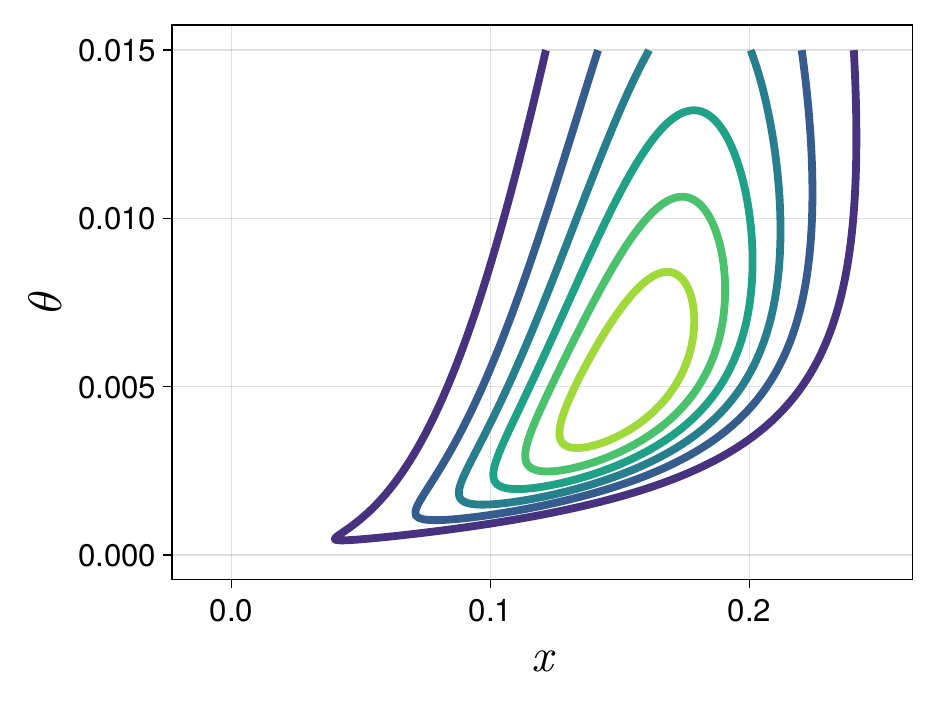} 
    		\caption{$r=-1/2$}
    		\label{fig:contour_posterior_toy_rm05}
  	\end{subfigure}%
	\begin{subfigure}[b]{0.245\textwidth}
		\includegraphics[width=\textwidth]{%
      		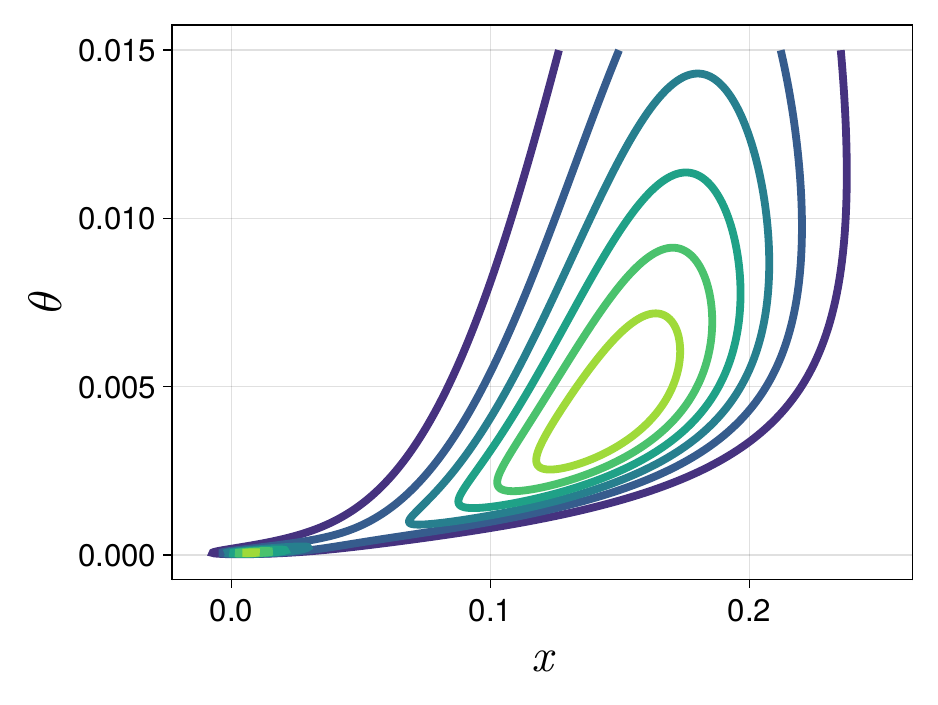} 
    		\caption{$r=-1$}
    		\label{fig:contour_posterior_toy_rm1}
  	\end{subfigure}%
  	\caption{ 
	Contour plots of the posterior density $\pi^y( x, \theta )$ in \cref{expl:toy_problem} with $y=0.2$ and $\sigma^2 = 10^{-2.8}$. 
	We use the parameter combinations $(r,\beta,\vartheta)$ as detailed in \cref{tab:parameters}. 
  	}
  	\label{fig:contour_posterior_toy}
\end{figure} 

\cref{fig:contour_posterior_toy} provides contour plots of the posterior density $\pi^y( x, \theta )$ in \cref{expl:toy_problem} for $y=0.2$ and $\sigma^2 = 10^{-2.8}$.
We use the parameter combinations $(r,\beta,\vartheta)$ as detailed in \cref{tab:parameters}.
Notably, \cref{fig:contour_posterior_toy_rm1} highlights two significant challenges for efficient MCMC sampling of SBL posteriors: the presence of high anisotropy, leading to strong correlations and poor mixing, and the existence of multiple modes that are difficult to traverse---one close to zero, resulting from the sparsity-promoting SBL prior, and a second at around $x \approx 0.15$, resulting from the observational data. 

\begin{remark}[Mixture models]
	We briefly address the connection between the hierarchical prior model \cref{eq:prior_model,eq:SBL_priors} and mixture formulations. 
	 If one is interested in inferring only $\mathbf{x}$---and not also the hyper-parameters $\boldsymbol{\theta}$---one can rewrite the joint posterior as $\pi^{y}( \mathbf{x}, \boldsymbol{\theta} ) = \pi( \mathbf{x} | \boldsymbol{\theta}, \mathbf{y} ) \, \pi( \boldsymbol{\theta} | \mathbf{y} )$ and marginalize out $\boldsymbol{\theta}$ to get the marginal posterior 
	 \begin{equation}\label{eq:marginal_posterior}
	 	\pi^{y}( \mathbf{x} ) 
			= \int \pi( \mathbf{x} | \boldsymbol{\theta}, \mathbf{y} ) \, \pi( \boldsymbol{\theta} | \mathbf{y} ) \intd \boldsymbol{\theta}.
	 \end{equation}
	 Notably, \cref{eq:marginal_posterior} is a mixture formulation of the marginal posterior $\pi^{y}( \mathbf{x} )$ with component density $\pi( \mathbf{x} | \boldsymbol{\theta}, \mathbf{y} ) \propto \pi( \mathbf{y} | \mathbf{x} ) \, \pi( \mathbf{x} | \boldsymbol{\theta} )$ and mixing density $\pi( \boldsymbol{\theta} | \mathbf{y} ) \propto \pi( \boldsymbol{\theta} ) \int \pi( \mathbf{y} | \mathbf{x}' ) \, \pi( \mathbf{x}' | \boldsymbol{\theta} ) \intd \mathbf{x}'$---assuming that $\mathbf{y}$ and $\boldsymbol{\theta}$ are conditional independent given $\mathbf{x}$.
	In principle, one can sample from \cref{eq:marginal_posterior} using, for instance, a marginal-then-conditional (MTC) strategy \cite{fox2016fast}: first draw $\boldsymbol{\theta}$-samples from the mixing density $\pi( \boldsymbol{\theta} | \mathbf{y} )$, then draw $\mathbf{x}$ from the component density $\pi( \mathbf{x} | \boldsymbol{\theta}, \mathbf{y} )$ conditioned on the $\boldsymbol{\theta}$-samples. 
	While the conditional density $\pi(\mathbf{x}\mid\boldsymbol{\theta},\mathbf{y})$ is Gaussian for linear data models with additive Gaussian noise, it can be difficult to explore for the more general likelihoods considered in this work.
	More importantly, sampling from the mixing density $\pi(\boldsymbol{\theta}\mid \mathbf{y})$ is itself highly challenging---even in the linear–Gaussian case.
	For example, \cite{flock2024continuous} addressed this difficulty by resorting to dimension-reduced approximations of $\pi(\boldsymbol{\theta}\mid \mathbf{y})$, as the exact distribution was too costly to explore directly.
	In contrast, the hierarchical representation \cref{eq:prior_model,eq:SBL_priors} retains the joint structure of $(\mathbf{x},\boldsymbol{\theta})$, making the dependence between these variables explicit.
	This structure is precisely what enables the proposed prior-normalization approach, which operates on the full joint prior.
\end{remark}
\section{Prior-normalizing TMs} 
\label{sec:transport} 

We use TMs \cite{villani2009optimal,marzouk2016sampling,parno2018transport} to transform the challenging SBL posterior \cref{eq:posterior} into a simpler prior-normalized posterior \cref{eq:posterior_reference} with a standard normal prior.

\subsection{The basic idea}
\label{sub:idea}

The performance of MCMC algorithms is primarily affected by the posterior landscape, which the prior can largely control.  
We thus employ prior-normalizing TMs to transform the SBL prior $\pi^0(\mathbf{x},\boldsymbol{\theta})$ into a standard normal density $\phi^0(\mathbf{u},\boldsymbol{\tau}) \propto \exp( -[ \| \mathbf{u} \|_2^2 + \| \boldsymbol{\tau} \|_2^2 ] / 2 )$. 
Specifically, we seek a diffeomorphism $S: (\mathbf{x},\boldsymbol{\theta}) \mapsto ( \mathbf{u}, \boldsymbol{\tau} )$ such that $( \mathbf{u}, \boldsymbol{\tau} ) = S(\mathbf{x},\boldsymbol{\theta})$ follows $\phi^0$ whenever $(\mathbf{x},\boldsymbol{\theta})$ follows $\pi^0$.  
In terms of densities, this means that the \emph{pushforward} of $\pi^0$ under $S$, $(S_{\sharp} \pi^0)( \mathbf{u}, \boldsymbol{\tau} ) := (\pi^0 \circ S^{-1})( \mathbf{u}, \boldsymbol{\tau} ) | \det \nabla S^{-1}( \mathbf{u}, \boldsymbol{\tau} ) |$, satisfies 
\begin{equation}\label{eq:PN_TM} 
	\left( S_{\sharp} \pi^0 \right)( \mathbf{u}, \boldsymbol{\tau} ) = \phi^0( \mathbf{u}, \boldsymbol{\tau} ),
\end{equation}
where $\nabla S^{-1}$ denotes the Jacobian of $S^{-1}$.
Since $S$ is a bijection, we can push forward the challenging SBL posterior $\pi^y(\mathbf{x},\boldsymbol{\theta})$ in \cref{eq:posterior} to obtain a simpler prior-normalized posterior $\phi^y = S_{\sharp} \pi^y$ as in \cref{eq:posterior_reference}. 
One can generate independent samples from the SBL posterior by pulling back samples of the prior-normalized posterior through the TM.

\subsection{Deriving prior-normalizing TMs}
\label{sub:deriving}

We analytically derive a prior-normalizing TM satisfying \cref{eq:PN_TM} by extending the approach from \cite{cui2022prior} to SBL priors. 
While \cite{cui2022prior} explored prior normalization for dimension reduction, they did not consider its effect on the posterior's geometry nor its application to hierarchical prior models.
Observe that \cref{eq:prior_model} can be expressed in a product-like form as 
\begin{equation}\label{eq:prod_form}
	\pi^0( \mathbf{x}, \boldsymbol{\theta} ) 
		= \prod_{i=1}^n \pi_i^0( x_i, \theta_i ),
\end{equation}
where $\pi_i^{0}( x_i, \theta_i ) = \pi_i^{x|\theta}( x_i | \theta_i ) \, \pi_i^{\theta}( \theta_i )$ with univariate densities $\pi_i^{x|\theta}( x | \theta ) = \mathcal{N}( x | 0, \theta )$ and $\pi_i^{\theta}( \theta ) = \mathcal{GG}( \theta | r, \beta, \vartheta )$. 
We therefore restrict our search to block-diagonal TMs\footnote{The TM \cref{eq:diag_transform} is called ``bock-diagonal'' because its $i$th component depends only on $x_i$ and $\theta_i$, resulting in its Jacobian being a block diagonal matrix.} of the form
\begin{equation}\label{eq:diag_transform}
	S( \mathbf{x}, \boldsymbol{\theta} ) 
		=  \begin{bmatrix} s_1( x_1, \theta_1 ) \\ \vdots \\ s_n( x_n, \theta_n ) \end{bmatrix}
\end{equation}
with $s_i$ pushing forward $\pi_i^{0}( x_i, \theta_i )$ to the 2D standard normal density $\phi_i^{0}( u_i, \tau_i )$. 
Importantly, we have replaced the problem of finding a high-dimensional TM $S$ with the easier task of finding the decoupled 2D TMs $s_1,\dots,s_n$.
Moreover, if we make the common choice that all $\theta_1,\dots,\theta_n$ follow the same generalized gamma distribution, then the $\pi^0_i$'s in \cref{eq:prod_form} are the same. 
In this case, we only have to find a \emph{single} 2D TM.
For this reason and ease of notation, we will omit the subindex ``$i$." 

\begin{figure}[tb]
	\centering 
	\begin{subfigure}[b]{0.45\textwidth}
		\includegraphics[width=\textwidth]{%
      		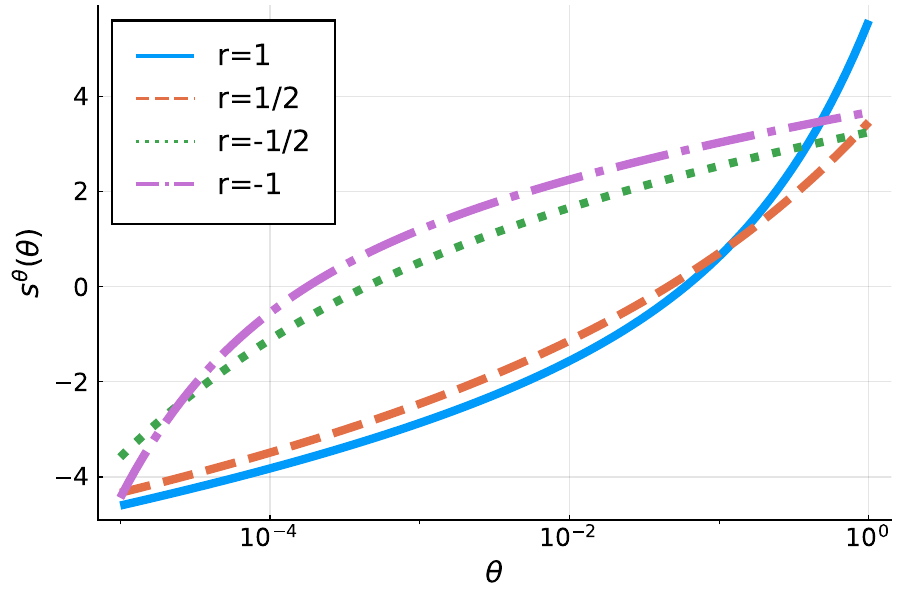} 
    		\caption{First component of the KR map $s$}
    		\label{fig:KR_map_s1}
  	\end{subfigure}%
	~
	\begin{subfigure}[b]{0.45\textwidth}
		\includegraphics[width=\textwidth]{%
      		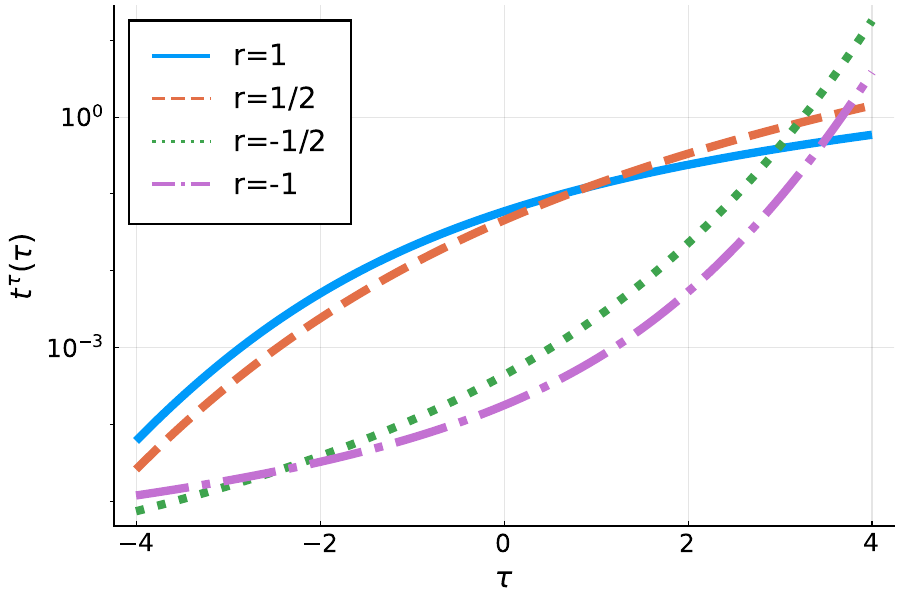} 
    		\caption{First component of the inverse KR map $t$}
    		\label{fig:KR_map_t1}
  	\end{subfigure}%
  	\caption{ 
  		First components of the KR map $s$ and its inverse $t = s^{-1}$ for the parameters in \cref{tab:parameters} 
  	}
  	\label{fig:KR_map}
\end{figure}

\begin{figure}[tb]
	\centering 
	\begin{subfigure}[b]{0.245\textwidth}
		\includegraphics[width=\textwidth]{%
      		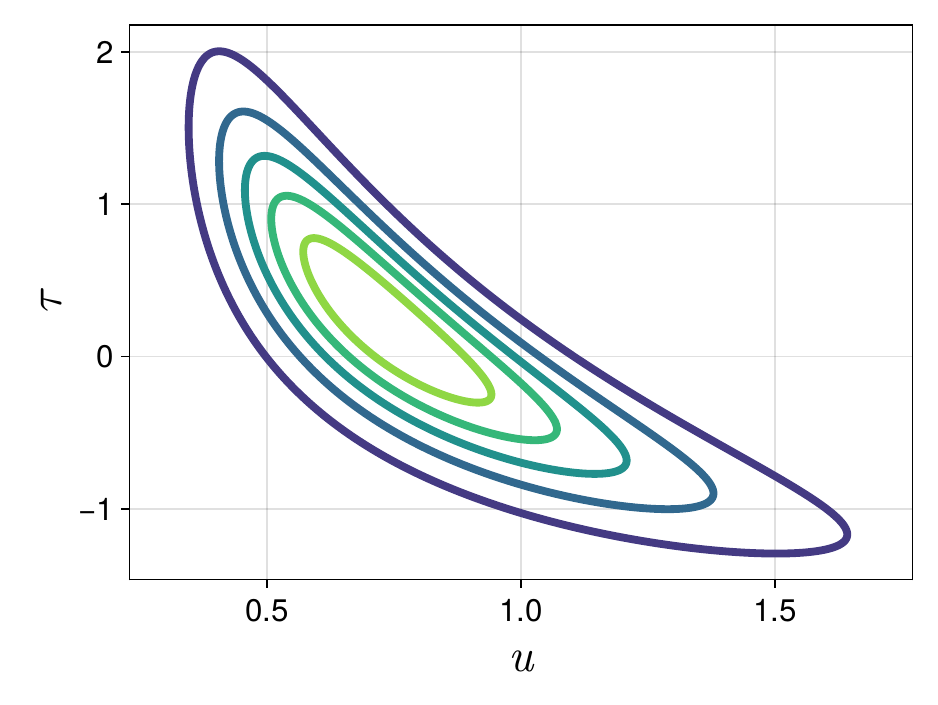} 
    		\caption{$r=1$}
    		\label{fig:contour_posterior_toy_transformed_rp1}
  	\end{subfigure}%
	\begin{subfigure}[b]{0.245\textwidth}
		\includegraphics[width=\textwidth]{%
      		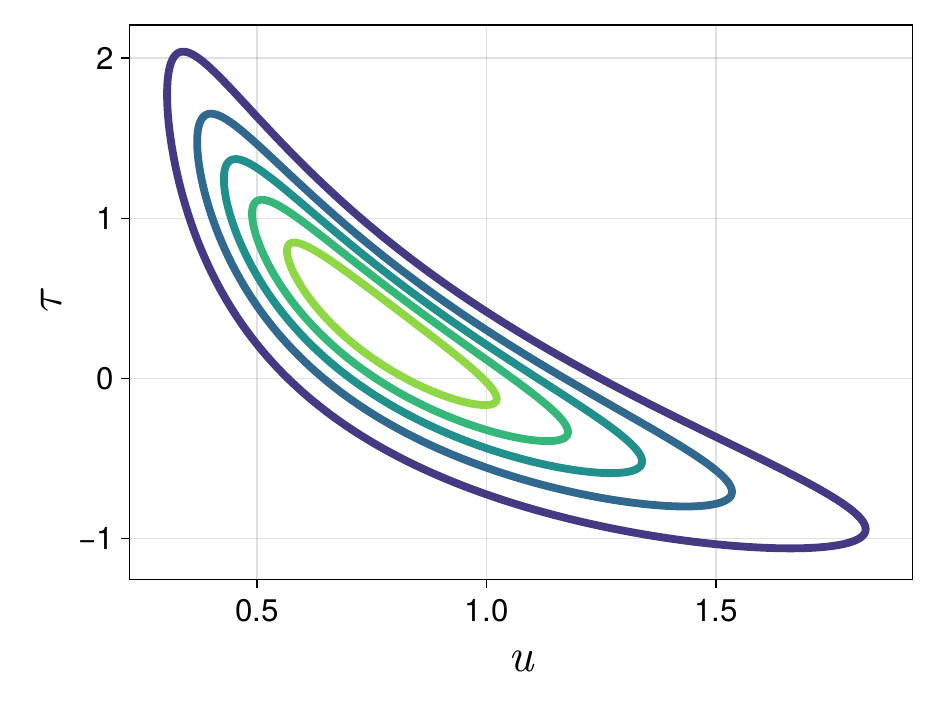} 
    		\caption{$r=1/2$}
    		\label{fig:contour_posterior_toy_transformed_rp05}
  	\end{subfigure}%
	\begin{subfigure}[b]{0.245\textwidth}
		\includegraphics[width=\textwidth]{%
      		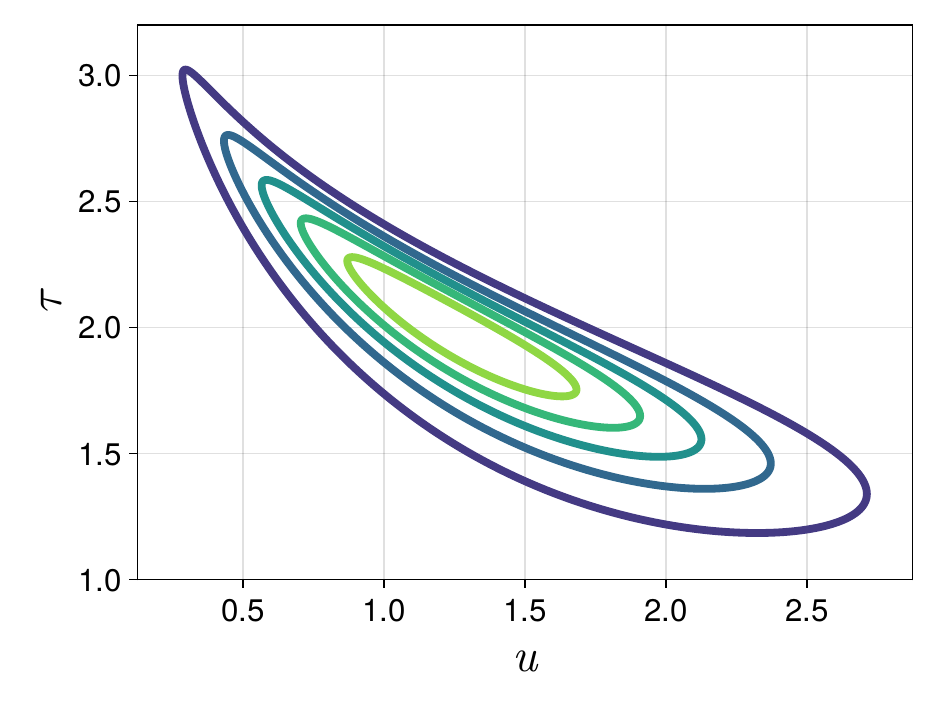} 
    		\caption{$r=-1/2$}
    		\label{fig:contour_posterior_toy_transformed_rm05}
  	\end{subfigure}%
	\begin{subfigure}[b]{0.245\textwidth}
		\includegraphics[width=\textwidth]{%
      		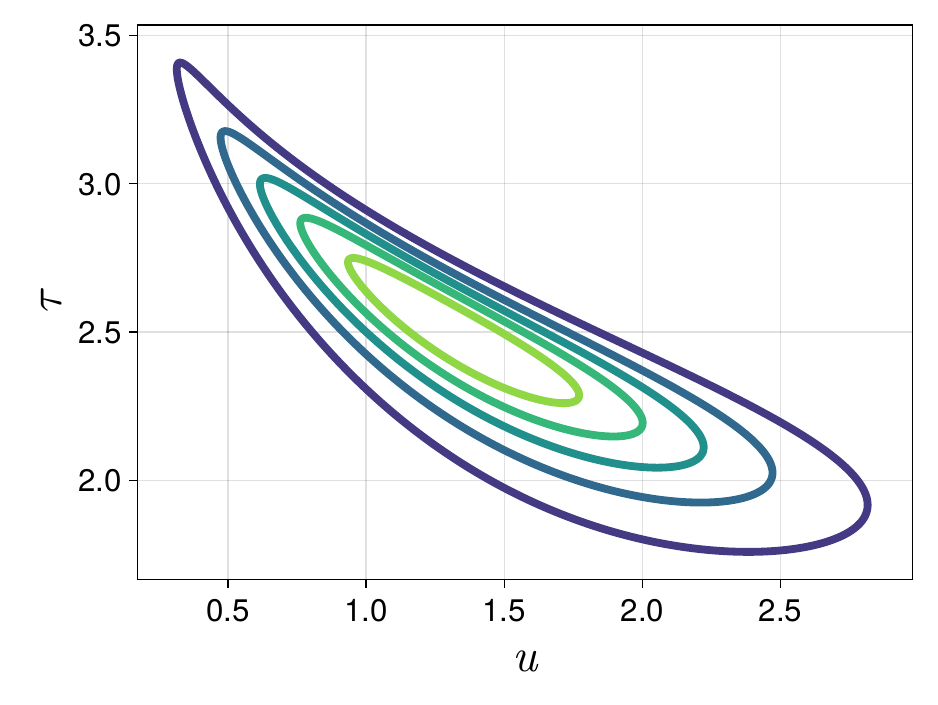} 
    		\caption{$r=-1$}
    		\label{fig:contour_posterior_toy_transformed_rm1}
  	\end{subfigure}%
  	\caption{
	Contour plots of the prior-normalized posterior for the univariate toy problem in \cref{expl:toy_problem}. 
	We use the same parameters as in \cref{fig:contour_posterior_toy} for the original posterior. 
	}
  	\label{fig:contour_posterior_toy_transformed}
\end{figure}

There can still be infinitely many TMs $s: \R \times \R_{>0} \to \R^2$ pushing forward $\pi^{0}( x, \theta ) = \pi^{x|\theta}( x | \theta ) \pi^{\theta}( \theta )$ to the standard normal density $\phi^{0}( u, \tau )$. 
One way of regularizing the problem and finding a unique map is to restrict the search to KR rearrangement \cite{rosenblatt1952remarks,knothe1957contributions,bogachev2005triangular}: 
\begin{equation}\label{eq:KR_TM} 
	s( x, \theta ) =  
	\begin{bmatrix} s^{\theta}( \theta ) \\ s^{x | \theta}( x; \theta ) \end{bmatrix}, 
\end{equation} 
where $s^{\theta}$ and $s^{x | \theta}$ push forward $\pi^{\theta}$ and $\pi^{x|\theta}$ to the univariate standard normal density, respectively \cite[Section 2.3]{santambrogio2015optimal}. 
Specifically, the second component of the KR map \cref{eq:KR_TM} is  
\begin{equation}\label{eq:s2} 
	s^{x|\theta}(x; \theta) = \frac{x}{\sqrt{\theta}},  
\end{equation}
as it transforms $\mathcal{N}(0|\theta)$ into $\mathcal{N}(0|1)$.
The first component is formally given by 
\begin{equation}\label{eq:s1} 
	s^{\theta} = \left( \Phi^0 \right)^{-1} \circ \mathcal{P}^{\theta},  
\end{equation}
where $\Phi^0(z) = \int_{-\infty}^z \phi^0(t) \intd t$ is the cumulative distribution function (CDF) of the univariate standard normal density $\phi^0$ and $\mathcal{P}^{\theta}(\theta) = \int_{0}^{\theta} \pi^0( t ) \intd t$ is the CDF of the generalized gamma hyper-prior $\pi^{\theta}$. 
Notably, \cref{eq:s1} corresponds to the univariate \emph{optimal TM} \cite{villani2009optimal}. 
\cref{fig:KR_map} illustrates the first components of $s^{\theta}$ and its inverse, $t^{\tau} = (s^{\theta})^{-1}$. 
Moreover, \cref{fig:contour_posterior_toy_transformed} shows the resulting prior-normalized posteriors for the toy problem in \cref{expl:toy_problem}. 

\begin{remark} 
	Alternative approaches for (approximating) prior-normalizing TMs include (multivariate) optimal TMs \cite{villani2009optimal,ambrosio2013user,santambrogio2015optimal} and normalizing flows \cite{rezende2015variational,baptista2023representation}. 
	However, unlike the KR map \cref{eq:KR_TM}, these TMs do not yield exact solutions and/or are not analytically tractable.
\end{remark}

\begin{remark}
	We note that the prior-normalized posterior $\phi^y = S_{\sharp} \pi^y$ in \cref{eq:posterior_reference} is, in general, neither guaranteed to be log-concave nor unimodal.
	For example, the scatter plots in \cref{fig:deblurring_scatter_AM_prior} show that the prior-normalized posterior can exhibit multi-modality---and hence a lack of log-concavity---even in the case of a linear signal deblurring problem.
	A systematic geometric characterization of the prior-normalized posterior remains an open question and is left for future work.
\end{remark}

\subsection{Computing the inverse prior-normalizing map}

Sampling from the prior-normalized posterior \cref{eq:posterior_reference} requires evaluating the transformed likelihood $g( \mathbf{u}, \boldsymbol{\tau}; \mathbf{y} ) = f( S^{-1}( \mathbf{u}, \boldsymbol{\tau} ); \mathbf{y} )$, which depends on computing the inverse prior-normalizing TM, $T = S^{-1}$. 
The map $T$ transforms samples from the prior-normalized posterior into samples from the original posterior. 
To derive $T$, recall that $S$ pushes forward the joint prior $\pi^0( \mathbf{x}, \boldsymbol{\theta} )$ to the standard normal density $\phi^0(\mathbf{u},\boldsymbol{\tau})$, i.e., $S_{\sharp} \pi^0 = \phi^0$.
The map $T$ has the reverse effect: 
It pushes forward $\phi^0$ to $\pi^0$, i.e., $T_{\#} \phi^0 = \pi^0$.
It also has a block diagonal form 
\begin{equation}\label{eq:diag_transform_T}
	T( \mathbf{u}, \boldsymbol{\tau} ) 
		=  \begin{bmatrix} t_1( u_1, \tau_1 ) \\ \vdots \\ t_n( u_n, \tau_n ) \end{bmatrix}
\end{equation}
with $t_i = s_i^{-1}$ pushing forward $\phi_i^{0}( u_i, \tau_i ) \propto \exp( -[ u_i^2 + \tau_i^2 ]/2 )$ to $\pi_i^{0}( x_i, \theta_i )$.  
Again omitting the subindex ``$i$," the inverse 2D KR map $t: \R^2 \to \R \times \R_{>0}$ is 
\begin{equation}\label{eq:KR_TM_inverse} 
	t( u, \tau ) =  
	\begin{bmatrix} 
		t^{\tau}( \tau ) \\ t^{u | \tau}( u; \tau ) 
	\end{bmatrix}, 
\end{equation} 
where $t^{\tau} = (s^{\theta})^{-1}$ and $t^{u | \theta} = (s^{x | \theta})^{-1}$ push forward the univariate standard normal density to $\pi^{\theta}$ and $\pi^{x|\theta}$, respectively. 
The first component is formally given by 
\begin{equation}\label{eq:t1} 
	t^{\tau} = \left( \mathcal{P}^{\theta} \right)^{-1} \circ \Phi^0,  
\end{equation}
where $\Phi^0(z)$ and $\mathcal{P}^{\theta}(\theta)$ are again the CDF of $\phi^0$ and $\pi^{\theta}$, respectively. 
We comment on the stable implementation of \cref{eq:t1} in \Cref{app:implementation}.
The second component is $t^{u | \tau}( u; \tau ) = \sqrt{t^{\tau}(\tau)} \, u$, as it transforms $\mathcal{N}(0|1)$ into $\mathcal{N}(0|\theta)$ with $\theta = t^{\tau}(\tau)$.

\subsection{Comparison with an existing transformation}

\begin{figure}[tb]
	\centering 
	\begin{subfigure}[b]{0.245\textwidth}
		\includegraphics[width=\textwidth]{%
      		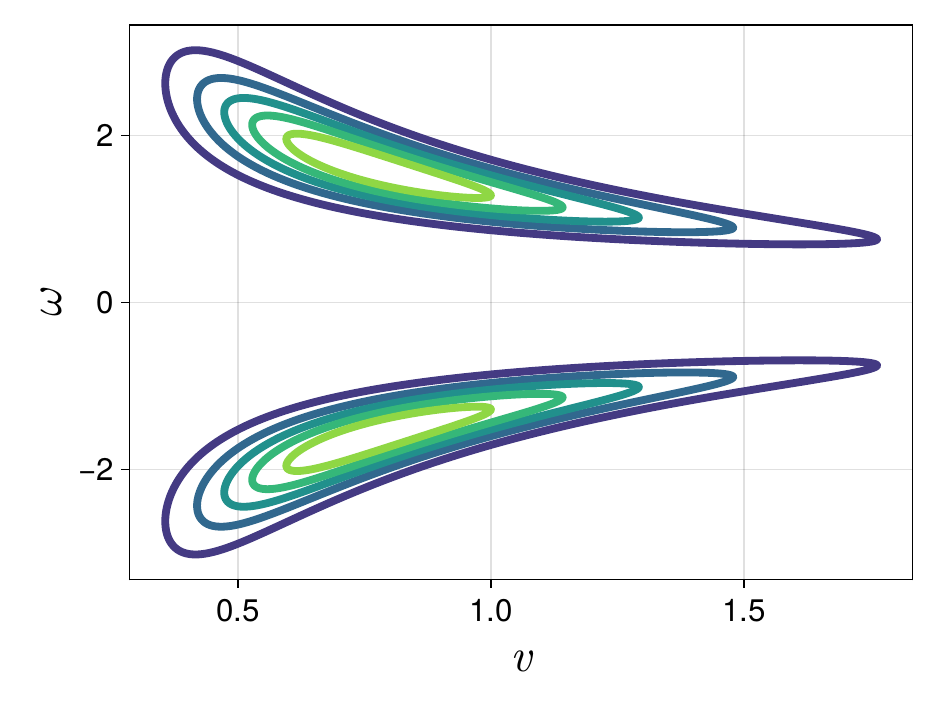} 
    		\caption{$r=1$}
    		\label{fig:contour_posterior_toy_Calvetti2024_rp1}
  	\end{subfigure}%
	\begin{subfigure}[b]{0.245\textwidth}
		\includegraphics[width=\textwidth]{%
      		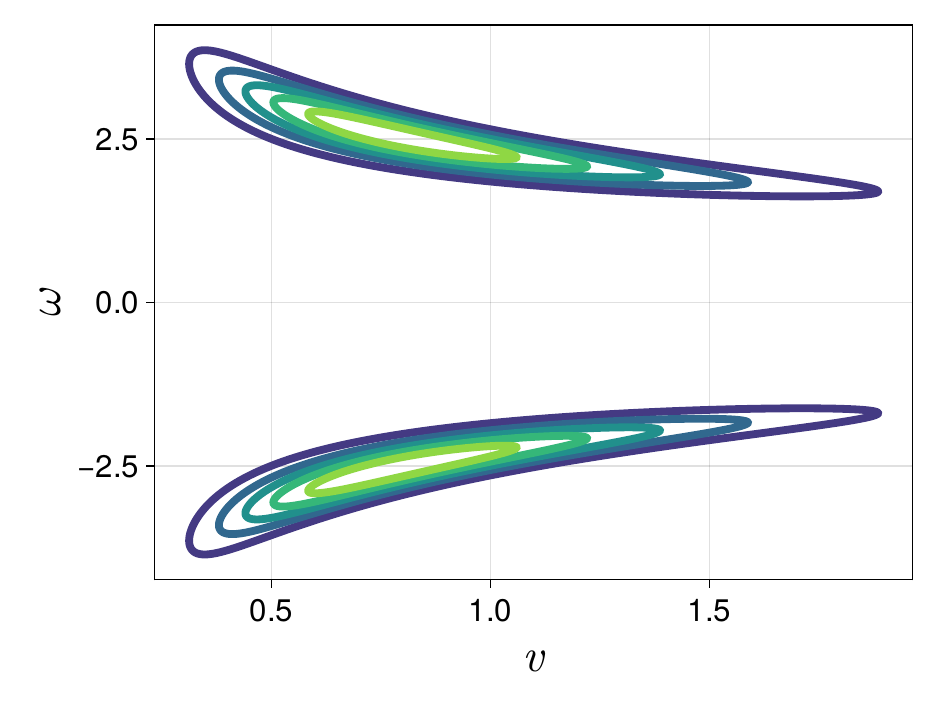} 
    		\caption{$r=1/2$}
    		\label{fig:contour_posterior_toy_Calvetti2024_rp05}
  	\end{subfigure}%
	\begin{subfigure}[b]{0.245\textwidth}
		\includegraphics[width=\textwidth]{%
      		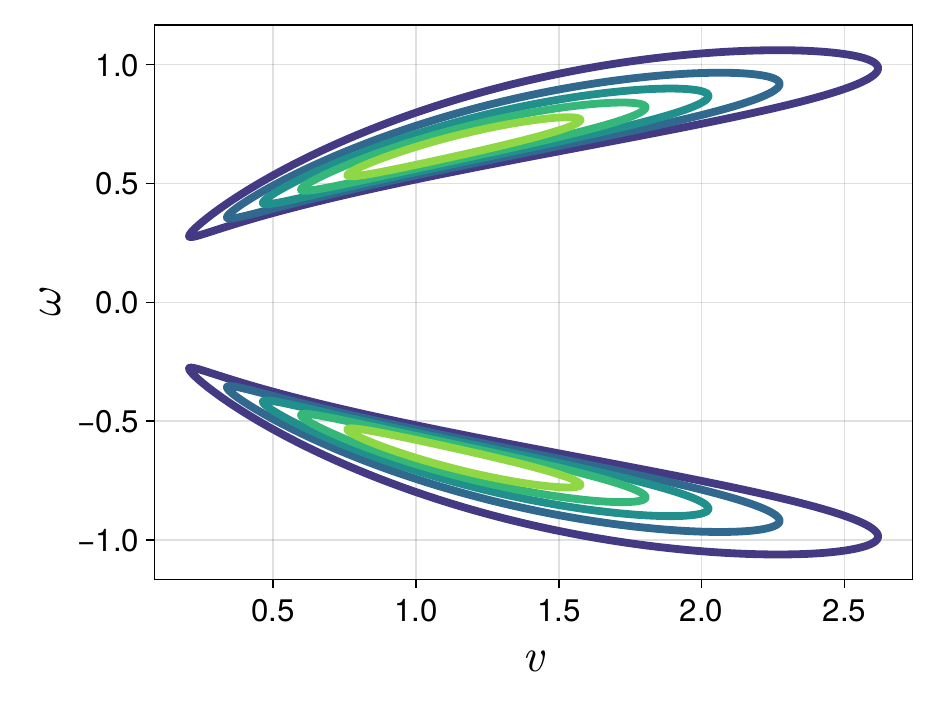} 
    		\caption{$r=-1/2$}
    		\label{fig:contour_posterior_toy_Calvetti2024_rm05}
  	\end{subfigure}%
	\begin{subfigure}[b]{0.245\textwidth}
		\includegraphics[width=\textwidth]{%
      		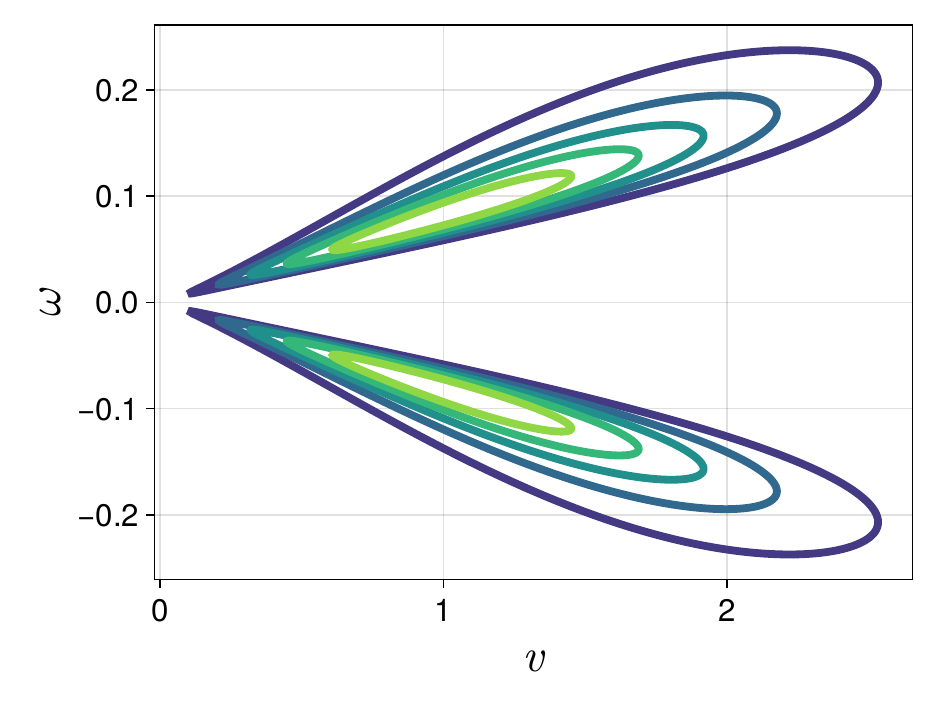} 
    		\caption{$r=-1$}
    		\label{fig:contour_posterior_toy_Calvetti2024_rm1}
  	\end{subfigure}%
  	\caption{ 
		Contour plots of the transformed posterior density $\psi^y( v, \omega )$ using the approach in \cite{calvetti2024computationally} for the univariate toy problem in \cref{expl:toy_problem} with $y=0.2$ and $\sigma^2 = 10^{-2.8}$. 
		We use the parameter combinations $(r,\beta,\vartheta)$ as detailed in \cref{tab:parameters}. 
  	}
  	\label{fig:contour_posterior_toy_Calvetti2024}
\end{figure}

We briefly summarize the transformation proposed in \cite{calvetti2024computationally} and compare it with the prior-normalization approach introduced here.
Consider the SBL posterior density in \cref{eq:posterior_expl}.
The transformation in \cite{calvetti2024computationally} is based on the reparametrization \cref{eq:posterior_expl} using $v_i^2 = x_i^2 / \theta_i$ and $\omega_i^2 = 2 ( \theta_i / \vartheta )^r$ for $i=1,\dots,n$.
This re-parametrization transforms the original SBL posterior into one with density 
\begin{equation}\label{eq:Calvetti2024_density}
	\psi^y( \mathbf{v}, \boldsymbol{\omega} ) 
		\propto e^{- p(\mathbf{v}, \boldsymbol{\omega})} \mathcal{N}( \mathbf{v}, \boldsymbol{\omega} | \mathbf{0}, I ),
\end{equation}
where $\mathcal{N}( \mathbf{v}, \boldsymbol{\omega} | \mathbf{0}, I ) \propto \exp( -[ \| \mathbf{v} \|_2^2 + \| \boldsymbol{\omega} \|_2^2 ] / 2 )$ denotes the density of the $2n$-dimensional standard normal distribution and $p(\mathbf{v}, \boldsymbol{\omega})$ is the potential defined as 
\begin{equation}\label{eq:Calvetti2024_potential}
	p(\mathbf{v}, \boldsymbol{\omega}) 
		= \frac{1}{2} \norm{ \Gamma_{\rm noise}^{-1/2} \left[ F\left( 2^{-1/2r} \vartheta^{1/2} \mathbf{v} | \boldsymbol{\omega} |^{1/r} \right) - \mathbf{y} \right] }_2^2 
			- (r \beta - 3/2) \sum_{i=1}^n \log |\omega_i|.
\end{equation}
Here, $\mathbf{v} | \boldsymbol{\omega} |^{1/r}$ is understood componentwise as $( \mathbf{v} | \boldsymbol{\omega} |^{1/r} )_i = v_i | \omega_i |^{1/r}$ for $i=1,\dots,n$.
This transformation may be interpreted as a partial normalization of the SBL prior: it converts the exponential terms of the conditional Gaussian prior and the generalized gamma hyperprior into a standard normal density, but it does not absorb the corresponding power function factors.
\Cref{fig:contour_posterior_toy_Calvetti2024} shows the resulting transformed posteriors for the toy problem from \cref{expl:toy_problem}.
These posteriors exhibit mirror symmetry with respect to the line ${(v,\omega) \mid \omega = 0}$; consequently, inference may be restricted to either half-space.
Nevertheless, the transformed density in \cref{fig:contour_posterior_toy_Calvetti2024} displays more pronounced anisotropy than the fully prior-normalized posterior proposed here (see \cref{fig:contour_posterior_toy_transformed}), particularly in the relevant case $r=-1$.
Additional numerical comparisons between the transformation of \cite{calvetti2024computationally} and our prior-normalization approach are provided in \Cref{sub:tests_deconvolution}.

\section{Computational experiments} 
\label{sec:tests} 

We assess the impact of the proposed hierarchical prior normalization approach on the efficiency of standard MCMC samplers.
To this end, we compare the performance of different samplers when applied to both the original posterior and its corresponding prior-normalized counterpart. 
For a better comparison, we transform all results into the $(\mathbf{x},\boldsymbol{\theta})$-coordinates. 
Our computational experiments cover a range of problems, including the toy example in \cref{expl:toy_problem} (\Cref{sub:tests_toy}), an undersampled signal deconvolution problem (\Cref{sub:tests_deconvolution}), inference of the initial data for the nonlinear inviscid Burgers equation from solution measurements (\Cref{sub:tests_Burgers}), and the recovery of an impulse image from noisy DCT data (\Cref{sub:tests_impulse}). 
To evaluate the performance of standard MCMC algorithms on both the original and the prior-normalized posteriors, we employed four different schemes: the adaptive Metropolis (AM) algorithm \cite{haario2001adaptive,andrieu2008tutorial,atchade2010limit}, the Metropolis-adjusted Langevin algorithm (MALA) \cite{atchade2006adaptive,marshall2012adaptive}, the Gibbs sampler \cite{agapiou2014analysis,ascolani2024scalability}, and the elliptical slice sampler (ES) \cite{murray2010elliptical}. 
For simplicity, we set the target mean acceptance rates to 23.4\% for AM and 57.4\% for MALA. 
We emphasize that these values are not intended to be optimal for the computational experiments considered here. 
Classical results such as \cite{roberts1998optimal,beskos2013optimal,yang2020optimal} justify these acceptance rates only in high-dimensional diffusion limits, and do not claim optimality in small dimensions, where higher acceptance rates are generally recommended. 
Our choices should therefore be viewed purely as pragmatic, fixed tuning parameters rather than theoretically optimal ones.
The Julia code to reproduce our results is available at \url{https://github.com/jglaubitz/paper-2025-SBL-priorNormalization}.

\subsection{Toy problem}
\label{sub:tests_toy}

Consider the setup of \cref{expl:toy_problem} with parameters $y = 0.2$ and $\sigma^2 = 10^{-2.8}$. 
We generate $J = 4$ independent chains, each consisting of $10^6$ samples. 
All chains are initialized at the MAP estimates corresponding to the original and prior-normalized posteriors.
Furthermore, via several experiments, we found the following MALA step sizes, $\varepsilon_{\rm o}$ for the original and $\varepsilon_{\rm pn}$ for prior-normalized posterior, to approximately achieve a mean acceptance rate of approximately 57.4\%: 
\begin{itemize}
	\item 
	For $r=1$, we use $\varepsilon_{\rm o} = 3 \cdot 10^{-3}$ and $\varepsilon_{\rm pn} = 7 \cdot 10^{-2}$.

	\item 
	For $r=1/2$, we use $\varepsilon_{\rm o} = 2 \cdot 10^{-3}$ and $\varepsilon_{\rm pn} = 7 \cdot 10^{-2}$.

	\item 
	For $r=-1/2$, we use $\varepsilon_{\rm o} = 2 \cdot 10^{-4}$ and $\varepsilon_{\rm pn} = 8 \cdot 10^{-2}$.

	\item 
	For $r=-1$, we use $\varepsilon_{\rm o} = 10^{-4}$ and $\varepsilon_{\rm pn} = 8 \cdot 10^{-2}$.
\end{itemize}
The associated values for the hyper-prior parameters $\beta$ and $\vartheta$ are listed in \cref{tab:parameters}. 
Notably, for the original posterior, the MALA step size $\varepsilon_{\rm o}$ must be reduced by an order of magnitude as $r$ decreases from $1$ to $-1$ to maintain the desired acceptance rate. 
In contrast, $\varepsilon_{\rm pn}$ for the prior-normalized posterior remains consistent across all values of $r$. 
This behavior can be explained by the posterior landscape of the prior-normalized posterior (see \cref{fig:contour_posterior_toy_transformed}) changing less compared to that of the original posterior (see \cref{fig:contour_posterior_toy}) across different values of $r$.

\begin{table}[h]
	\renewcommand{\arraystretch}{1.1}
	\centering
	\begin{adjustbox}{width=0.99\textwidth}
	\begin{tabular}{ l | c c c | c c c | c c c | c c c }
  		\multicolumn{13}{c}{AM} \\
		\toprule 
		 & \multicolumn{3}{c}{$r=1$} &
		 \multicolumn{3}{c}{$r=1/2$} &
		 \multicolumn{3}{c}{$r=-1/2$} &
		 \multicolumn{3}{c}{$r=-1$} \\
		posterior & 
		time & MPSRF & ESS & 
		time & MPSRF & ESS &  
		time & MPSRF & ESS &  
		time & MPSRF & ESS \\
		\midrule
    		original & 
		2.8 & 7.0e-5 & 1.3e+5 & 
		2.9 & 8.5e-5 & 1.0e+5 &  
		2.6 & 3.7e-4 & 4.3e+4 & 
		2.6 & 7.6e-4 & 3.0e+4 \\
   		prior-norm. & 
		3.3 & 1.8e-4 & 1.0e+5 & 
		3.6 & 6.0e-5 & 9.4e+4 &  
		3.3 & 1.3e-4 & 1.0e+5 & 
		3.3 & 2.3e-4 & 9.0e+4 \\
		\bottomrule
	\end{tabular}
	\end{adjustbox}
	\\ \vspace{.4cm}
	\begin{adjustbox}{width=0.99\textwidth}
	\begin{tabular}{ l | c c c | c c c | c c c | c c c }
  		\multicolumn{13}{c}{MALA} \\
		\toprule 
		 & \multicolumn{3}{c}{$r=1$} &
		 \multicolumn{3}{c}{$r=1/2$} &
		 \multicolumn{3}{c}{$r=-1/2$} &
		 \multicolumn{3}{c}{$r=-1$} \\
		posterior & 
		time & MPSRF & ESS & 
		time & MPSRF & ESS & 
		time & MPSRF & ESS & 
		time & MPSRF & ESS \\
		\midrule
    		original & 
		9.1 & 1.8e-3 & 5.9e+4 & 
		10 & 4.1e-3 & 2.1e+4 &  
		10 & 1.6e-2 & 3.3e+3 & 
		8.9 & 1.5e-1 & 2.3e+3 \\
   		prior-norm. & 
		11 & 2.8e-3 & 4.1e+4 & 
		12 & 3.5e-3 & 2.0e+4 &  
		12 & 1.1e-3 & 2.3e+4 & 
		13 & 4.8e-3 & 1.4e+4 \\
		\bottomrule
	\end{tabular}
	\end{adjustbox}
	\caption{
		Computational wall time in seconds (``time"), MPSRF minus one multiplied by the time (``MPSRF"), and the ESS per second (``ESS'') for for the original and prior-normalized posterior for \cref{expl:toy_problem}. 
		A smaller MPSRF and higher ESS indicate better sampler performance.
	}
  	\label{tab:toy_MPSRF}
\end{table}

\cref{tab:toy_MPSRF} measure the samplers' performance using the multivariate potential scale reduction factor (MPSRF) and the effective sample size (ESS).
The MPSRF is commonly used to assess the convergence of multiple MCMC chains \cite{brooks1998general}.
One typically has $\operatorname{MPSRF} \geq 1$ (assuming the chains have overdispersed starting points that cause the inter-chain variance to be larger than the within-chain variance). 
When the MPSRF approaches $1$, the variance within each sequence approaches the variance across sequences, thus indicating that each chain has converged to the target distribution. 
The literature contains several recommendations for the MPSRF values that indicate convergence. 
For example, \cite{gelman1992inference} suggests $\operatorname{MPSRF} < 1.1$ while \cite{vehtari2021rank} argues for the more conservative threshold $\operatorname{MPSRF} < 1.01$.
The ESS of a Markov chain is defined as the number of independent samples that are needed to estimate $\mathbb{E}[ \mathcal{G} ]$ (for some quantity of interest $\mathcal{G}(\mathbf{z})$) with the same statistical accuracy as an estimate from the Markov chain.
It is a measure of the amount of information contained in the MCMC chain and should be as large as possible; see \cite{gelman2003bayesian,wolff2004monte}. 
We return to \cref{tab:toy_MPSRF}, which presents the computational wall time in seconds (``time"), MPSRF minus one multiplied by the time (``MPSRF"), and the ESS for the AM and MALA samplers applied to the original and prior-normalized posterior of the in \cref{expl:toy_problem}.
For $r=1$, corresponding to a log-concave original posterior, the MPSRF and ESS values indicate better performance for the original posterior. 
However, as $r$ decreases from $1$ to $-1$, the MPSRF and ESS values demonstrate superior performance for the prior-normalized posterior. 
We also observe a faster decay of autocorrelations and improved mixing in the trace plots for the prior-normalized posterior compared to the original.
These plots are not included here but are available in the code repository.

\subsection{Signal deconvolution}
\label{sub:tests_deconvolution}

\begin{figure}[tb]
	\centering 
	\includegraphics[width=0.6\textwidth]{%
      		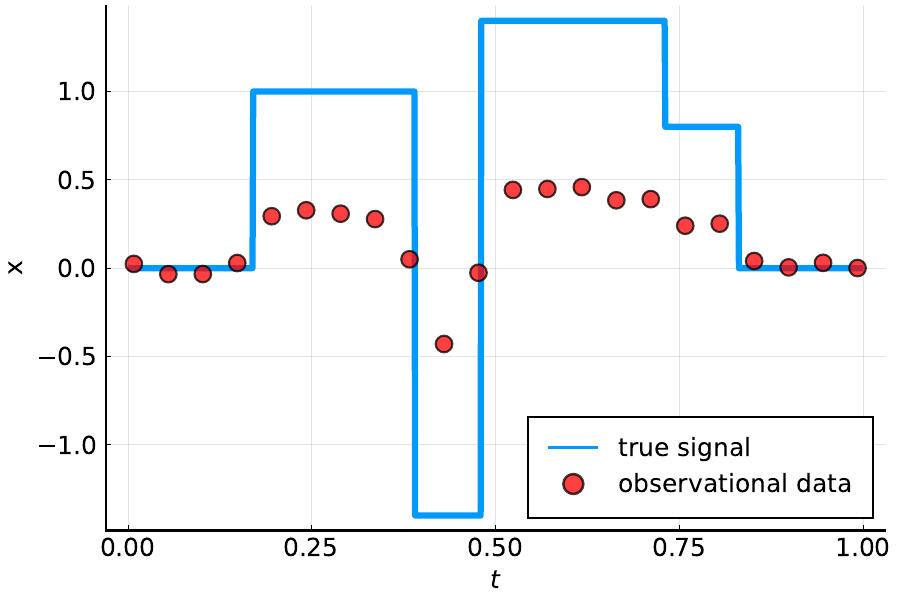} 
    	\caption{ 
		True signal and noisy, blurred observational data 
	}
  	\label{fig:deblurring_data}
\end{figure}

We consider a signal deblurring problem from \cite{calvetti2024computationally}. 
The goal is to estimate the nodal values of a piecewise constant signal $x: [0,1] \to \mathbb{R}$ from noisy observations of its convolution with a Gaussian kernel:
\begin{equation}\label{eq:signal_convolution}
	y_j = \int_0^1 k( t_j - s ) x(s) \, \mathrm{d}s + e_j, \quad 
	k(t) = A \exp\left( -\frac{t^2}{2 \omega^2} \right), \quad 
	j = 1, \dots, m, 
\end{equation}
where the kernel has amplitude $A = 6.2$ and width $\omega = 2 \cdot 10^{-2}$. 
The integral is discretized using $n = 128$ equidistant points. 
The observation points coincide with every sixth discretization node, yielding $m=22$ observations. 
The above setup leads to the data model
\begin{equation}\label{eq:signal_model}
	\mathbf{y} = F \mathbf{x} + \mathbf{e},
\end{equation}
where $F \in \mathbb{R}^{m \times n}$ is the matrix representation of the discretized convolution in \cref{eq:signal_convolution}, $\mathbf{x} = [x_1, \dots, x_n]$ contains the nodal values of the signal, and $e_j \sim \mathcal{N}(0, \sigma^2)$ with $\sigma = 3 \cdot 10^{-2}$. 
To avoid the ``inverse crime" \cite{kaipio2007statistical}, where the same model is used for both data generation and the inverse problem, we generate the observational data using a finer discretization with $10^3$ equidistant grid points. 
\Cref{fig:deblurring_data} illustrates the piecewise constant signal $x$ and the noisy, blurred, and undersampled data $\mathbf{y}$.

To formulate a sparse representation of $\mathbf{x}$, we introduce the finite difference matrix 
\begin{equation}\label{eq:FD_matrix}
	L = \begin{bmatrix}
		1 & & & \\ 
		-1 & 1 & & \\ 
		& \ddots & \ddots & \\ 
		& & -1 & 1
	\end{bmatrix} 
	\in \mathbb{R}^{n \times n} 
\end{equation}
The first row models a homogeneous Dirichlet boundary condition (the signal is zero at the left domain boundary), making $L$ invertible. 
We can then express the discrete signal $\mathbf{x}$ in terms of its increments $\mathbf{z}$ as $\mathbf{x} = L^{-1} \mathbf{z}$. 
This allows us to reformulate \cref{eq:signal_model} as 
\begin{equation}\label{eq:increments_model}
	\mathbf{y} = F L^{-1} \mathbf{z} + \mathbf{e}.   
\end{equation} 
Consequently, we model $\mathbf{x}$ being piecewise constant by assuming that its increments $\mathbf{z}$ are sparse, which is promoted by employing the sparsity-promoting SBL prior.
For brevity, we focus on assessing the performance of the AM algorithm  \cite{haario2001adaptive,andrieu2008tutorial,atchade2010limit} with a target mean acceptance rate of 23.4\% for the original posterior, the transformed posterior from \cite{calvetti2024computationally}, and our prior-normalized posteriors with $r=\pm 1$ and $\beta, \vartheta$ as in \cref{tab:parameters}.
The results for $r=\pm 1/2$ are qualitatively comparable to those for $r=\pm 1$, respectively, and are thus omitted here.
We generate $J = 6$ independent chains, each comprising $10^7$ samples. 
To manage memory efficiently, we apply thinning, retaining only every $10^3$th sample from each chain.

\subsubsection{Initializing with the MAP estimate}
\label{subsub:tests_deconvolution_MAP}

\begin{figure}[tb]
	\centering 
	\begin{subfigure}[b]{0.32\textwidth}
		\includegraphics[width=\textwidth]{%
      		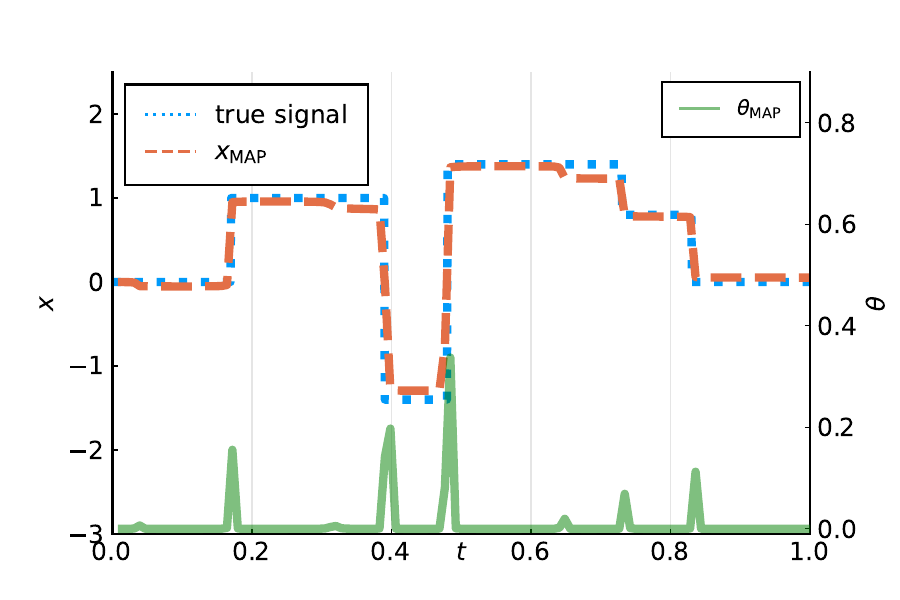} 
    		\caption{Original, $r=1$}
    		\label{fig:deblurring_model1_MAP_original}
  	\end{subfigure}%
	~
	\begin{subfigure}[b]{0.32\textwidth}
		\includegraphics[width=\textwidth]{%
      		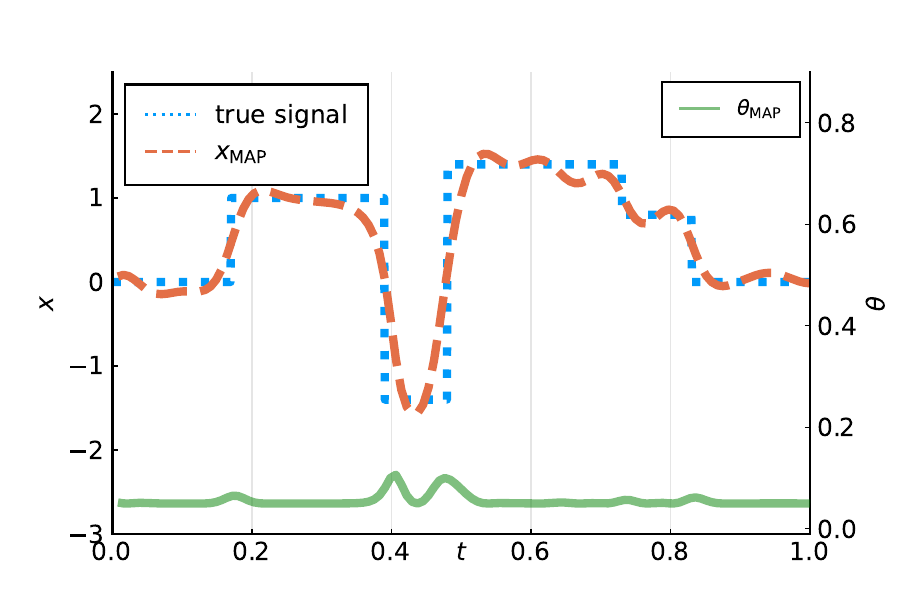} 
    		\caption{Approach from \cite{calvetti2024computationally}, $r=1$}
    		\label{fig:deblurring_model1_MAP_Calvetti2024}
  	\end{subfigure}%
	~
	\begin{subfigure}[b]{0.32\textwidth}
		\includegraphics[width=\textwidth]{%
      		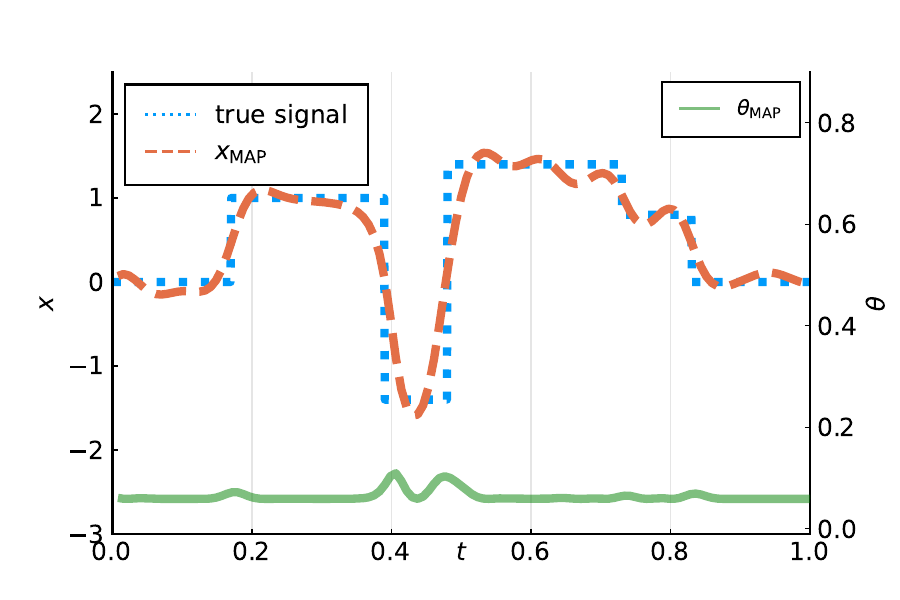} 
    		\caption{Prior-normalized, $r=1$}
    		\label{fig:deblurring_model1_MAP_priorNormalized}
  	\end{subfigure}%
	\\
	\begin{subfigure}[b]{0.32\textwidth}
		\includegraphics[width=\textwidth]{%
      		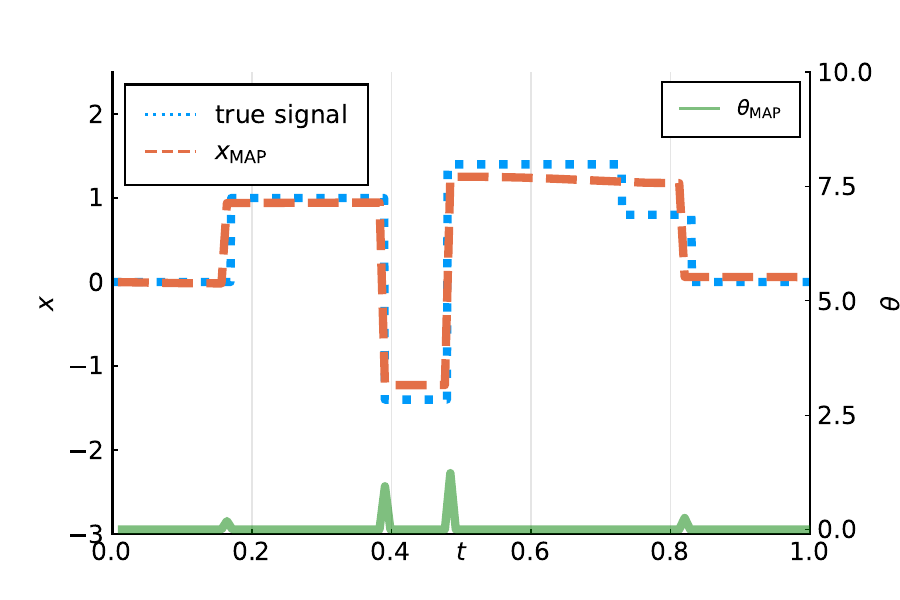} 
    		\caption{Original, $r=-1$}
    		\label{fig:deblurring_model4_MAP_original}
  	\end{subfigure}%
	~
	\begin{subfigure}[b]{0.32\textwidth}
		\includegraphics[width=\textwidth]{%
      		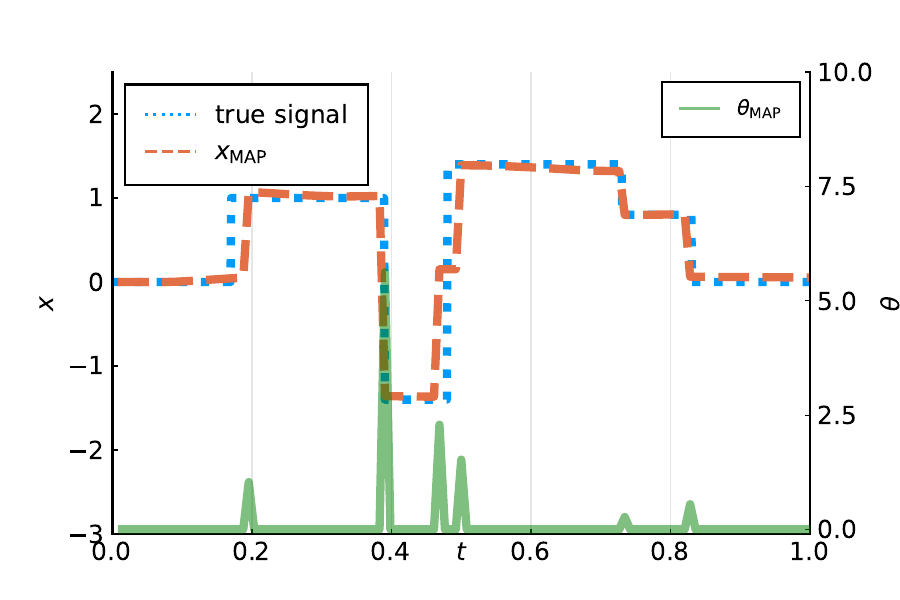} 
    		\caption{Approach from \cite{calvetti2024computationally}, $r=-1$}
    		\label{fig:deblurring_model4_MAP_Calvetti2024}
  	\end{subfigure}%
	~
	\begin{subfigure}[b]{0.32\textwidth}
		\includegraphics[width=\textwidth]{%
      		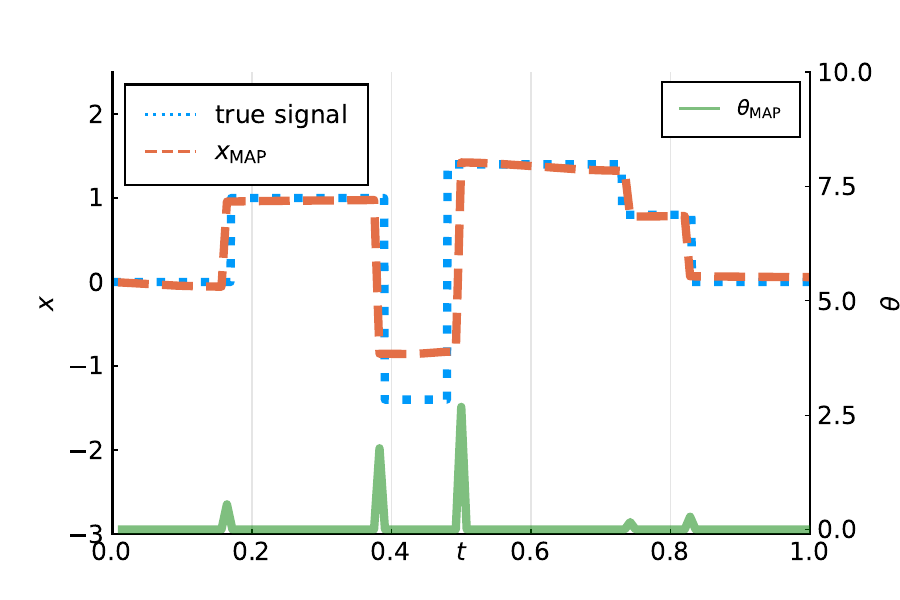} 
    		\caption{Prior-normalized, $r=-1$}
    		\label{fig:deblurring_model4_MAP_priorNormalized}
  	\end{subfigure}%
  	\caption{
		MAP estimates $(\mathbf{x}^{\rm MAP},\boldsymbol{\theta}^{\rm MAP})$ of the original posterior (left), pullback of the MAP estimates $(\mathbf{v}^{\rm MAP},\boldsymbol{\omega}^{\rm MAP})$ of the transformed posterior from \cite{calvetti2024computationally} (middle), and pullback of the MAP estimates $(\mathbf{u}^{\rm MAP},\boldsymbol{\tau}^{\rm MAP})$ of our prior-normalized posterior (right).
	}
  	\label{fig:deblurring_MAP}
\end{figure}

We first initiate all chains using the respective MAP estimates, illustrated in \cref{fig:deblurring_MAP}. 
We compute these using the limited-memory Broyden--Fletcher--Goldfarb--Shanno (L-BFGS) algorithm \cite{nocedal2006numerical} to the negative log-posterior. 
Notably, we generally do not expect the pullback of $(\mathbf{u}^{\rm MAP}, \boldsymbol{\tau}^{\rm MAP} )$ to correspond to $(\mathbf{z}^{\rm MAP}, \boldsymbol{\theta}^{\rm MAP} )$. 

\begin{figure}[tb]
	\centering 
	\begin{subfigure}[b]{0.32\textwidth}
		\includegraphics[width=\textwidth]{%
      		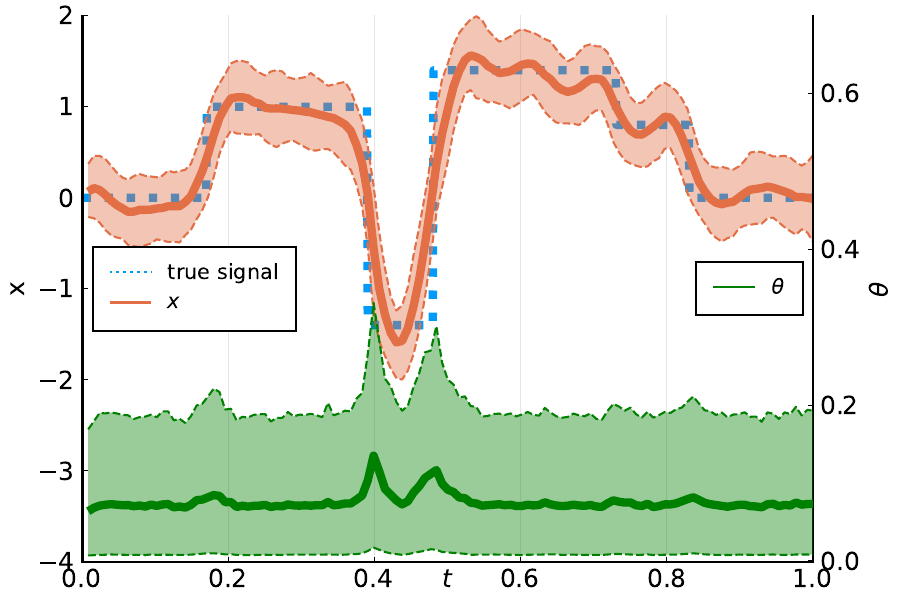} 
    		\caption{Original, $r=1$}
    		\label{fig:deblurring_model1_UQquantile_original_AM_MAP}
  	\end{subfigure}%
	~
	\begin{subfigure}[b]{0.32\textwidth}
		\includegraphics[width=\textwidth]{%
      		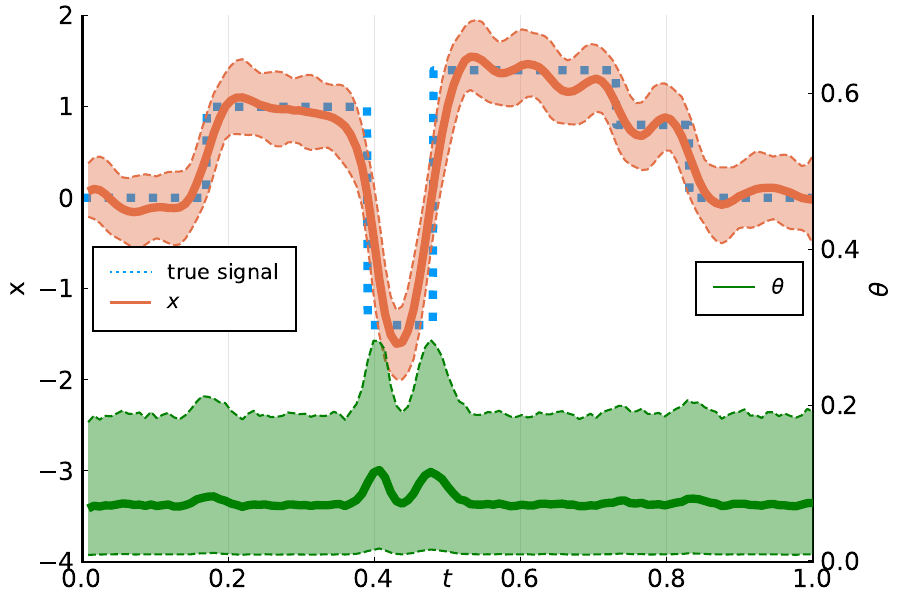} 
    		\caption{Approach from \cite{calvetti2024computationally}, $r=1$}
    		\label{fig:deblurring_model1_UQquantile_Calvetti2024_AM_MAP}
  	\end{subfigure}%
	~
	\begin{subfigure}[b]{0.32\textwidth}
		\includegraphics[width=\textwidth]{%
      		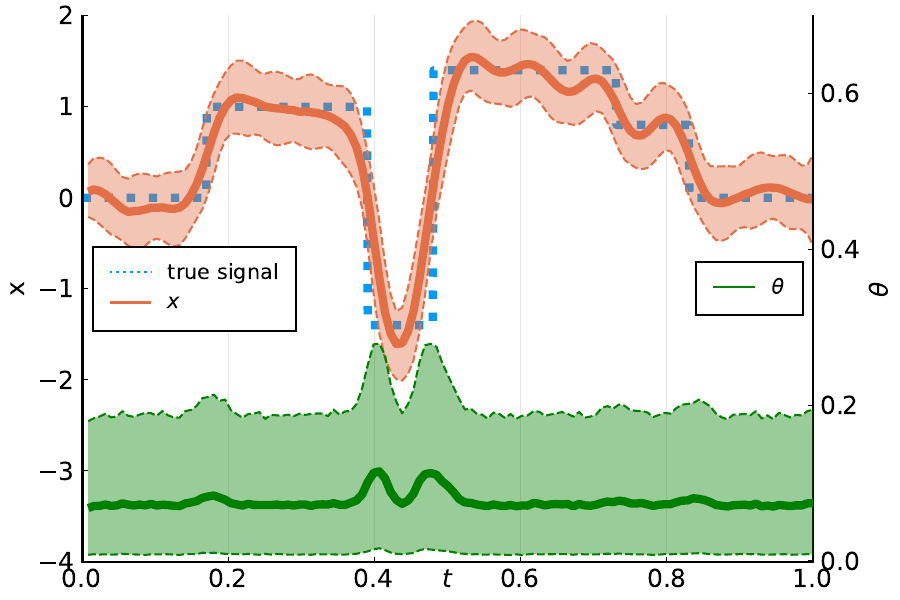} 
    		\caption{Prior-normalized (our), $r=1$}
    		\label{fig:deblurring_model1_UQquantile_priorNormalized_AM_MAP}
  	\end{subfigure}%
	\\
	\begin{subfigure}[b]{0.32\textwidth}
		\includegraphics[width=\textwidth]{%
      		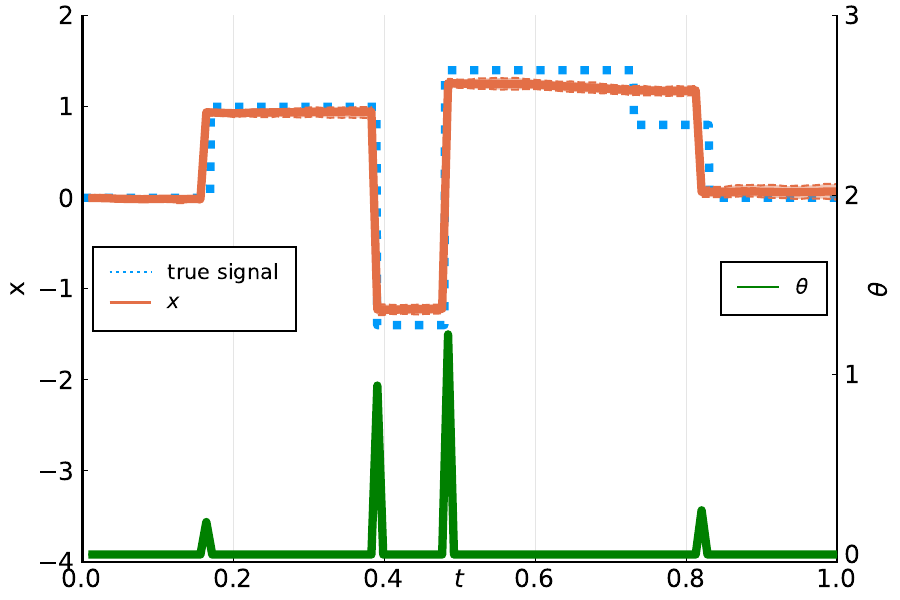} 
    		\caption{Original, $r=-1$}
    		\label{fig:deblurring_model4_UQquantile_original_AM_MAP}
  	\end{subfigure}%
	~
	\begin{subfigure}[b]{0.32\textwidth}
		\includegraphics[width=\textwidth]{%
      		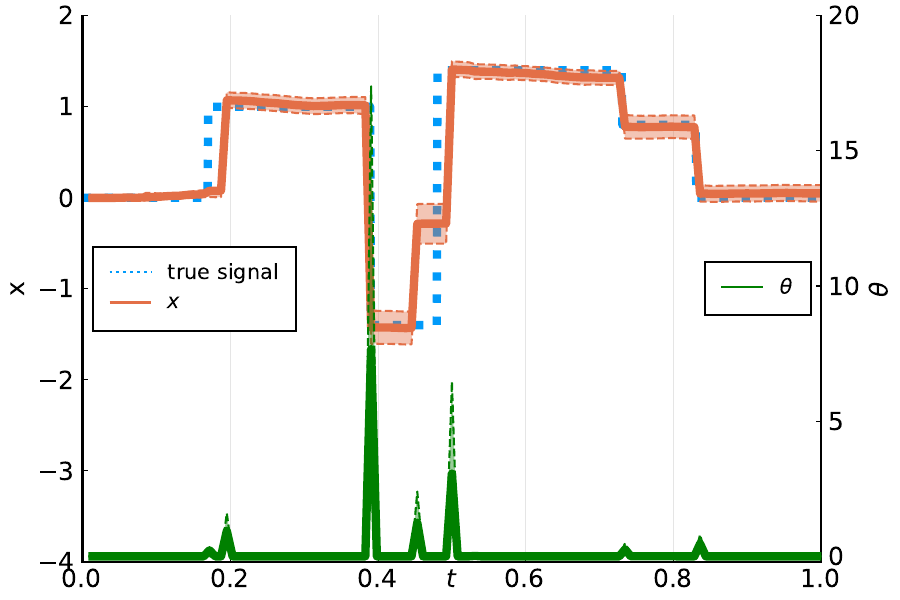} 
    		\caption{Approach from \cite{calvetti2024computationally}, $r=-1$}
    		\label{fig:deblurring_model4_UQquantile_Calvetti2024_AM_MAP}
  	\end{subfigure}%
	~
	\begin{subfigure}[b]{0.32\textwidth}
		\includegraphics[width=\textwidth]{%
      		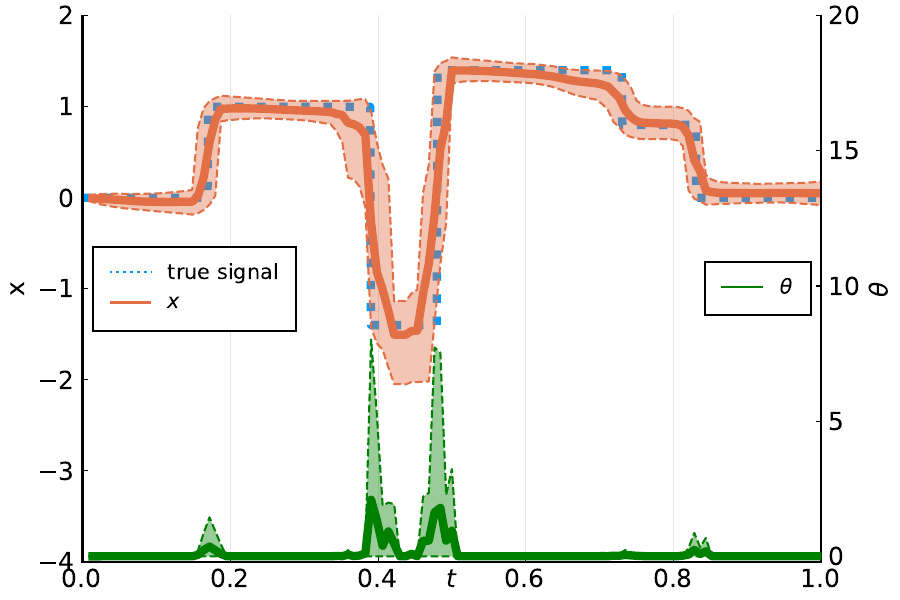} 
    		\caption{Prior-normalized (our), $r=-1$}
    		\label{fig:deblurring_model4_UQquantile_priorNormalized_AM_MAP}
  	\end{subfigure}%
	\caption{  
		Mean and $90\%$ quantile ranges of the $\mathbf{x}$- and $\boldsymbol{\theta}$-samples of the original posterior (left), the transformed posterior from \cite{calvetti2024computationally} (middle), and the proposed prior-normalized posterior (right) for $r = \pm 1$. 
		We used the AM algorithm initialized with the respective MAP estimate. 
	}
  	\label{fig:deblurring_UQ_AM_MAP}
\end{figure}

Next, \cref{fig:deblurring_UQ_AM_MAP}  presents the mean and $90\%$ quantile range for the $\mathbf{x}$- and $\boldsymbol{\theta}$-samples of the original posteriors, the transformed posteriors from \cite{calvetti2024computationally}, and our prior-normalized posteriors.  
The dashed lines correspond to the $5\%$ and $95\%$ quantiles. 
For $r=1$, corresponding to a log-concave original posterior, the AM method performs similarly on all three posteriors. 
In contrast, noticeable differences emerge for $r = -1$, where the sample means from the original posterior in \cref{fig:deblurring_model4_UQquantile_original_AM_MAP} and the transformed posterior from \cite{calvetti2024computationally} in \cref{fig:deblurring_model4_UQquantile_Calvetti2024_AM_MAP} closely match the MAP estimates shown in \cref{fig:deblurring_model4_MAP_original,fig:deblurring_model4_MAP_Calvetti2024}, which were used to initialize the MCMC chains.
Additionally, the $90\%$ quantile ranges in \cref{fig:deblurring_model4_UQquantile_original_AM_MAP,fig:deblurring_model4_UQquantile_Calvetti2024_AM_MAP} are fairly narrow, indicating an underestimation of uncertainty.
The above observations suggest that the AM method explores the original posterior and the transformed posterior from \cite{calvetti2024computationally} only locally, in a high-density region around the MAP estimate. 

\begin{figure}[tb]
	\centering 
	\begin{subfigure}[b]{0.32\textwidth}
		\includegraphics[width=\textwidth]{%
      		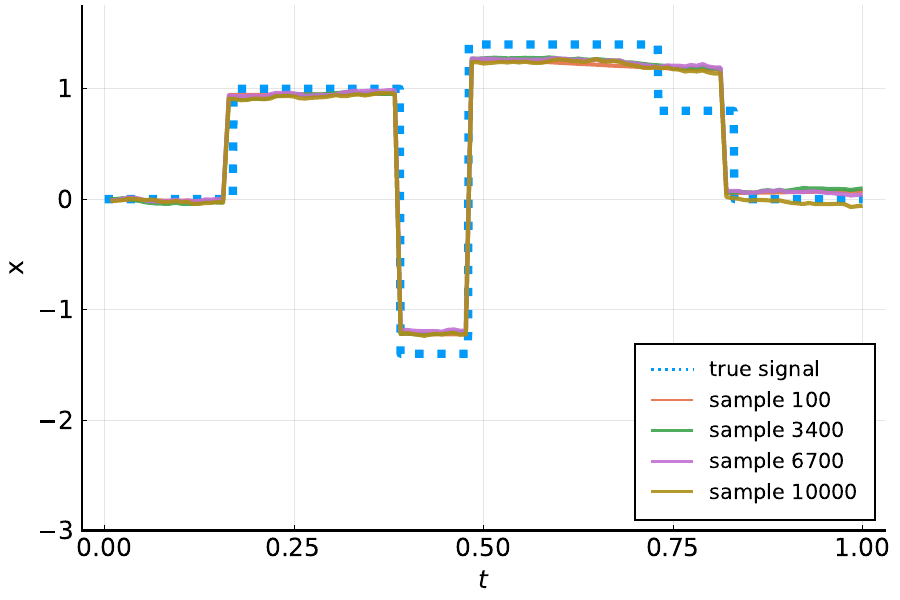} 
    		\caption{Original}
    		\label{fig:deblurring_model4_xSamples_original_AM_MAP}
  	\end{subfigure}
	~
	\begin{subfigure}[b]{0.32\textwidth}
		\includegraphics[width=\textwidth]{%
      		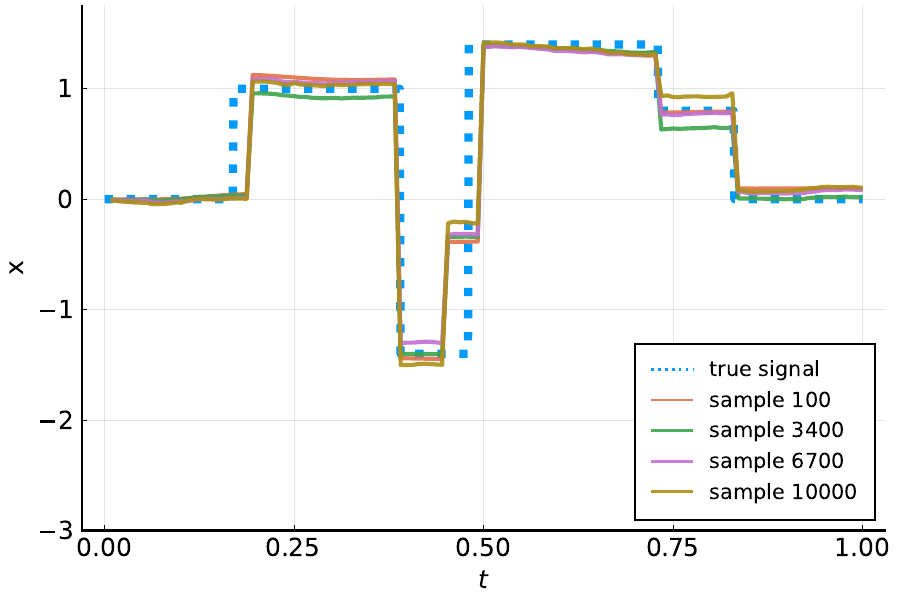} 
    		\caption{Approach from \cite{calvetti2024computationally}}
    		\label{fig:deblurring_model4_xSamples_Calvetti2024_AM_MAP}
  	\end{subfigure}%
	~
	\begin{subfigure}[b]{0.32\textwidth}
		\includegraphics[width=\textwidth]{%
      		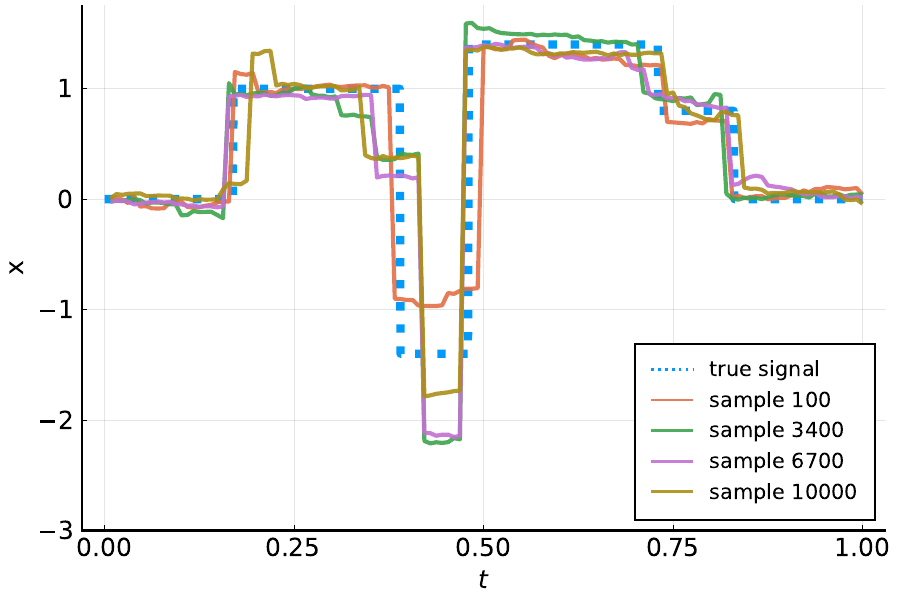} 
    		\caption{Prior-normalized}
    		\label{fig:deblurring_model4_xSamples_priorNormalized_AM_MAP}
  	\end{subfigure}%
  	\caption{ 
		Individual $x$-samples of the first MCMC chain from the original posterior, the transformed posterior from \cite{calvetti2024computationally}, and the prior-normalized posterior for $r=-1$. 
		We used the AM algorithm initialized with the respective MAP estimate.    
	}
  	\label{fig:deblurring_samples_AM_MAP}
\end{figure} 

Furthermore, \cref{fig:deblurring_samples_AM_MAP} illustrates four individual $x$-samples from the first MCMC chain for both the original posterior, the transformed posterior from \cite{calvetti2024computationally}, and our prior-normalized posterior for $r=-1$. 
We observe in \cref{fig:deblurring_model4_xSamples_original_AM_MAP} that each sample from the original posterior remains close to the MAP estimate, providing further evidence that the AM sampler explores the original posterior only locally around the MAP estimate. 
While the individual samples of the transformed posterior from \cite{calvetti2024computationally} show slightly more variability in \cref{fig:deblurring_model4_xSamples_Calvetti2024_AM_MAP}, they all have the same jump locations---again indicating that the AM sampler explores the transformed posterior from \cite{calvetti2024computationally} only locally around the MAP estimate.
In contrast, \cref{fig:deblurring_model4_xSamples_priorNormalized_AM_MAP} shows that the samples from the prior-normalized posterior differ visibly from one another, indicating that the same AM sampler explores the prior-normalized posterior more effectively.

\begin{figure}[tb]
	\centering 
	\begin{subfigure}[b]{0.32\textwidth}
		\includegraphics[width=\textwidth]{%
      		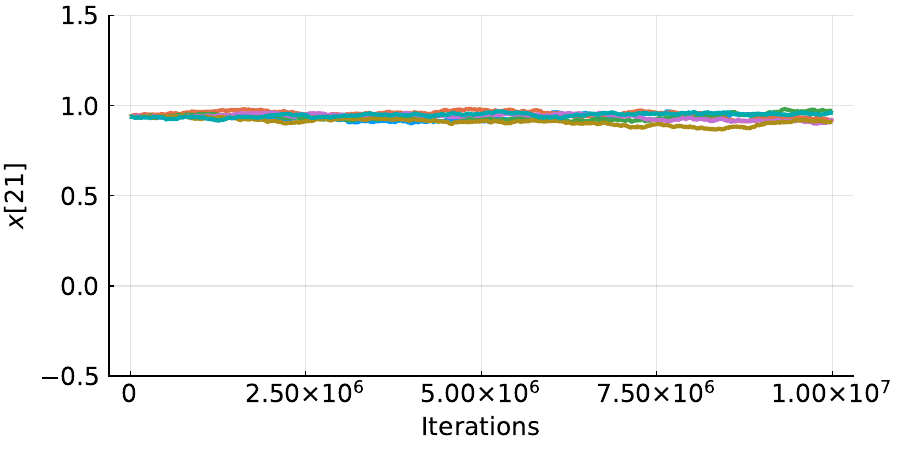} 
    		\caption{Original}
    		\label{fig:deblurring_model4_traces_x21_original_AM_MAP}
  	\end{subfigure}%
	~	 
	\begin{subfigure}[b]{0.32\textwidth}
		\includegraphics[width=\textwidth]{%
      		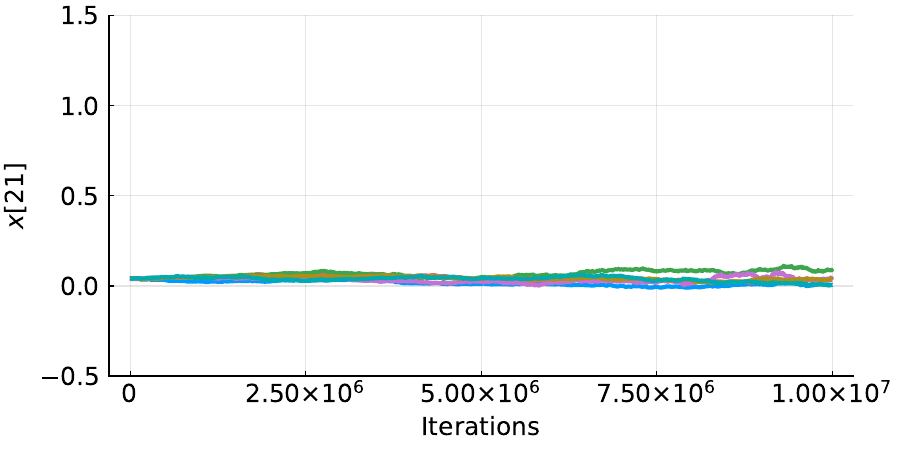} 
    		\caption{Approach from \cite{calvetti2024computationally}}
    		\label{fig:deblurring_model4_traces_x21_Calvetti2024_AM_MAP}
  	\end{subfigure}%
	~	 
	\begin{subfigure}[b]{0.32\textwidth}
		\includegraphics[width=\textwidth]{%
      		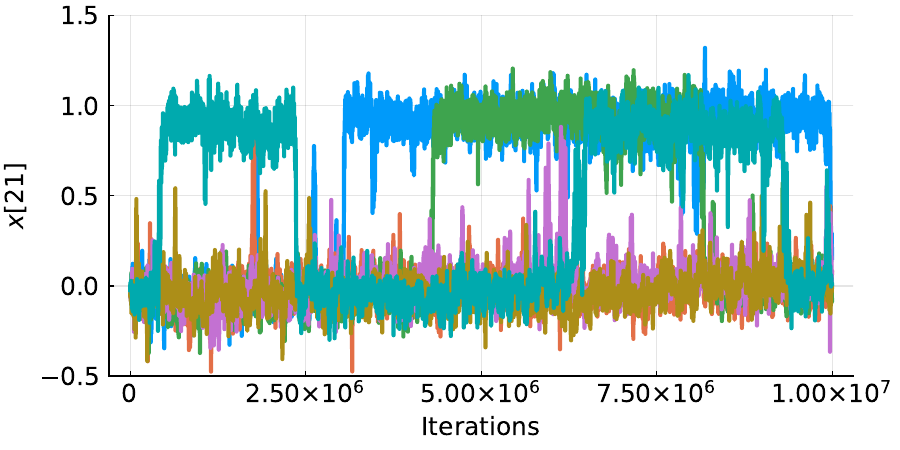} 
    		\caption{Prior-normalized}
    		\label{fig:deblurring_model4_traces_x21_priorNormalized_AM_MAP}
  	\end{subfigure}%
	\\
	\begin{subfigure}[b]{0.32\textwidth}
		\includegraphics[width=\textwidth]{%
      		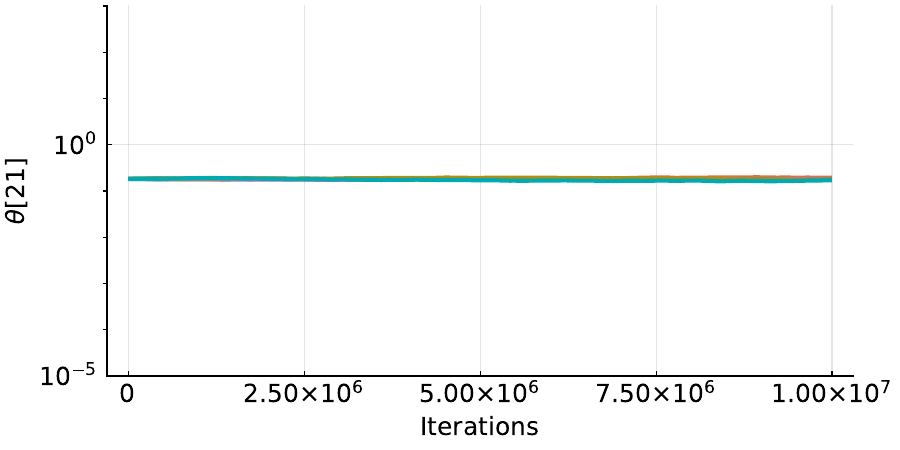} 
    		\caption{Original}
    		\label{fig:deblurring_model4_traces_theta21_original_AM_MAP}
  	\end{subfigure}%
	~	 
	\begin{subfigure}[b]{0.32\textwidth}
		\includegraphics[width=\textwidth]{%
      		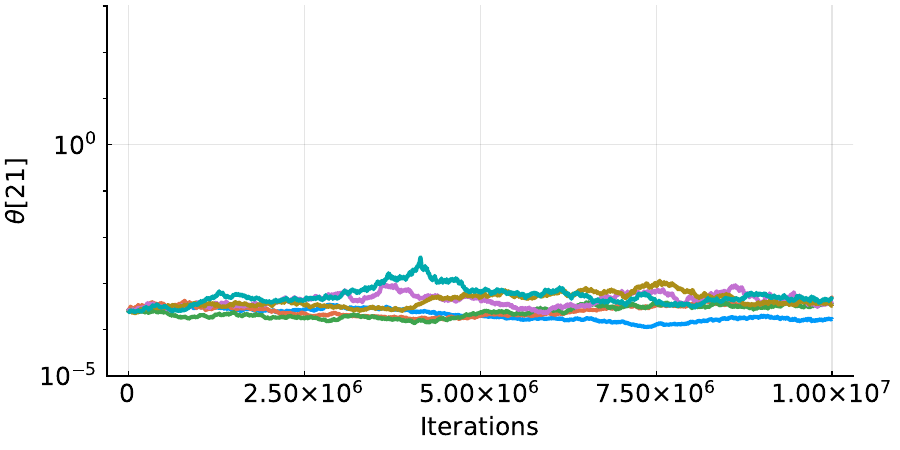} 
    		\caption{Approach from \cite{calvetti2024computationally}}
    		\label{fig:deblurring_model4_traces_theta21_Calvetti2024_AM_MAP}
  	\end{subfigure}%
	~	 
	\begin{subfigure}[b]{0.32\textwidth}
		\includegraphics[width=\textwidth]{%
      		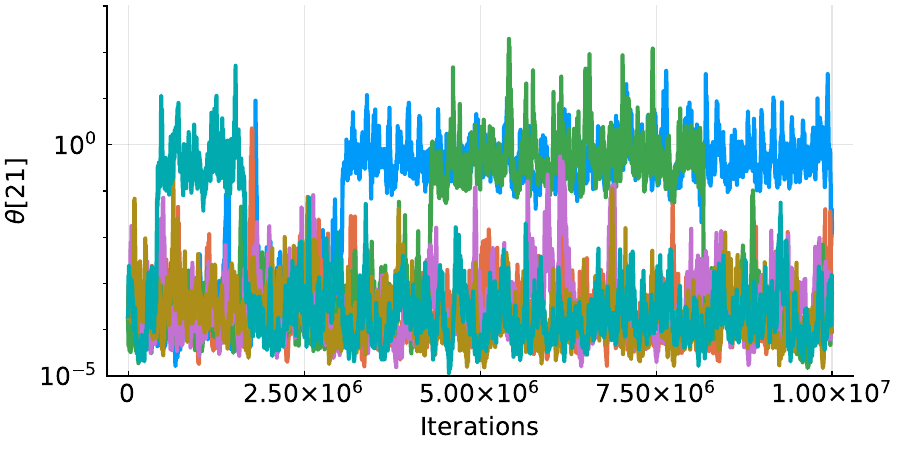} 
    		\caption{Prior-normalized}
    		\label{fig:deblurring_model4_traces_theta21_priorNormalized_AM_MAP}
  	\end{subfigure}%
  	\caption{ 
		Traces for the $x_{21}$- and $\theta_{21}$-samples and $r=-1$. 
		Note that $t_{21} \approx 0.164$ is slightly left of the jump at $t = 0.17$ and yields $g(t_{21}) = 0$.
		We initialized the AM algorithm with MAP estimates. 
	}
  	\label{fig:deblurring_traces_21_AM_MAP}
\end{figure} 

\begin{figure}[tb]
	\centering 
	\begin{subfigure}[b]{0.32\textwidth}
		\includegraphics[width=\textwidth]{%
      		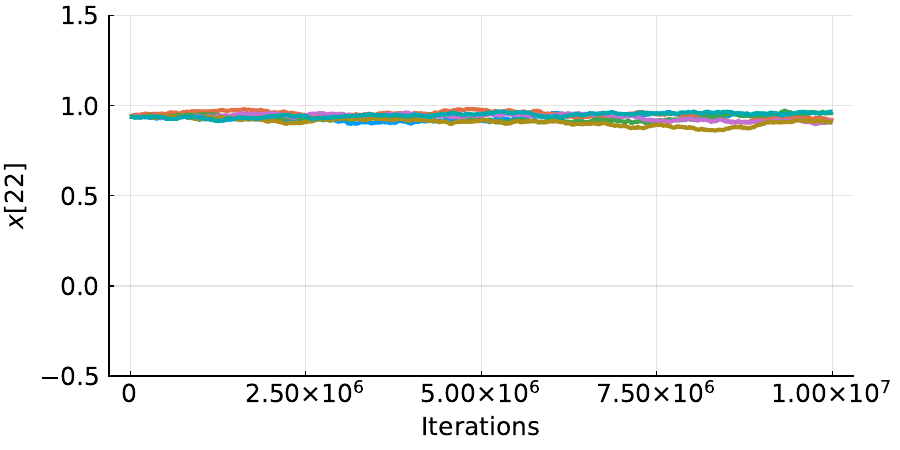} 
    		\caption{Original}
    		\label{fig:deblurring_model4_traces_x22_original_AM_MAP}
  	\end{subfigure}%
	~	 
	\begin{subfigure}[b]{0.32\textwidth}
		\includegraphics[width=\textwidth]{%
      		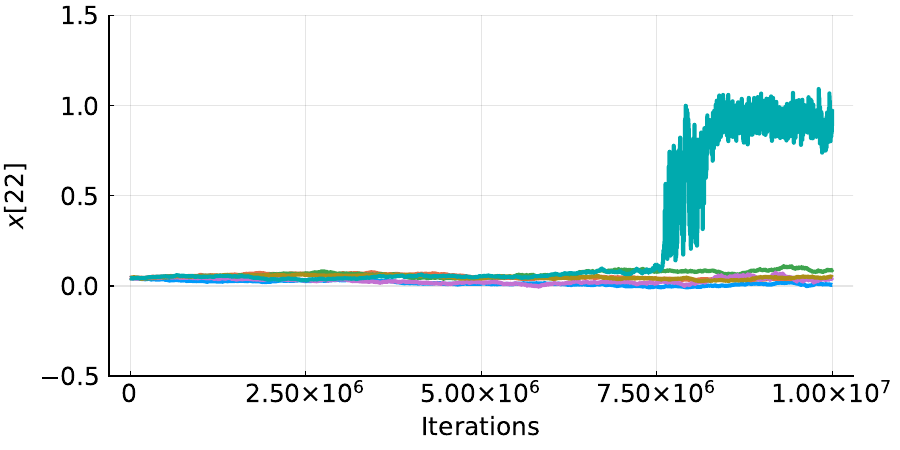} 
    		\caption{Approach from \cite{calvetti2024computationally}}
    		\label{fig:deblurring_model4_traces_x22_Calvetti2024_AM_MAP}
  	\end{subfigure}%
	~	 
	\begin{subfigure}[b]{0.32\textwidth}
		\includegraphics[width=\textwidth]{%
      		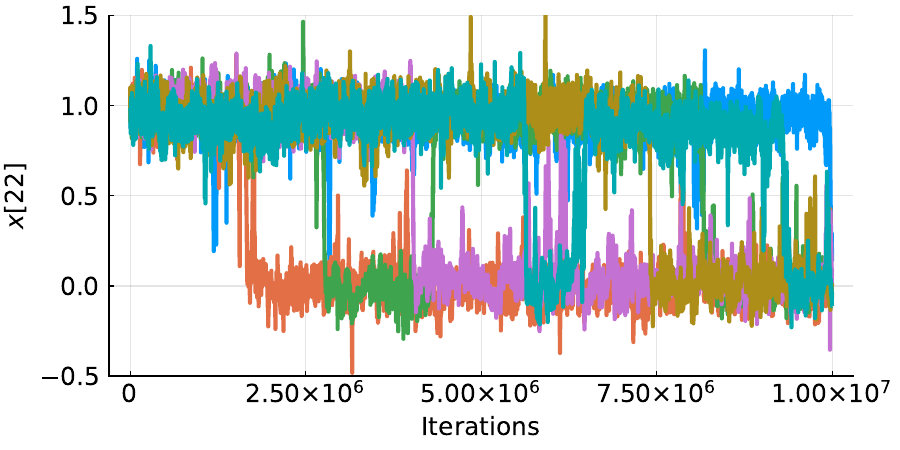} 
    		\caption{Prior-normalized}
    		\label{fig:deblurring_model4_traces_x22_priorNormalized_AM_MAP}
  	\end{subfigure}%
	\\
	\begin{subfigure}[b]{0.32\textwidth}
		\includegraphics[width=\textwidth]{%
      		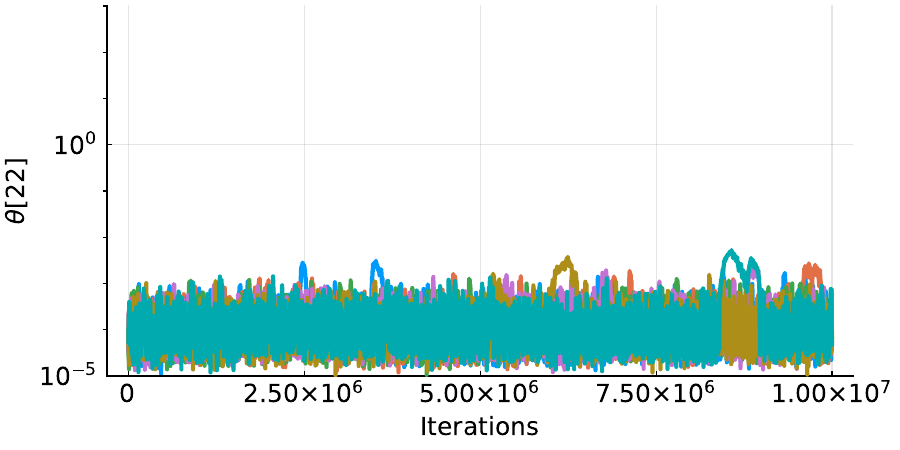} 
    		\caption{Original}
    		\label{fig:deblurring_model4_traces_theta22_original_AM_MAP}
  	\end{subfigure}%
	~	 
	\begin{subfigure}[b]{0.32\textwidth}
		\includegraphics[width=\textwidth]{%
      		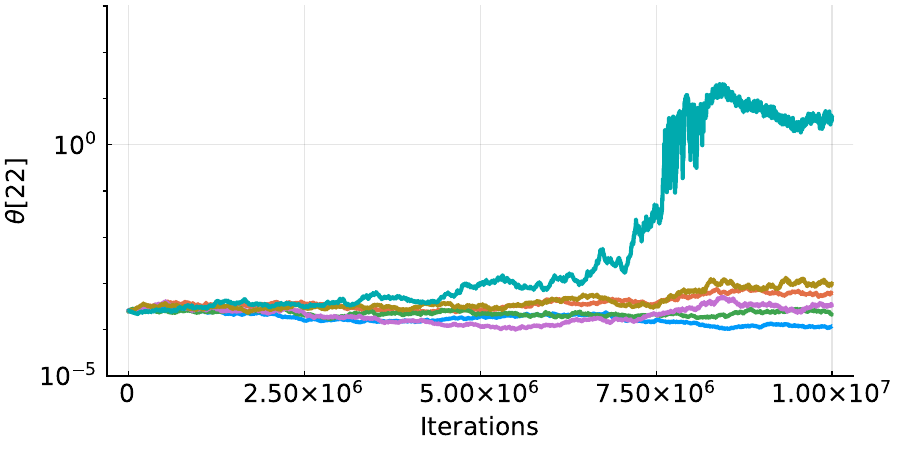} 
    		\caption{Approach from \cite{calvetti2024computationally}}
    		\label{fig:deblurring_model4_traces_theta22_Calvetti2024_AM_MAP}
  	\end{subfigure}%
	~	 
	\begin{subfigure}[b]{0.32\textwidth}
		\includegraphics[width=\textwidth]{%
      		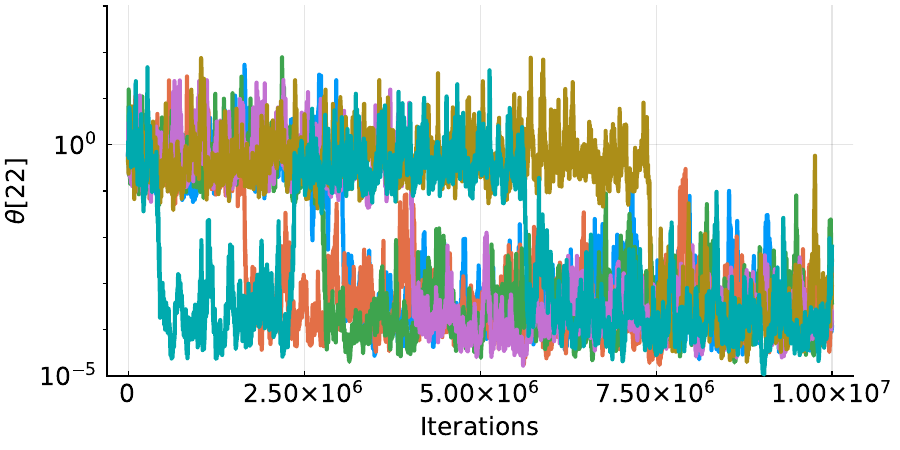} 
    		\caption{Prior-normalized}
    		\label{fig:deblurring_model4_traces_theta22_priorNormalized_AM_MAP}
  	\end{subfigure}%
  	\caption{ 
		Traces for the $x_{22}$- and $\theta_{22}$-samples and $r=-1$. 
		Note that $t_{22} \approx 0.172$ is slightly right of the jump at $t = 0.17$ and yields $g(t_{22}) = 1$.
		We initialized the AM algorithm with MAP estimates. 
	}
  	\label{fig:deblurring_traces_22_AM_MAP}
\end{figure}

Moreover, \cref{fig:deblurring_traces_21_AM_MAP,fig:deblurring_traces_22_AM_MAP} show the trace plots for the $x_{i}$- and $\theta_i$-samples for $i=21$ (\cref{fig:deblurring_traces_21_AM_MAP}), $i=22$ (\cref{fig:deblurring_traces_22_AM_MAP}), and $r = -1$. 
Note that $t_{21} \approx 0.164$ is slightly to the left of the first jump discontinuity at $t = 0.17$ and yields $g(t_{21}) = 0$.
At the same time, $t_{22} \approx 0.172$ is slightly to the right of the first jump discontinuity at $t = 0.17$ and yields $g(t_{22}) = 1$.
\cref{fig:deblurring_traces_21_AM_MAP,fig:deblurring_traces_22_AM_MAP} provides further evidence that the AM sampler explores the original posterior and the transformed posterior from \cite{calvetti2024computationally} only locally. 
Specifically, we see that the chains for the original posterior and the transformed posterior from \cite{calvetti2024computationally} exhibit minimal movement from their initialization. 
The only exception is one of the MCMC chains for the transformed posterior from \cite{calvetti2024computationally} illustrated in \cref{fig:deblurring_model4_traces_x22_Calvetti2024_AM_MAP,fig:deblurring_model4_traces_theta22_Calvetti2024_AM_MAP}. 
In contrast, the chains for the prior-normalized posterior better explore the two distinct high-density regions around zero and one for $x_{21}$ and $x_{22}$ as well as the corresponding small- and large-value regimes for $\theta_{21}$ and $\theta_{22}$.
This outcome is intuitively desirable, as $t_{21} \approx 0.164$ is slightly left of the first jump discontinuity at $t=0.17$. 
While $g(t_{21}) = 0$, the Bayesian model should reflect the uncertainty that the jump might occur at the present grid point, meaning both $x_{21} = 0$ and $x_{21} = 1$ should lie in high-density posterior regions.
Similarly, $t_{22} \approx 0.172$ is slightly right of the first jump discontinuity at $t=0.17$ with $g(t_{22}) = 1$. 
Hence, also in this case, both $x_{22} = 0$ and $x_{22} = 1$ should lie in high-density posterior regions.

\subsubsection{Initializing with prior samples}
\label{subsub:tests_deconvolution_prior}

\begin{figure}[tb]
	\centering 
	\begin{subfigure}[b]{0.32\textwidth}
		\includegraphics[width=\textwidth]{%
      		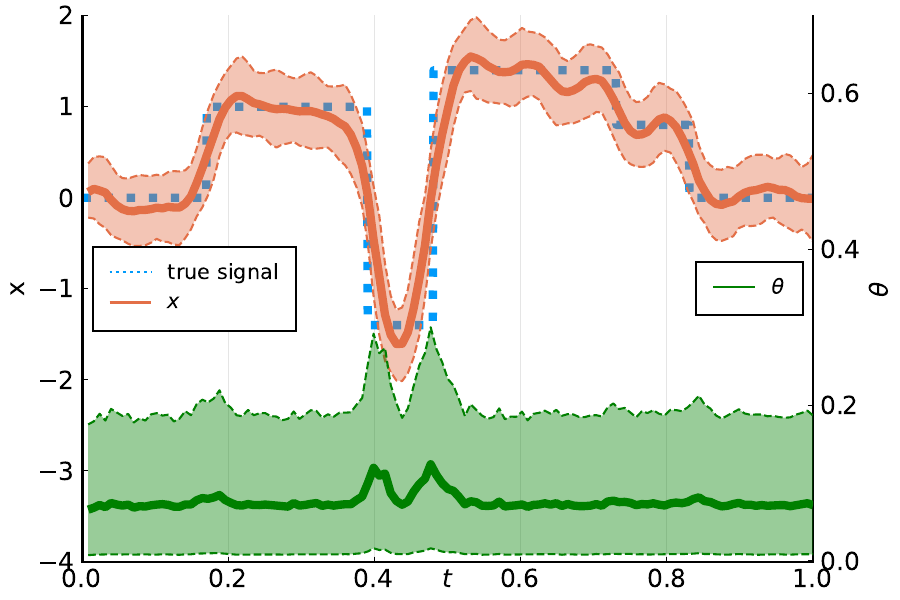} 
    		\caption{Original, $r=1$}
    		\label{fig:deblurring_model1_UQquantile_original_AM_prior}
  	\end{subfigure}%
	~
	\begin{subfigure}[b]{0.32\textwidth}
		\includegraphics[width=\textwidth]{%
      		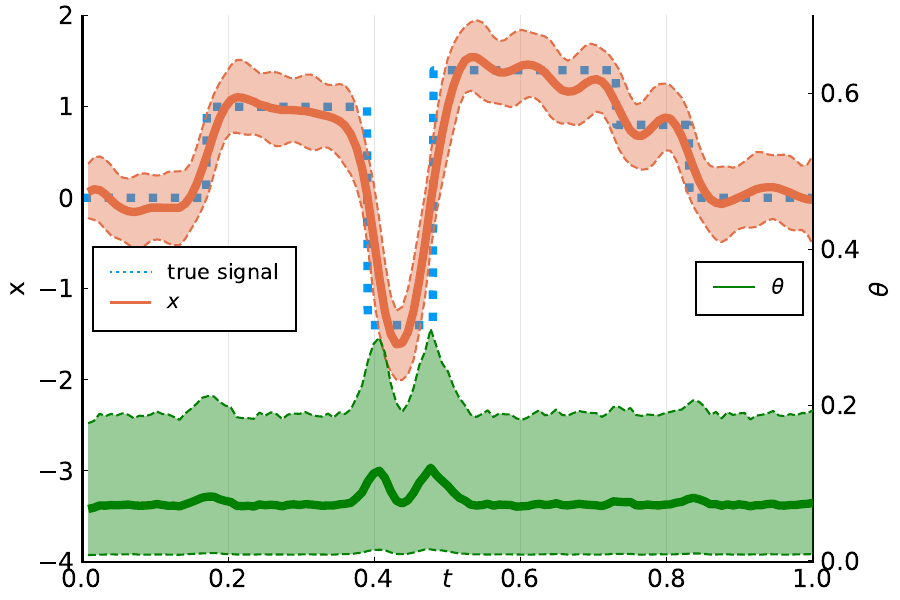} 
    		\caption{Approach from \cite{calvetti2024computationally}, $r=1$}
    		\label{fig:deblurring_model1_UQquantile_Calvetti2024_AM_prior}
  	\end{subfigure}%
	~
	\begin{subfigure}[b]{0.32\textwidth}
		\includegraphics[width=\textwidth]{%
      		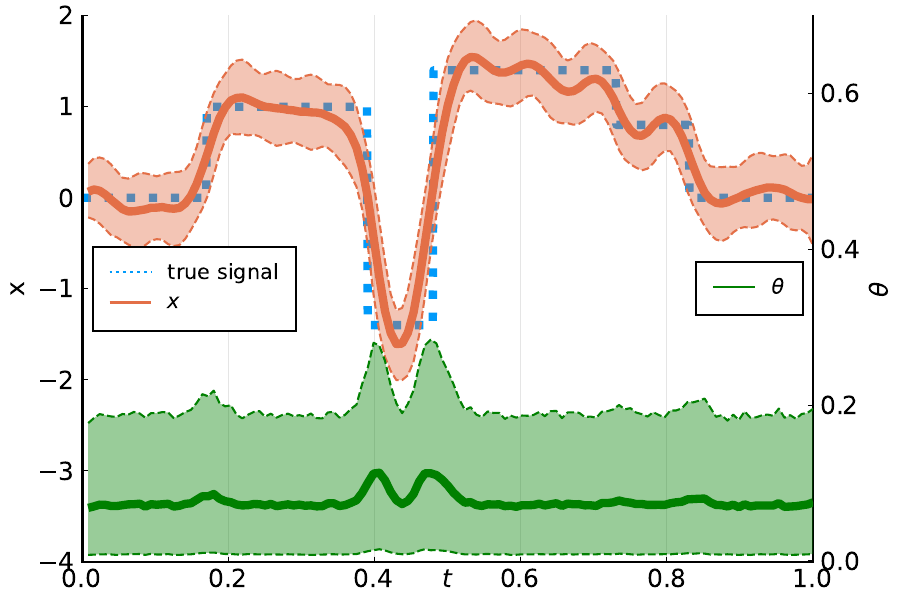} 
    		\caption{Prior-normalized (our), $r=1$}
    		\label{fig:deblurring_model1_UQquantile_priorNormalized_AM_prior}
  	\end{subfigure}%
	\\
	\begin{subfigure}[b]{0.32\textwidth}
		\includegraphics[width=\textwidth]{%
      		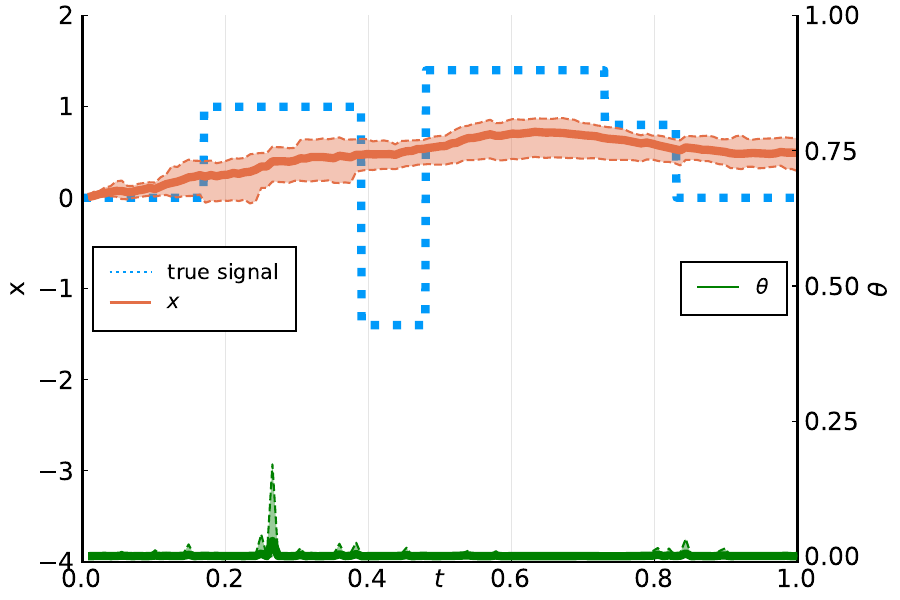} 
    		\caption{Original, $r=1$}
    		\label{fig:deblurring_model4_UQquantile_original_AM_prior}
  	\end{subfigure}%
	~
	\begin{subfigure}[b]{0.32\textwidth}
		\includegraphics[width=\textwidth]{%
      		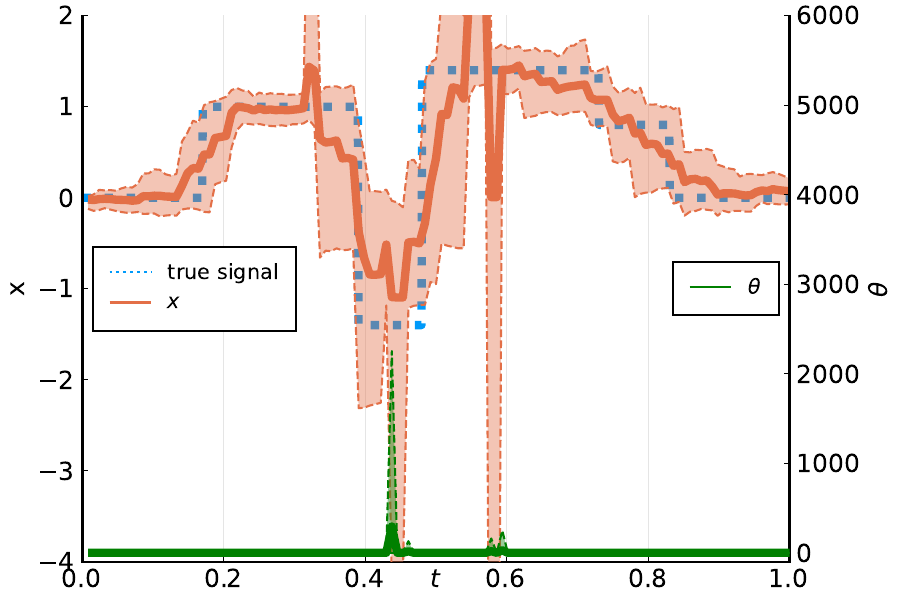} 
    		\caption{Approach from \cite{calvetti2024computationally}, $r=1$}
    		\label{fig:deblurring_model4_UQquantile_Calvetti2024_AM_prior}
  	\end{subfigure}%
	~
	\begin{subfigure}[b]{0.32\textwidth}
		\includegraphics[width=\textwidth]{%
      		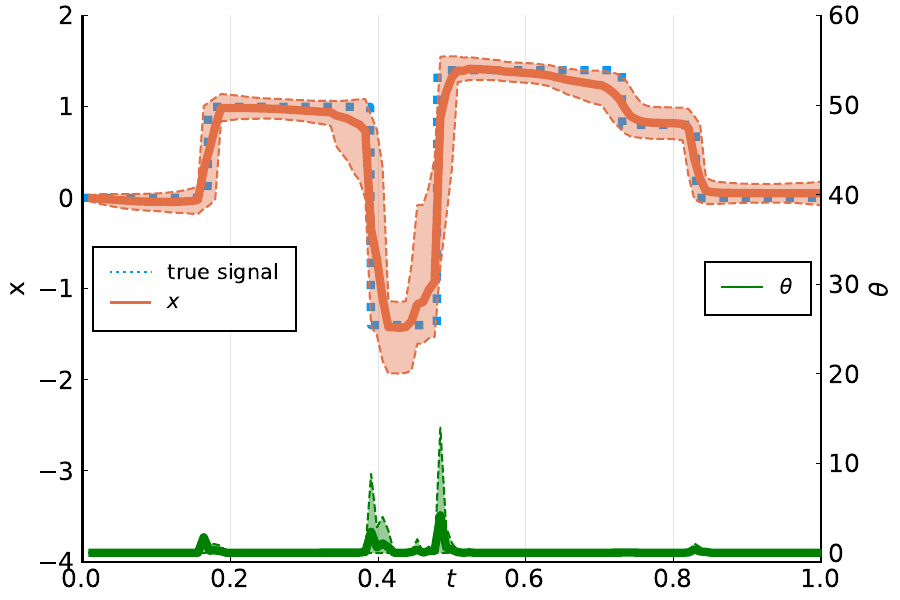} 
    		\caption{Prior-normalized (our), $r=1$}
    		\label{fig:deblurring_model4_UQquantile_priorNormalized_AM_prior}
  	\end{subfigure}%
	\caption{ 
		Mean and $90\%$ quantile ranges of the $\mathbf{x}$- and $\boldsymbol{\theta}$-samples of the original posterior (left), the transformed posterior from \cite{calvetti2024computationally} (middle), and the proposed prior-normalized posterior (right) for $r = \pm 1$. 
		We used the AM algorithm initialized with random draws from the respective prior.   
	}
  	\label{fig:deblurring_UQ_AM_prior}
\end{figure} 

We replicate the experiments from 
\Cref{subsub:tests_deconvolution_MAP}, but this time, we independently initialize the MCMC chains using randomly generated samples from the SBL and standard normal prior, respectively.
\cref{fig:deblurring_UQ_AM_prior} presents the sample means and $90\%$ quantile ranges from the original posterior, the transformed posterior from \cite{calvetti2024computationally}, and our prior-normalized posterior for $r=\pm1$.
For $r=1$, the AM method again performs similarly on all three posteriors, yielding comparable estimates for the means and quantile ranges. 
Noticeable differences appear for $r = -1$, however. 
This time, the $x$-mean for the original posterior in \cref{fig:deblurring_model4_UQquantile_original_AM_prior} deviates significantly from the true underlying signal. 
The $x$-mean for the transformed posterior from \cite{calvetti2024computationally} follows the true underlying signal only roughly, with visible differences. 
In contrast, the $x$-mean and -quantile ranges for the prior-normalized posterior in \cref{fig:deblurring_model4_UQquantile_priorNormalized_AM_prior} follow the true underlying signal significantly more accurately. 
Furthermore, they are qualitatively similar to those in \cref{fig:deblurring_model4_UQquantile_priorNormalized_AM_MAP}, indicating that the prior-normalized posterior is effectively explored by the AM sampler, regardless of how the chains are initialized. 

\begin{figure}[tb]
	\centering 
	\begin{subfigure}[b]{0.32\textwidth}
		\includegraphics[width=\textwidth]{%
      		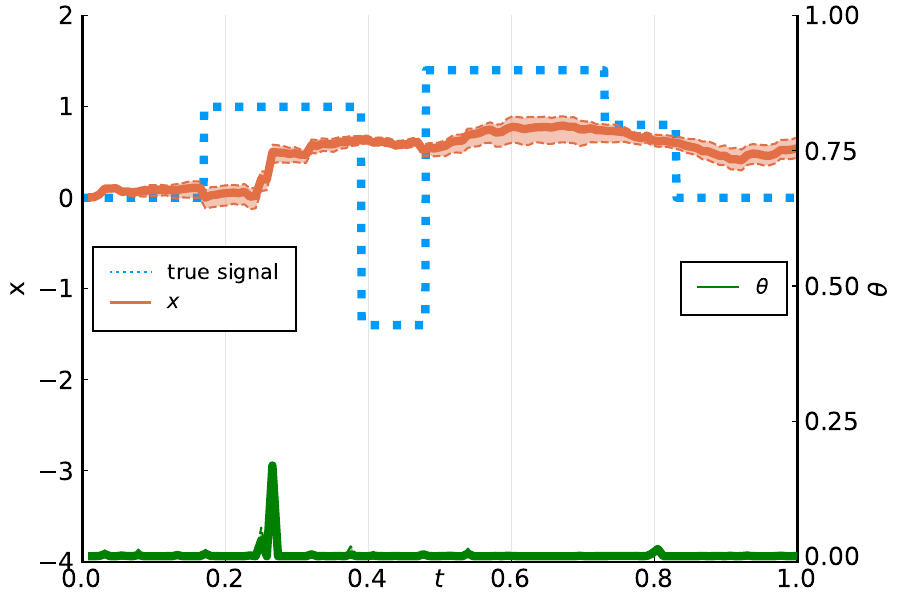} 
    		\caption{Original, $1$st chain}
    		\label{deblurring_model4_UQquantile_original_AM_prior_chain1}
	\end{subfigure}%
	~
	\begin{subfigure}[b]{0.32\textwidth}
		\includegraphics[width=\textwidth]{%
      		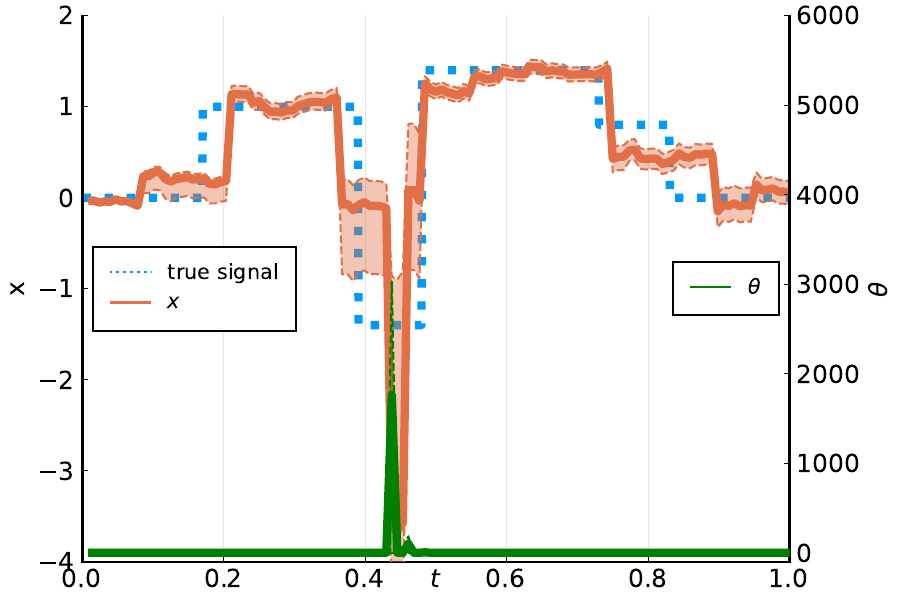} 
    		\caption{Approach from \cite{calvetti2024computationally}, $1$st chain}
    		\label{deblurring_model4_UQquantile_Calvetti2024_AM_prior_chain1}
	\end{subfigure}%
	~
	\begin{subfigure}[b]{0.32\textwidth}
		\includegraphics[width=\textwidth]{%
      		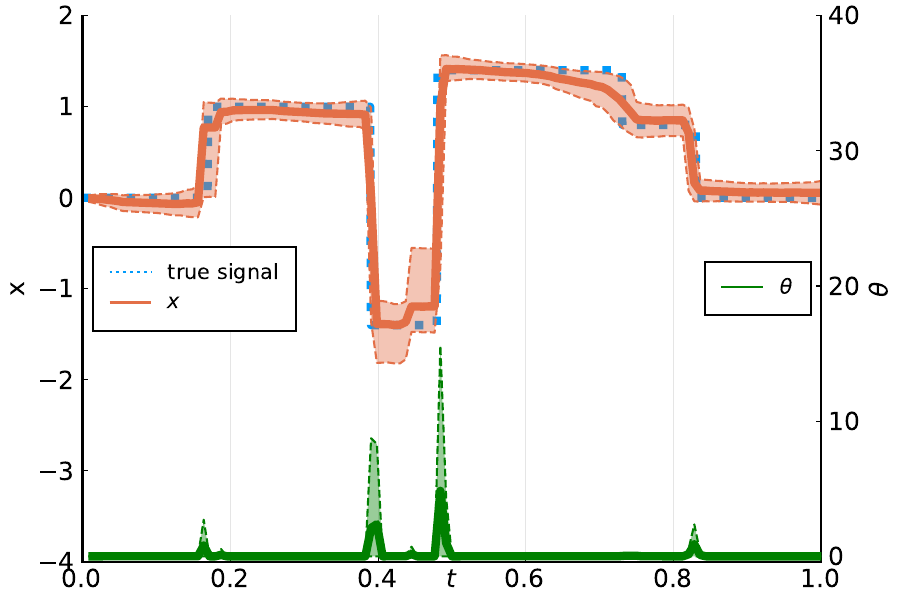} 
    		\caption{Prior-normalized, $1$st chain}
    		\label{deblurring_model4_UQquantile_priorNormalized_AM_prior_chain1}
	\end{subfigure}%
	\\ 
	\begin{subfigure}[b]{0.32\textwidth}
		\includegraphics[width=\textwidth]{%
      		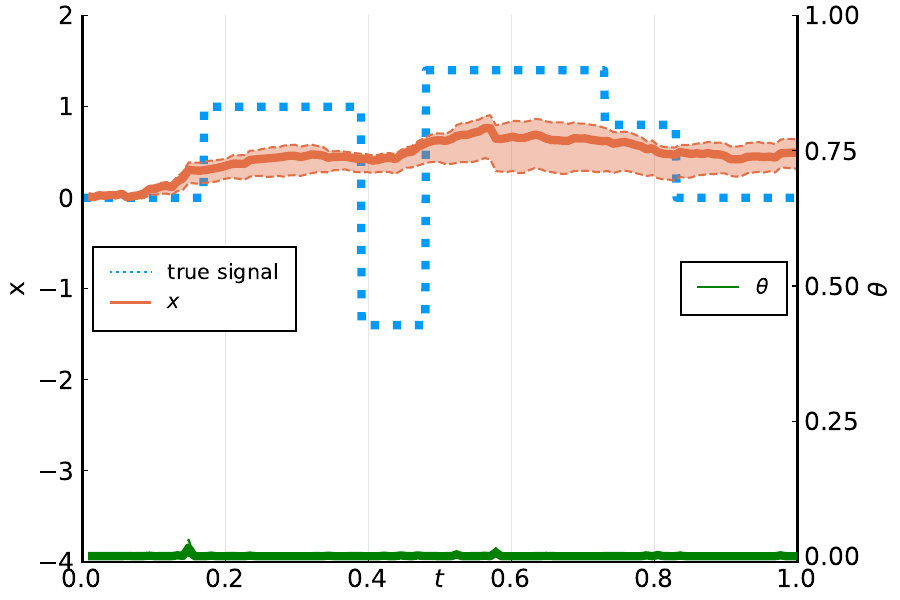} 
    		\caption{Original, $2$nd chain}
    		\label{deblurring_model4_UQquantile_original_AM_prior_chain2}
	\end{subfigure}%
	~
	\begin{subfigure}[b]{0.32\textwidth}
		\includegraphics[width=\textwidth]{%
      		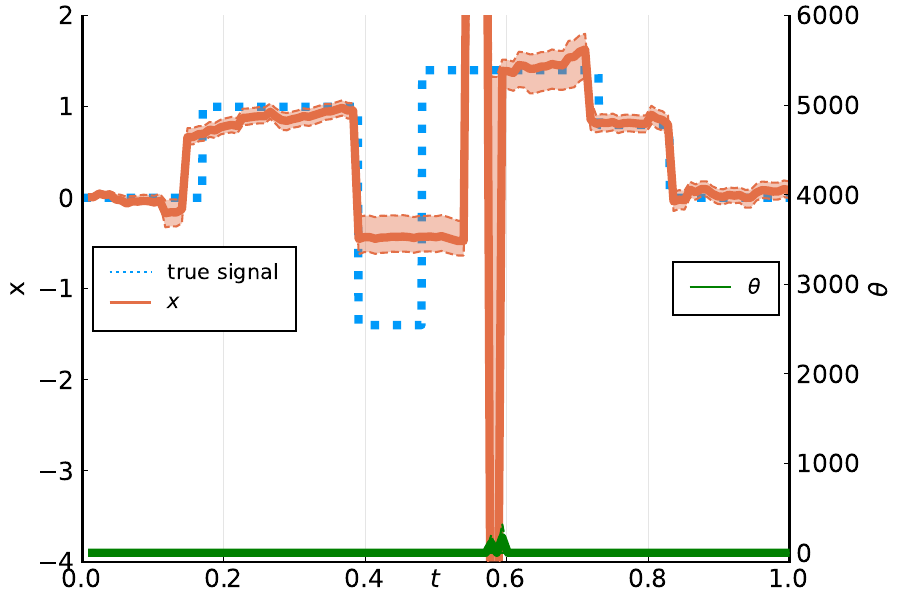} 
    		\caption{Approach from \cite{calvetti2024computationally}, $2$nd chain}
    		\label{deblurring_model4_UQquantile_Calvetti2024_AM_prior_chain2}
	\end{subfigure}%
	~
	\begin{subfigure}[b]{0.32\textwidth}
		\includegraphics[width=\textwidth]{%
      		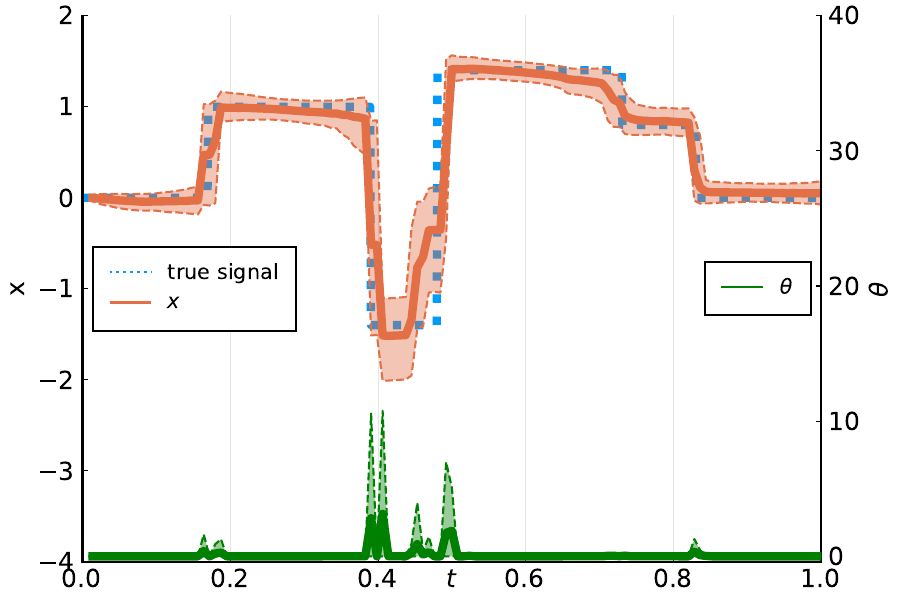} 
    		\caption{Prior-normalized, $2$nd chain}
    		\label{deblurring_model4_UQquantile_priorNormalized_AM_prior_chain2}
	\end{subfigure}%
	\\
	\begin{subfigure}[b]{0.32\textwidth}
		\includegraphics[width=\textwidth]{%
      		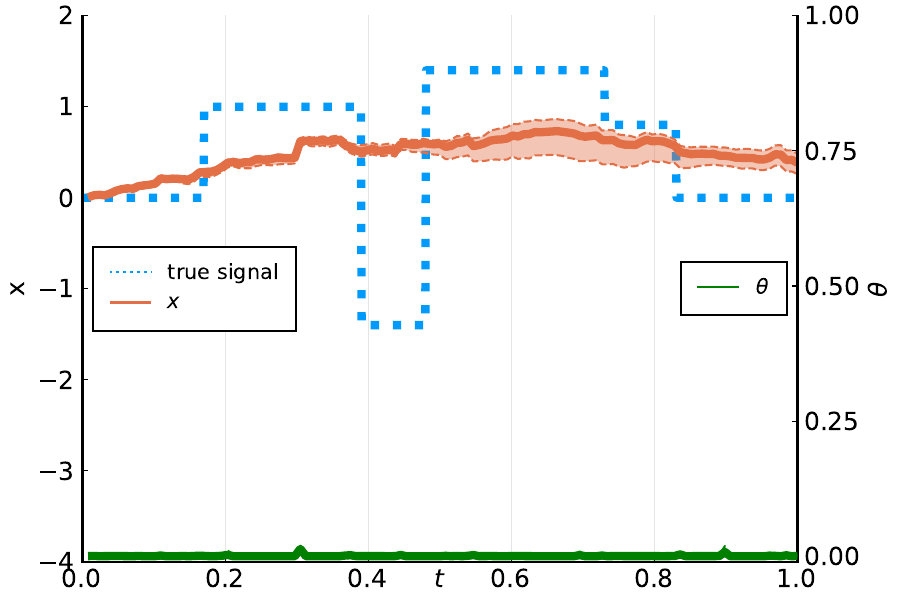} 
    		\caption{Original, $3$rd chain}
    		\label{deblurring_model4_UQquantile_original_AM_prior_chain3}
	\end{subfigure}%
	~
	\begin{subfigure}[b]{0.32\textwidth}
		\includegraphics[width=\textwidth]{%
      		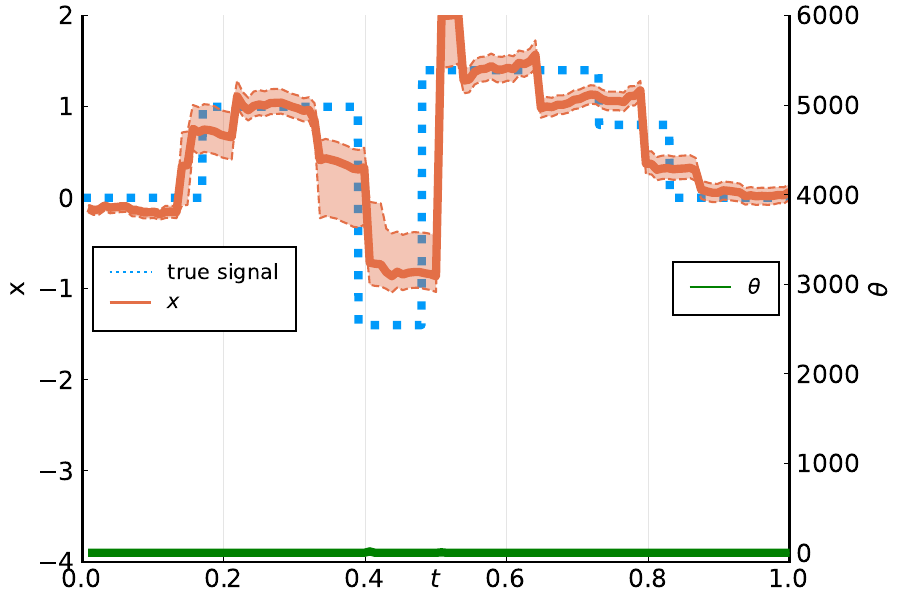} 
    		\caption{Approach from \cite{calvetti2024computationally}, $3$rd chain}
    		\label{deblurring_model4_UQquantile_Calvetti2024_AM_prior_chain3}
	\end{subfigure}%
	~
	\begin{subfigure}[b]{0.32\textwidth}
		\includegraphics[width=\textwidth]{%
      		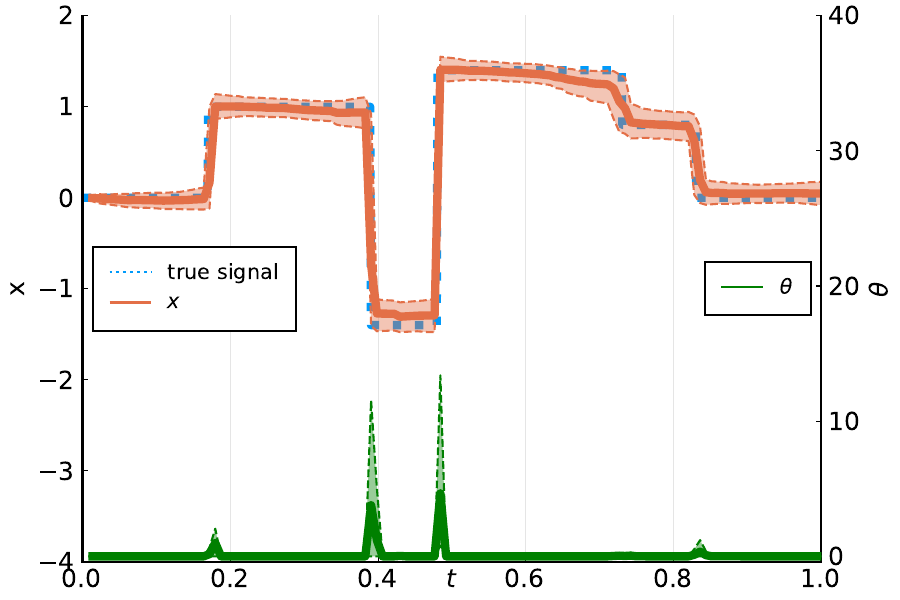} 
    		\caption{Prior-normalized, $3$rd chain}
    		\label{deblurring_model4_UQquantile_priorNormalized_AM_prior_chain3}
	\end{subfigure}%
  	\caption{ 
		Mean and $90\%$ quantile range of the $\mathbf{x}$- and $\boldsymbol{\theta}$-samples for the first three MCMC chains of the original posterior, the transformed posterior from \cite{calvetti2024computationally}, and the prior-normalized posterior for $r=-1$.
		We used the AM algorithm initialized with random draws from the respective prior. 
	}
  	\label{fig:deblurring_chains_AM_prior}
\end{figure}

\cref{fig:deblurring_chains_AM_prior} additionally presents the means and quantile ranges for the first three chains across all posteriors for $r =- 1$. 
Notably, each individual chain for the original posterior fails to produce reasonable estimates. 
At the same time, each individual chain for the transformed posterior from \cite{calvetti2024computationally} follows the true underlying signal only fairly inaccurately, showing visible differences.
In contrast, each individual chain from the prior-normalized posterior provides reasonable estimates. 

\begin{figure}[tb]
	\centering 
	\begin{subfigure}[b]{0.32\textwidth}
		\includegraphics[width=\textwidth]{%
      		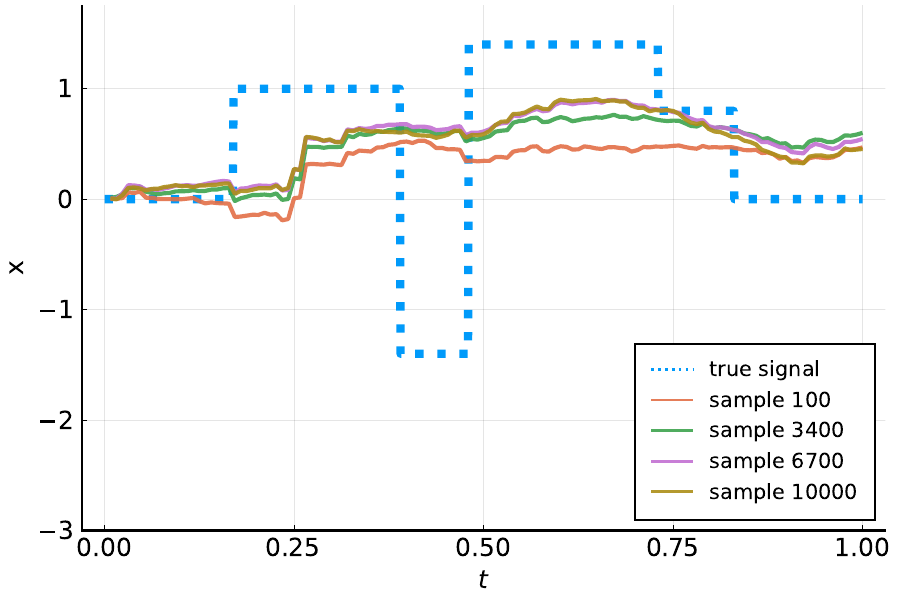} 
    		\caption{Original}
    		\label{fig:deblurring_model4_xSamples_original_AM_prior}
  	\end{subfigure}
	~
	\begin{subfigure}[b]{0.32\textwidth}
		\includegraphics[width=\textwidth]{%
      		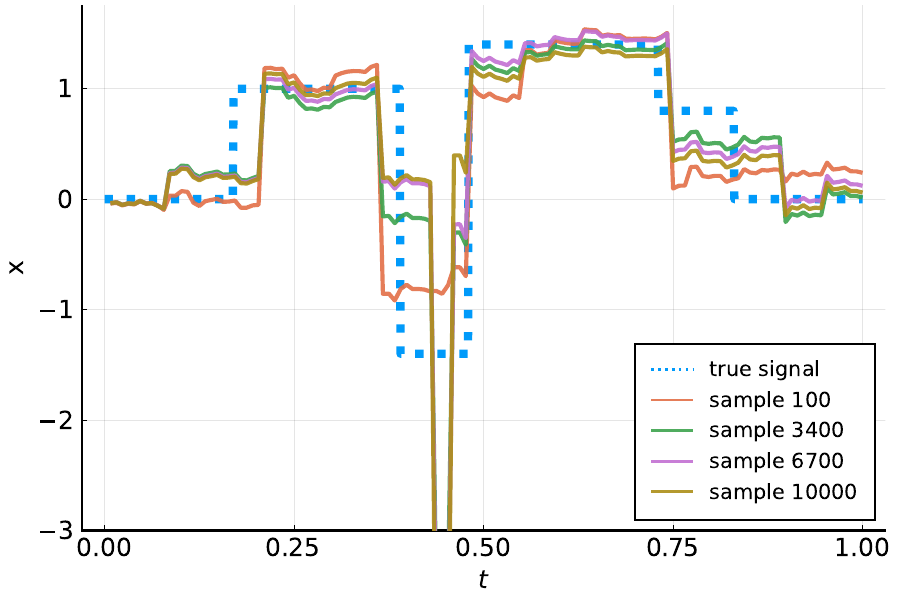} 
    		\caption{Approach from \cite{calvetti2024computationally}}
    		\label{fig:deblurring_model4_xSamples_Calvetti2024_AM_prior}
  	\end{subfigure}%
	~
	\begin{subfigure}[b]{0.32\textwidth}
		\includegraphics[width=\textwidth]{%
      		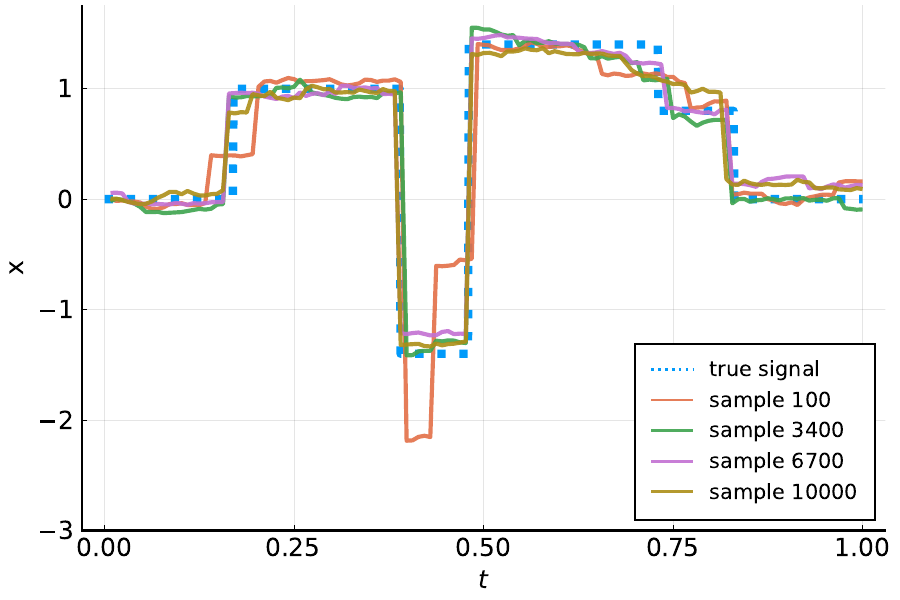} 
    		\caption{Prior-normalized}
    		\label{fig:deblurring_model4_xSamples_priorNormalized_AM_prior}
  	\end{subfigure}%
  	\caption{ 
		Individual $x$-samples of the first MCMC chain from the original posterior, the transformed posterior from \cite{calvetti2024computationally}, and the prior-normalized posterior for $r=-1$. 
		We used the AM algorithm initialized with random prior draws. 
	}
  	\label{fig:deblurring_samples_AM_prior}
\end{figure} 

\begin{figure}[tb]
	\centering 
	\begin{subfigure}[b]{0.32\textwidth}
		\includegraphics[width=\textwidth]{%
      		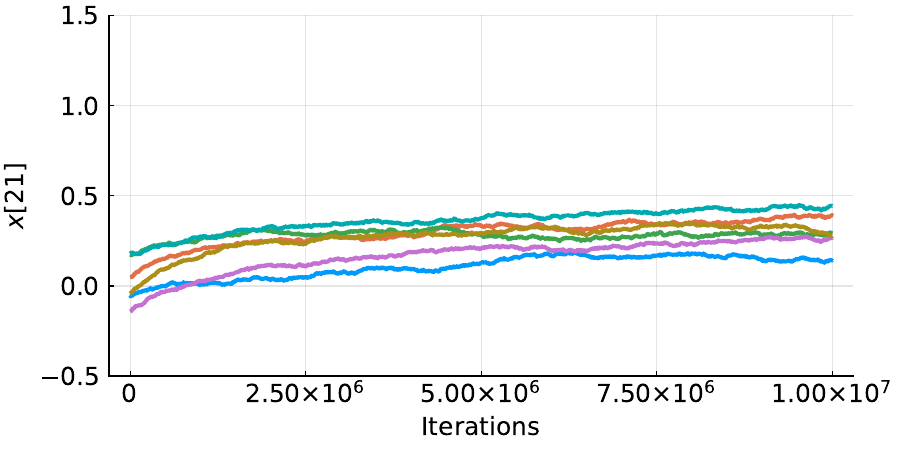} 
    		\caption{Original}
    		\label{fig:deblurring_model4_traces_x21_original_AM_prior}
  	\end{subfigure}%
	~	 
	\begin{subfigure}[b]{0.32\textwidth}
		\includegraphics[width=\textwidth]{%
      		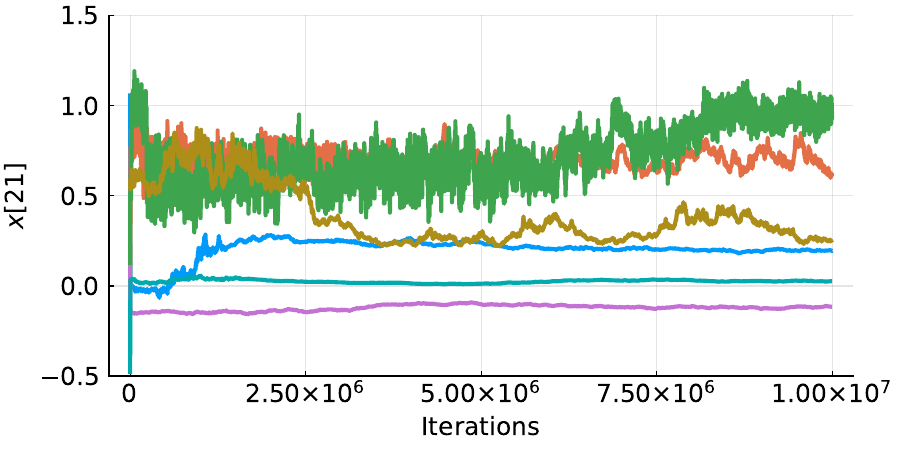} 
    		\caption{Approach from \cite{calvetti2024computationally}}
    		\label{fig:deblurring_model4_traces_x21_Calvetti2024_AM_prior}
  	\end{subfigure}%
	~	 
	\begin{subfigure}[b]{0.32\textwidth}
		\includegraphics[width=\textwidth]{%
      		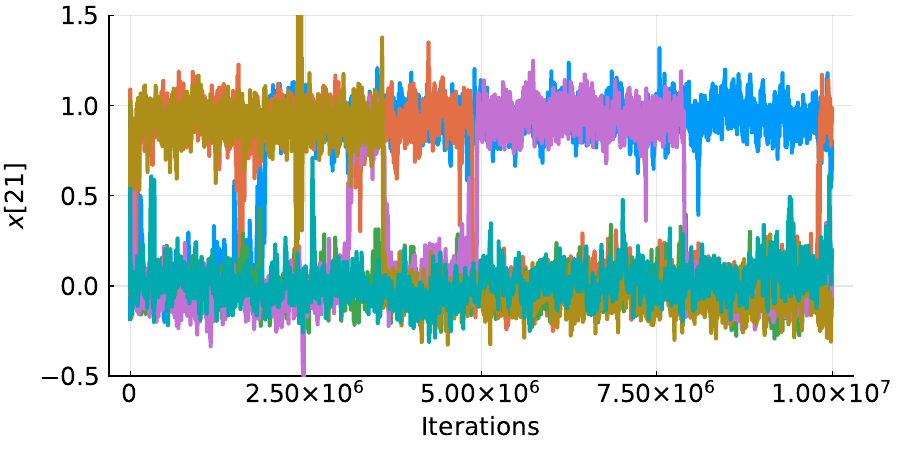} 
    		\caption{Prior-normalized}
    		\label{fig:deblurring_model4_traces_x21_priorNormalized_AM_prior}
  	\end{subfigure}%
	\\
	\begin{subfigure}[b]{0.32\textwidth}
		\includegraphics[width=\textwidth]{%
      		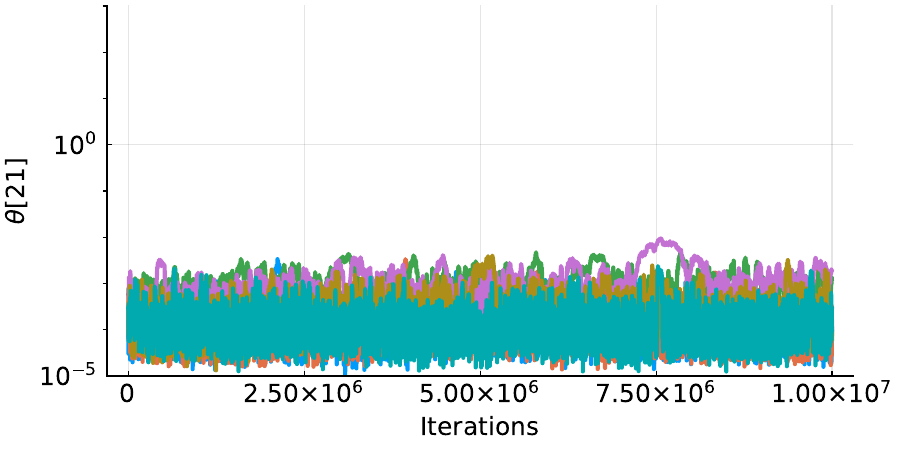} 
    		\caption{Original}
    		\label{fig:deblurring_model4_traces_theta21_original_AM_prior}
  	\end{subfigure}%
	~	 
	\begin{subfigure}[b]{0.32\textwidth}
		\includegraphics[width=\textwidth]{%
      		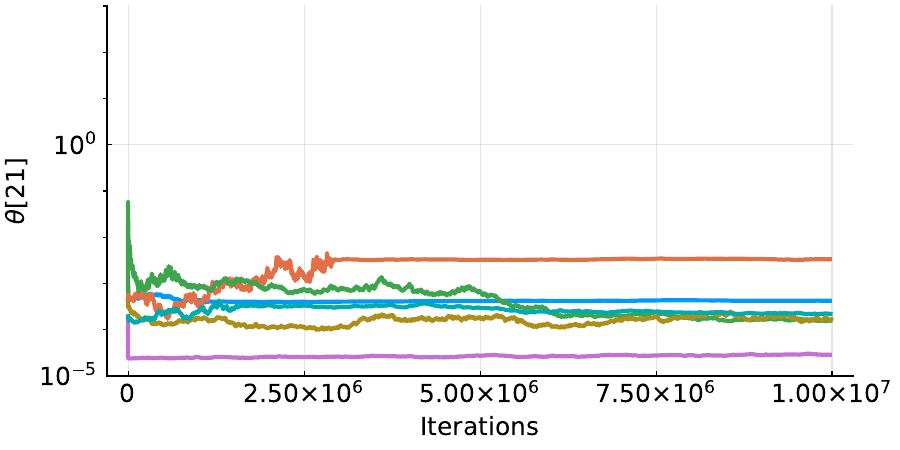} 
    		\caption{Approach from \cite{calvetti2024computationally}}
    		\label{fig:deblurring_model4_traces_theta21_Calvetti2024_AM_prior}
  	\end{subfigure}%
	~	 
	\begin{subfigure}[b]{0.32\textwidth}
		\includegraphics[width=\textwidth]{%
      		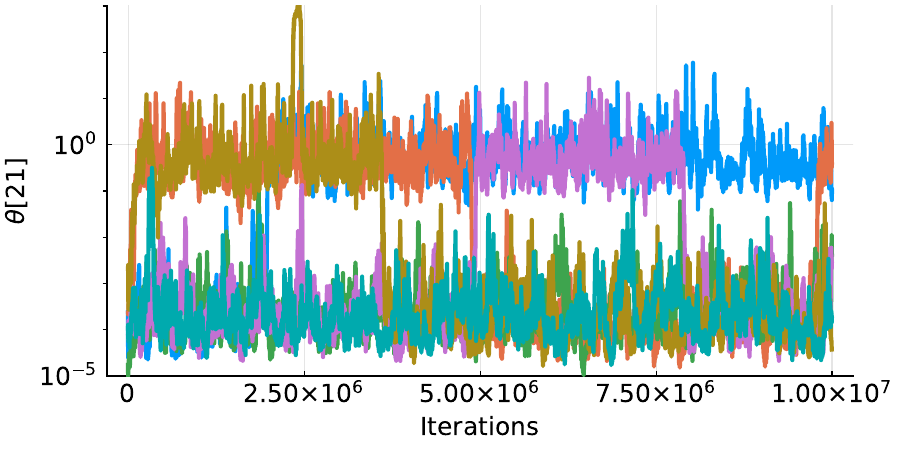} 
    		\caption{Prior-normalized}
    		\label{fig:deblurring_model4_traces_theta21_priorNormalized_AM_prior}
  	\end{subfigure}%
  	\caption{ 
		Traces for the $x_{21}$- and $\theta_{21}$-samples and $r=-1$. 
		Note that $t_{21} \approx 0.164$ is slightly to the left of the first jump discontinuity at $t = 0.17$ and yields $g(t_{21}) = 0$.
		We used the AM algorithm initialized with random draws from the respective prior. 
	}
  	\label{fig:deblurring_traces_21_AM_prior}
\end{figure} 

\begin{figure}[tb]
	\centering 
	\begin{subfigure}[b]{0.32\textwidth}
		\includegraphics[width=\textwidth]{%
      		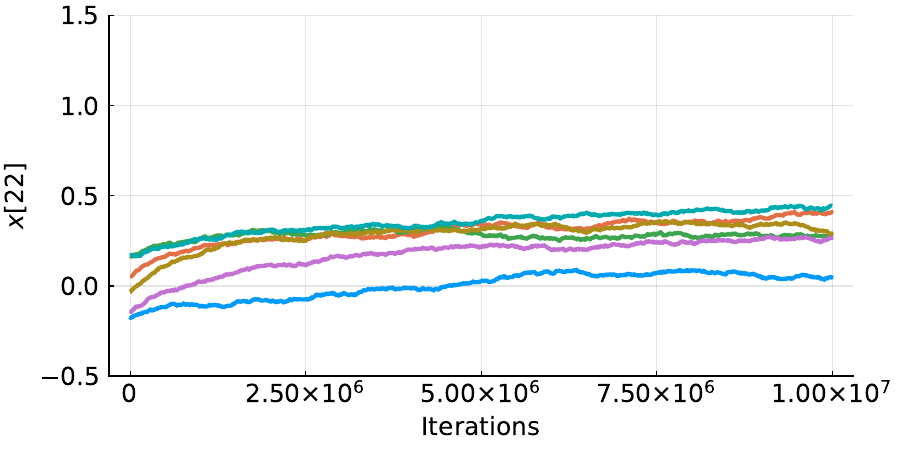} 
    		\caption{Original}
    		\label{fig:deblurring_model4_traces_x22_original_AM_prior}
  	\end{subfigure}%
	~	 
	\begin{subfigure}[b]{0.32\textwidth}
		\includegraphics[width=\textwidth]{%
      		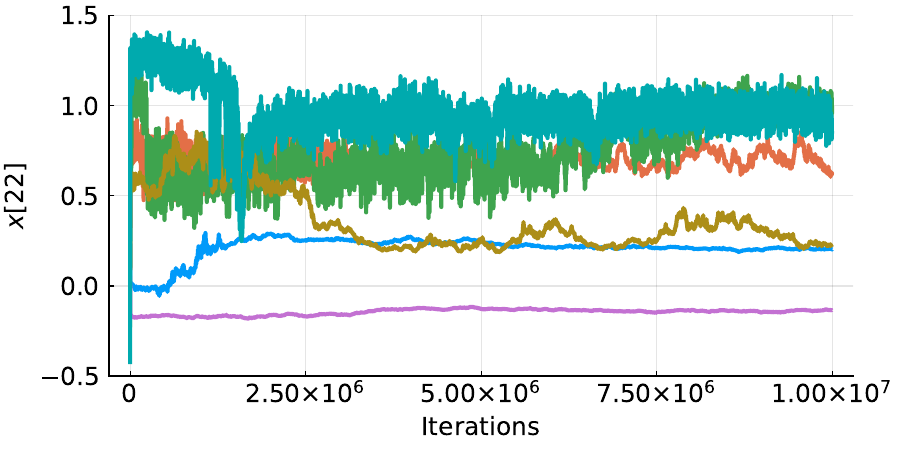} 
    		\caption{Approach from \cite{calvetti2024computationally}}
    		\label{fig:deblurring_model4_traces_x22_Calvetti2024_AM_prior}
  	\end{subfigure}%
	~	 
	\begin{subfigure}[b]{0.32\textwidth}
		\includegraphics[width=\textwidth]{%
      		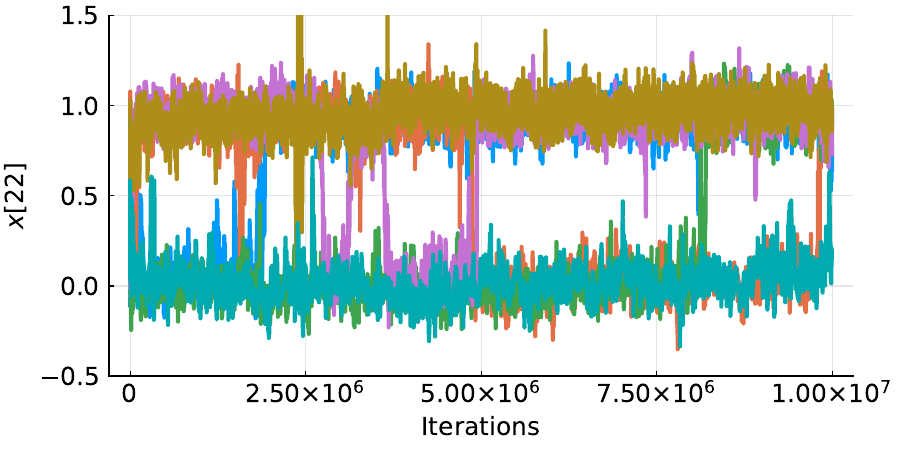} 
    		\caption{Prior-normalized}
    		\label{fig:deblurring_model4_traces_x22_priorNormalized_AM_prior}
  	\end{subfigure}%
	\\
	\begin{subfigure}[b]{0.32\textwidth}
		\includegraphics[width=\textwidth]{%
      		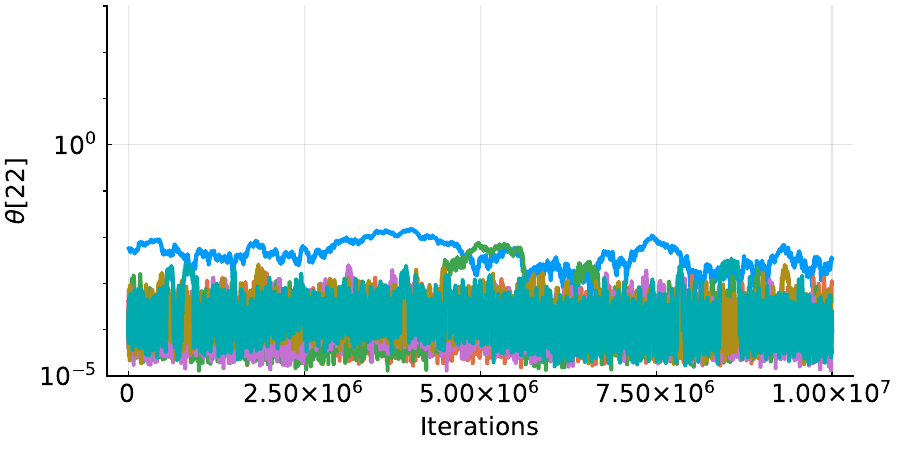} 
    		\caption{Original}
    		\label{fig:deblurring_model4_traces_theta22_original_AM_prior}
  	\end{subfigure}%
	~	 
	\begin{subfigure}[b]{0.32\textwidth}
		\includegraphics[width=\textwidth]{%
      		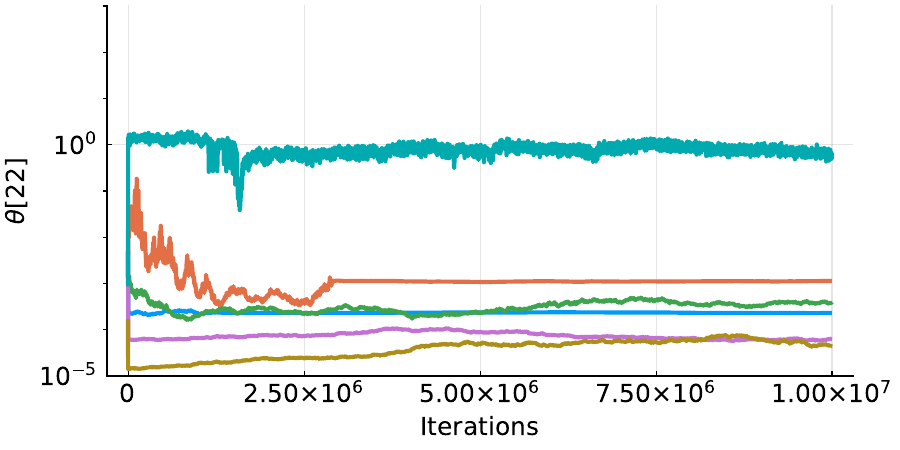} 
    		\caption{Approach from \cite{calvetti2024computationally}}
    		\label{fig:deblurring_model4_traces_theta22_Calvetti2024_AM_prior}
  	\end{subfigure}%
	~	 
	\begin{subfigure}[b]{0.32\textwidth}
		\includegraphics[width=\textwidth]{%
      		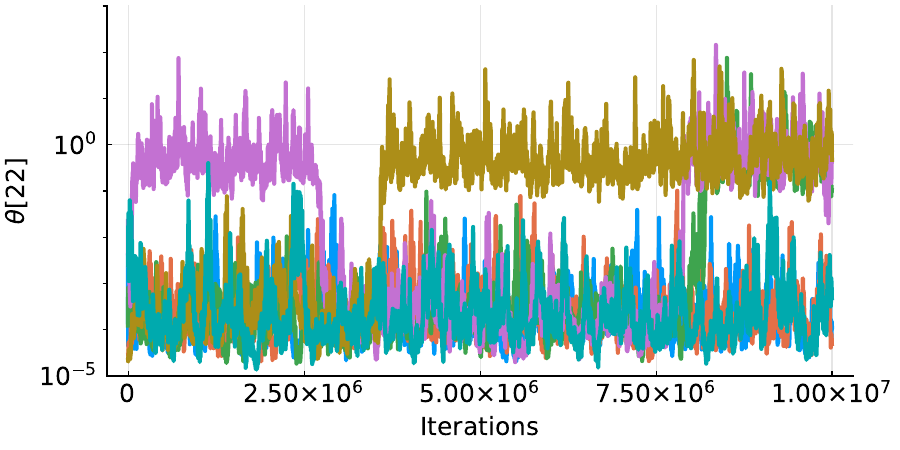} 
    		\caption{Prior-normalized}
    		\label{fig:deblurring_model4_traces_theta22_priorNormalized_AM_prior}
  	\end{subfigure}%
  	\caption{ 
		Traces for the $x_{22}$- and $\theta_{22}$-samples and $r=-1$. 
		Note that $t_{22} \approx 0.172$ is slightly to the right of the first jump discontinuity at $t = 0.17$ and yields $g(t_{22}) = 1$.
		We used the AM algorithm initialized with random draws from the respective prior. 
	}
  	\label{fig:deblurring_traces_22_AM_prior}
\end{figure}

\cref{fig:deblurring_samples_AM_prior} illustrates four individual $x$-samples from the first chain for the three posteriors with $r=-1$. 
We observe in \cref{fig:deblurring_model4_xSamples_original_AM_prior} that the individual samples from the original posterior qualitatively do not change significantly. 
This confirms that the AM method only locally explores the original posterior in a high-density region near the initialization. 
In contrast, \cref{fig:deblurring_model4_xSamples_Calvetti2024_AM_prior,fig:deblurring_model4_xSamples_priorNormalized_AM_prior} show that the samples for the transformed posterior from \cite{calvetti2024computationally} and our prior-normalized posterior vary more visibly. 
Yet, we observe that the samples from our prior-normalized posterior better capture the true underlying signal than those from the transformed posterior in \cite{calvetti2024computationally}.
The above findings reaffirm that our prior-normalized posterior is explored more effectively than the other two posteriors.

Next, \cref{fig:deblurring_traces_21_AM_prior,fig:deblurring_traces_22_AM_prior} display the trace plots for the $x_{i}$- and $\theta_i$-samples for $i=21$ (\cref{fig:deblurring_traces_21_AM_prior}), $i=22$ (\cref{fig:deblurring_traces_22_AM_prior}), and $r = -1$. 
The results in \cref{fig:deblurring_traces_21_AM_prior,fig:deblurring_traces_22_AM_prior} provide further evidence that the AM sampler explores the original posterior and the transformed posterior from \cite{calvetti2024computationally} only locally, while exploring the prior-normalized posterior more effectively. 
Specifically, in \cref{fig:deblurring_model4_traces_x22_priorNormalized_AM_prior}, the chains for the prior-normalized posterior again better explore the two distinct high-density regions.

Next, \cref{tab:signal_MPSRF} reports the computational wall times in seconds (``time"), the MPSRF minus one multiplied by the time (``MPSRF"), and the ESS per second (``ESS'') for the original posterior, the transformed posterior from \cite{calvetti2024computationally}, and our prior-normalized posterior.
For $r=1$, corresponding to a log-concave original posterior, the MPSRF and ESS values for the original and our prior-normalized posterior are comparable: the prior-normalized posterior exhibits a slightly lower MPSRF, while the original posterior achieves a higher ESS. 
At the same time, the transformed posterior from \cite{calvetti2024computationally} has an ESS between the two other posteriors, but a MPSRF of roughly two magnitudes higher. 
For $r = -1$, although the prior-normalized posterior again yields a lower MPSRF compared to the other two posteriors---indicating improved convergence speed---it also results in a slightly lower ESS compared to the original posterior.
More critically, despite the sampler failing to adequately explore the original posterior and the transformed posterior from \cite{calvetti2024computationally}---becoming trapped in local high-density regions, as illustrated in \cref{fig:deblurring_chains_AM_prior}---this deficiency is not reflected in the reported MPSRF and ESS values. 

\begin{table}[h]
	\renewcommand{\arraystretch}{1.1}
	\centering
	\begin{tabular}{ l | c c c | c c c }
		\toprule 
		 & \multicolumn{3}{c}{$r=1$} &
		 \multicolumn{3}{c}{$r=-1$} \\
		posterior & 
		time & MPSRF & ESS & 
		time & MPSRF & ESS \\
		\midrule
    		original & 
		5.4e+01 & 5.8e+01 & 2.5e+01 &  
		7.5e+01 & 5.4e+06 & 3.1e+00 \\
		approach from \cite{calvetti2024computationally} & 
		1.0e+02 & 5.2e+03 & 1.3e+01 &  
		9.3e+01 & 1.1e+07 & 1.3e+00 \\
   		prior-normalized & 
		3.7e+02 & 5.5e+01 & 7.7e+00 &  
		3.5e+02 & 6.8e+03 & 1.6e+00 \\
		\bottomrule
	\end{tabular}
	\caption{
		The computational wall time in seconds (``time"), the MPSRF minus one multiplied by the time (``MPSRF"), and the ESS per second (``ESS'') for the signal deconvolution problem's original and prior-normalized posterior. 
		We use the AM algorithm initialized with random prior draws.
	}
  	\label{tab:signal_MPSRF}
\end{table}

The above observations highlight a key limitation: \emph{neither MPSRF nor ESS reliably reflects the quality of posterior exploration in multi-modal settings}.
In fact, many MCMC diagnostics fail in the presence of multi-modality or heavy tails; among them, MPSRF is arguably ``not too bad''---provided that one can initialize the chains from sufficiently over-dispersed starting points.
Nonetheless, as with any diagnostic, these measures can only certify non-convergence: they raise a flag when sampling behaves poorly, but the absence of such a flag does not guarantee that the chains have mixed well.
For this reason, we refrain from reporting MPSRF and ESS values in the subsequent experiments.

\begin{figure}[tb]
	\centering 
	\begin{subfigure}[b]{0.32\textwidth}
		\includegraphics[width=\textwidth]{%
      		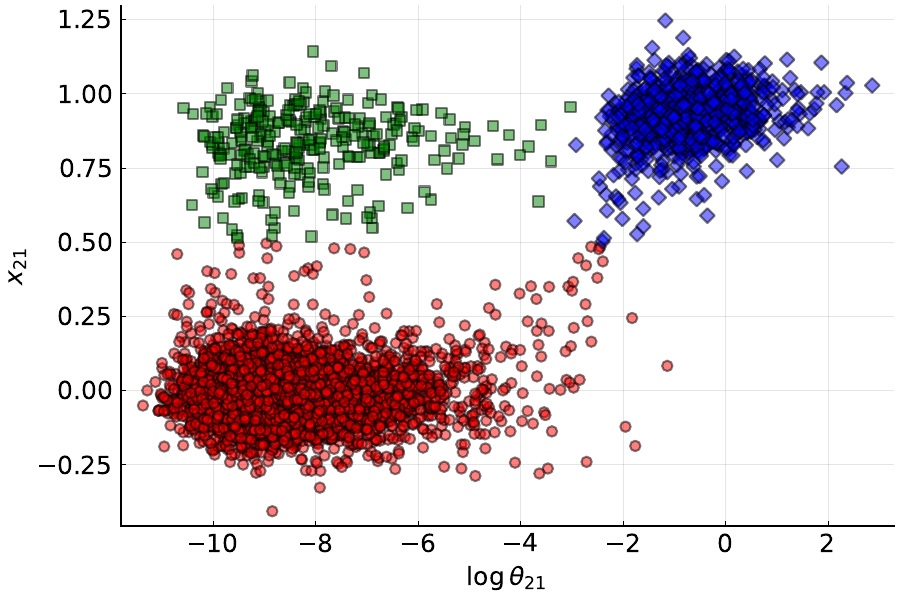} 
    		\caption{Original}
    		\label{fig:deblurring_model4_scatter_x_theta_21_AM_prior}
  	\end{subfigure}%
	~	 
	\begin{subfigure}[b]{0.32\textwidth}
		\includegraphics[width=\textwidth]{%
      		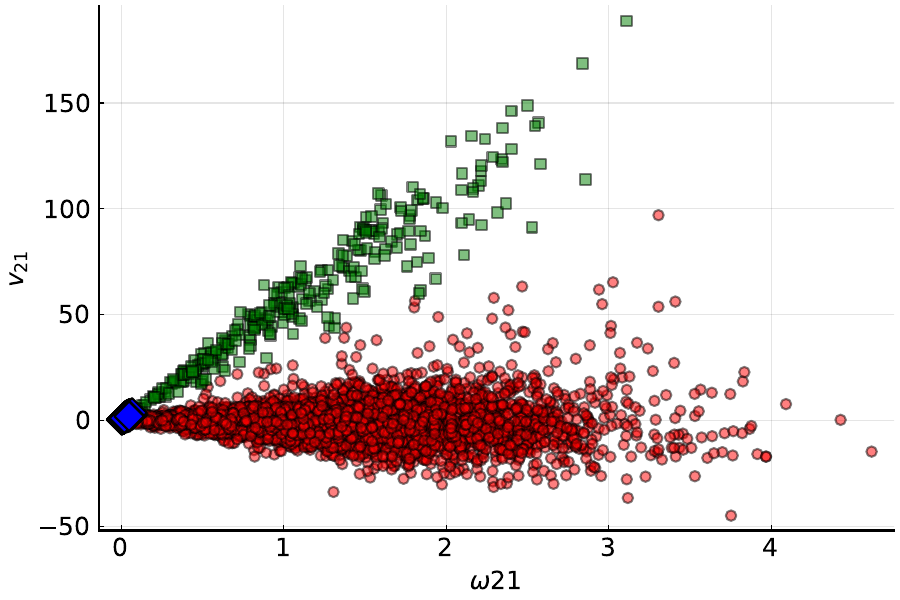} 
    		\caption{Approach from \cite{calvetti2024computationally}}
    		\label{fig:deblurring_model4_scatter_v_omega_21_AM_prior}
  	\end{subfigure}%
	~	 
	\begin{subfigure}[b]{0.32\textwidth}
		\includegraphics[width=\textwidth]{%
      		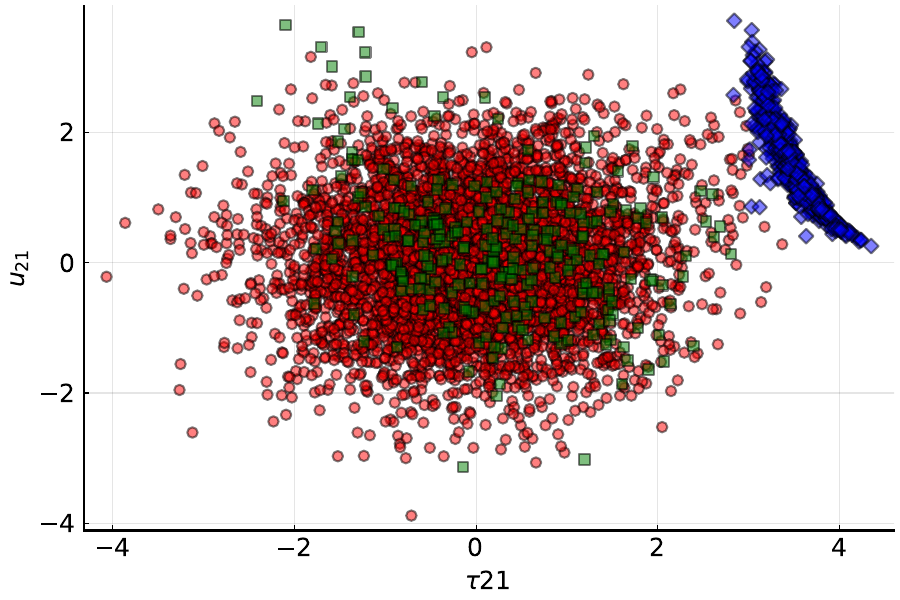} 
    		\caption{Prior-normalized}
    		\label{fig:deblurring_model4_scatter_u_tau_21_AM_prior}
  	\end{subfigure}%
	\\
	\begin{subfigure}[b]{0.32\textwidth}
		\includegraphics[width=\textwidth]{%
      		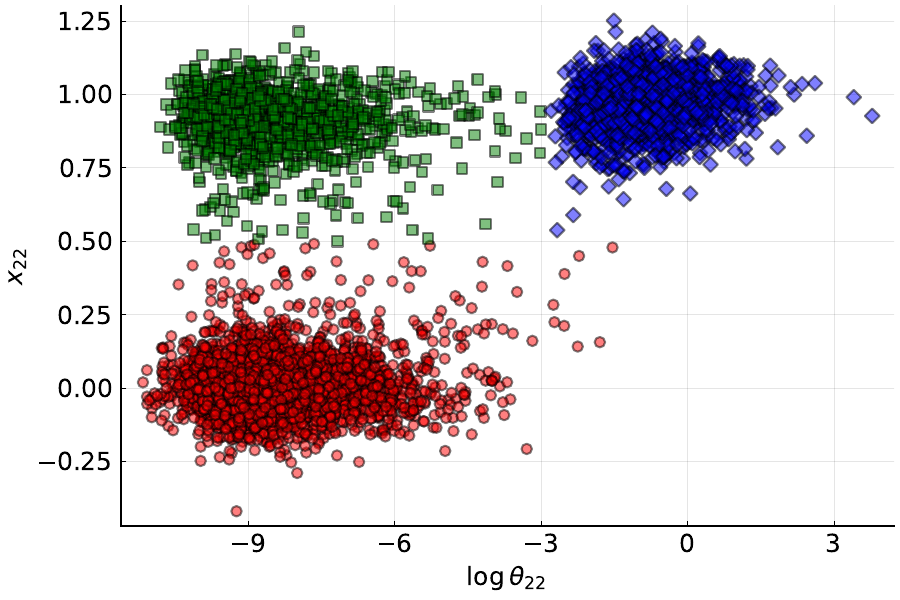} 
    		\caption{Original}
    		\label{fig:deblurring_model4_scatter_x_theta_22_AM_prior}
  	\end{subfigure}%
	~	 
	\begin{subfigure}[b]{0.32\textwidth}
		\includegraphics[width=\textwidth]{%
      		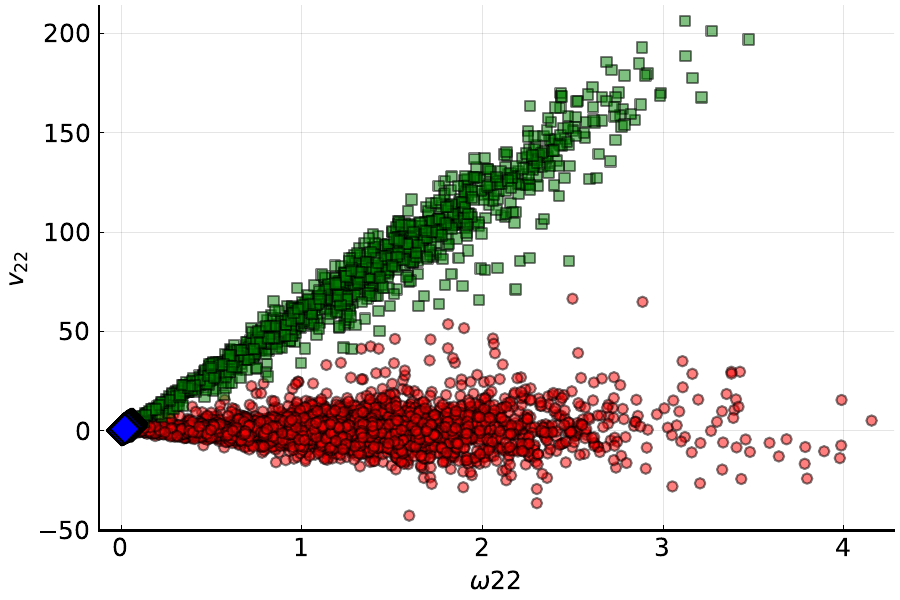} 
    		\caption{Approach from \cite{calvetti2024computationally}}
    		\label{fig:deblurring_model4_scatter_v_omega_22_AM_prior}
  	\end{subfigure}%
	~	 
	\begin{subfigure}[b]{0.32\textwidth}
		\includegraphics[width=\textwidth]{%
      		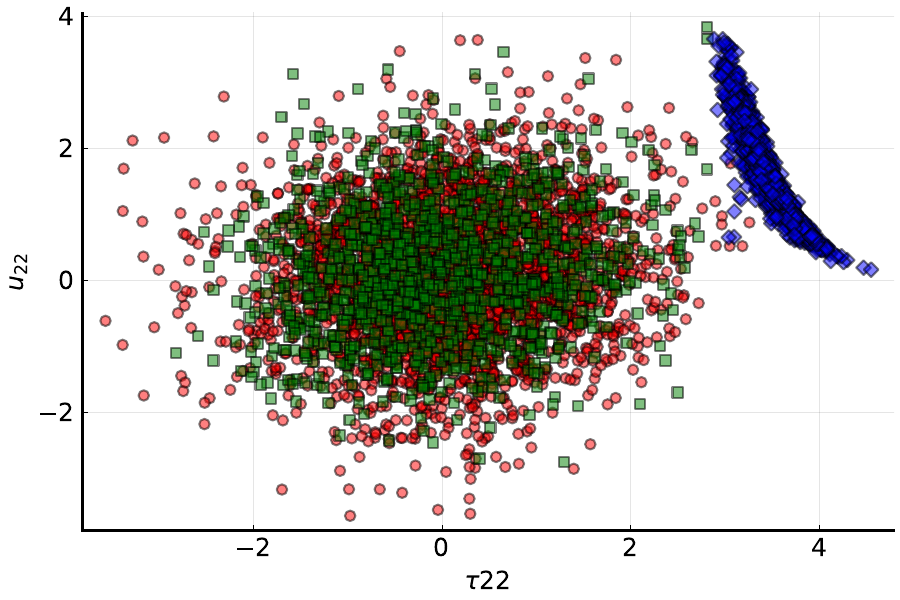} 
    		\caption{Prior-normalized}
    		\label{fig:deblurring_model4_scatter_u_tau_22_AM_prior}
  	\end{subfigure}%
  	\caption{ 
		Scatter plots for the $(x_{i},\theta_{i})$-, $(v_{i},\omega_{i})$-, and $(u_{i},\tau_{i})$-samples for $i=21$ (top), $i=22$ (bottom), and $r=-1$. 
		The $(u_{i},\tau_{i})$-samples were generated from the prior-normalized posterior using the AM algorithm. 
		These were then mapped into the $(\mathbf{x},\boldsymbol{\theta})$- and $(\mathbf{v},\boldsymbol{\tau})$-coordinate system to generate the $(x_{i},\theta_{i})$- and $(u_{i},\tau_{i})$-samples. 
	}
  	\label{fig:deblurring_scatter_AM_prior}
\end{figure}

We conclude this subsection by examining in more detail the geometry of the various posteriors arising in the signal deblurring problem.
To this end, \cref{fig:deblurring_scatter_AM_prior} displays scatter plots of the $(x_i,\theta_i)$-, $(v_i,\omega_i)$-, and $(u_i,\tau_i)$-samples for indices $i=21$ and $i=22$ with $r=-1$.
The $(u_i,\tau_i)$-samples were obtained by exploring our prior-normalized posterior using the AM algorithm; these samples were subsequently mapped into the $(\mathbf{x},\boldsymbol{\theta})$- and $(\mathbf{v},\boldsymbol{\omega})$-coordinates to produce the corresponding $(x_i,\theta_i)$- and $(v_i,\omega_i)$-samples.
Importantly, neither the $(x_i,\theta_i)$-samples for the original posterior nor the $(v_i,\omega_i)$-samples for the transformation of \cite{calvetti2024computationally} were generated by directly sampling their respective posteriors, since---as demonstrated extensively above---these distributions cannot be efficiently explored using the AM sampler.
The scatter plots in \cref{fig:deblurring_model4_scatter_x_theta_21_AM_prior,fig:deblurring_model4_scatter_x_theta_22_AM_prior} reveal that the original posterior exhibits three dominant modes, visible as three distinct sample clusters.
We highlight these clusters as follows:
\begin{enumerate}
	\item[(i)] red circles for samples with $x_i \le 0.5$;

	\item[(ii)] green squares for samples with $x_i > 0.5$ and $\log \theta_i \le -3$;

	\item[(iii)] blue diamonds for samples with $x_i > 0.5$ and $\log \theta_i > -3$.
\end{enumerate}
For the original posterior, cluster (iii) is clearly separated from clusters (i) and (ii), reflecting a difference of roughly seven orders of magnitude in their corresponding $\theta$-values.
In contrast, for the prior-normalized posterior shown in \cref{fig:deblurring_model4_scatter_u_tau_21_AM_prior,fig:deblurring_model4_scatter_u_tau_22_AM_prior}, clusters (i) and (ii) appear to merge and resemble draws from the normalized prior distribution.
Cluster (iii), on the other hand, corresponds to samples from a curved, high-density region that lies close to the other two clusters.
This indicates that, although the prior-normalized posterior may still be multimodal, the modes are substantially closer together, which in turn facilitates movement between them for the AM sampler.
For the transformed posterior of \cite{calvetti2024computationally}, shown in \cref{fig:deblurring_model4_scatter_v_omega_21_AM_prior,fig:deblurring_model4_scatter_v_omega_22_AM_prior}, clusters (i) and (ii) remain well separated and lie along two distinct linear stripes that appear to intersect at $(0,0)$.
Furthermore, cluster (iii) has effectively collapsed into a tightly concentrated region near the origin.
This geometry is expected to pose considerable challenges for sampling.
Indeed, comparing \cref{fig:deblurring_model4_scatter_v_omega_21_AM_prior,fig:deblurring_model4_scatter_v_omega_22_AM_prior} with the scatter plot for $j=50$ in \cite[Figure 8]{calvetti2024computationally} suggests that the modified pCN algorithm employed in \cite{calvetti2024computationally} explored only a portion of the lower stripe corresponding to the red-circle cluster in our visualization.

\subsection{Recovering the initial data for Burgers' equation}
\label{sub:tests_Burgers}

We consider the initial value problem (IVP) for the non-linear inviscid Burgers equation in one spatial dimension: 
\begin{equation}\label{eq:Burgers}
	\partial_t u(x,t) + \partial_x u(x,t)^2 = 0, \quad t>0, \ x \in \Omega = (0,1),
\end{equation}  
with periodic boundary conditions and initial condition $u(x,0) = u_0(x)$. 
Our goal is to estimate the nodal values of the piecewise constant initial data 
\begin{equation}\label{eq:initial_data}
	u_0(x) = 
	\begin{cases} 
		1 & \text{ if } |x-1/2| < 0.25, \\ 
		0 & \text{ otherwise},
	\end{cases}
\end{equation} 
from noisy and undersampled observations of the solution to \cref{eq:Burgers} at $t=0.25$. 
To set up the forward operator, we spatially discretize \cref{eq:Burgers} via an upwind finite difference scheme (see \cite{ranocha2025robustness}) on $n=100$ equidistant points, implemented in Trixi.jl \cite{ranocha2022adaptive}. 
Time integration is carried out with the third-order, four-stage strong stability preserving (SSP) Runge--Kutta (RK) method ``SSPRK43'' from the OrdinaryDiffEq.jl package \cite{rackauckas2017differentialequations}. 
We then apply undersampling by observing the numerical solution only at every fifth grid point, yielding $m=20$ noisy observations.
The resulting data model is 
\begin{equation}\label{eq:Burgers_model}
	\mathbf{b} = \mathcal{F}(\mathbf{u}) + \mathbf{e}, 
\end{equation} 
where $\mathbf{u} = [u_1, \dots, u_N]$ are the nodal values of the initial data and $\mathcal{F}: \R^N \to \R^M$ is the nonlinear forward operator corresponding to numerically solving \cref{eq:Burgers} with initial data $\mathbf{u}$, followed by undersampling.
To avoid the ``inverse crime" \cite{kaipio2007statistical}, we generate the observational data from the exact solution of \cref{eq:Burgers} at time $t=0.25$. 
We then add i.i.d.\ Gaussian noise $e_j \sim \mathcal{N}(0, \sigma^2)$ with $\sigma = 10^{-2}$. 
\Cref{fig:deblurring_data} illustrates the problem setup.

\begin{figure}[tb]
	\centering 
	\includegraphics[width=0.6\textwidth]{%
      		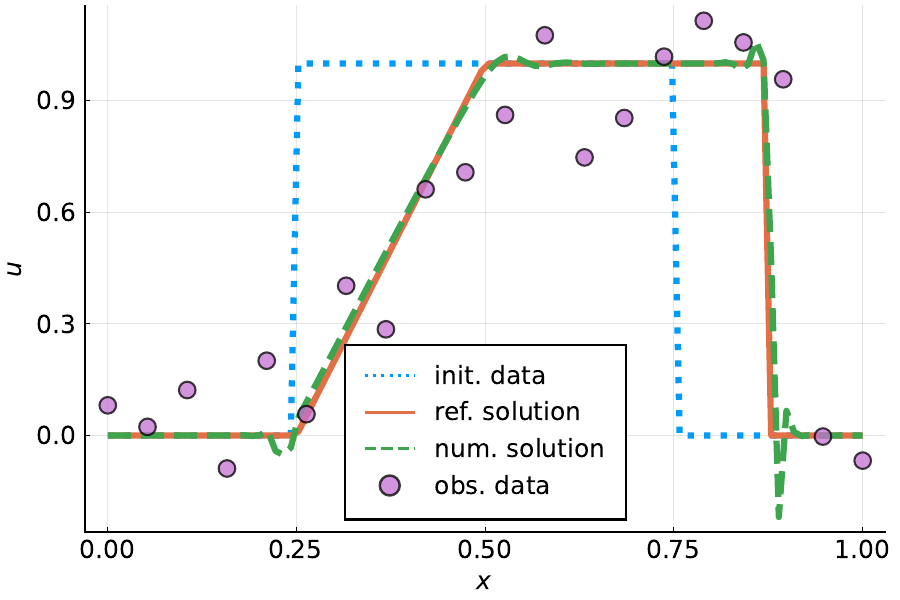} 
  	\caption{ 
		Setup for the inverse Burgers problem, including the initial data (dotted blue), the reference (solid orange) and numerical (dashed green) solutions, and the observational data (purple dots). 
		The observational data is generated from the reference solution, whereas the forward operator uses the numerical solution.
  	}
  	\label{fig:Burgers_data}
\end{figure} 

To sparsely represent the initial data $\mathbf{u}$, we again let $L \in \mathbb{R}^{n \times n}$ be the finite difference matrix \cref{eq:FD_matrix} and write $\mathbf{u}$ in terms of its increments $\mathbf{v}$ as $\mathbf{u} = L^{-1} \mathbf{v}$.
Substituting into \cref{eq:Burgers_model} yields $\mathbf{b} = \mathcal{F}( L^{-1} \mathbf{v} ) + \mathbf{e}$. 
We model $\mathbf{u}$ being piecewise constant by assuming $\mathbf{v}$ is sparse.

\begin{figure}[tb]
	\centering 
	\begin{subfigure}[b]{0.45\textwidth}
		\includegraphics[width=\textwidth]{%
      		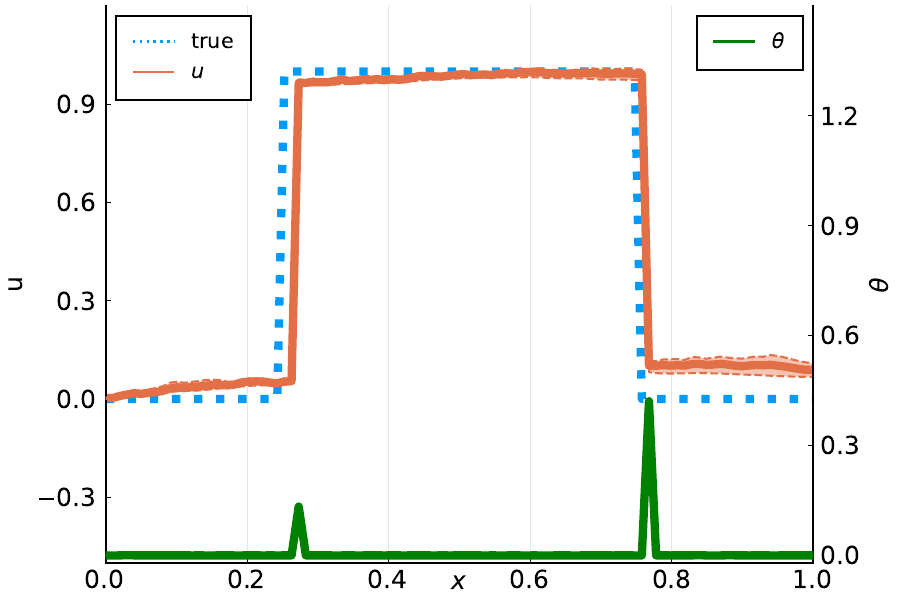} 
    		\caption{Original, MAP initialized}
    		\label{fig:Burgers_model4_UQquantile_original_AM_MAP}
  	\end{subfigure}%
	~
	\begin{subfigure}[b]{0.45\textwidth}
		\includegraphics[width=\textwidth]{%
      		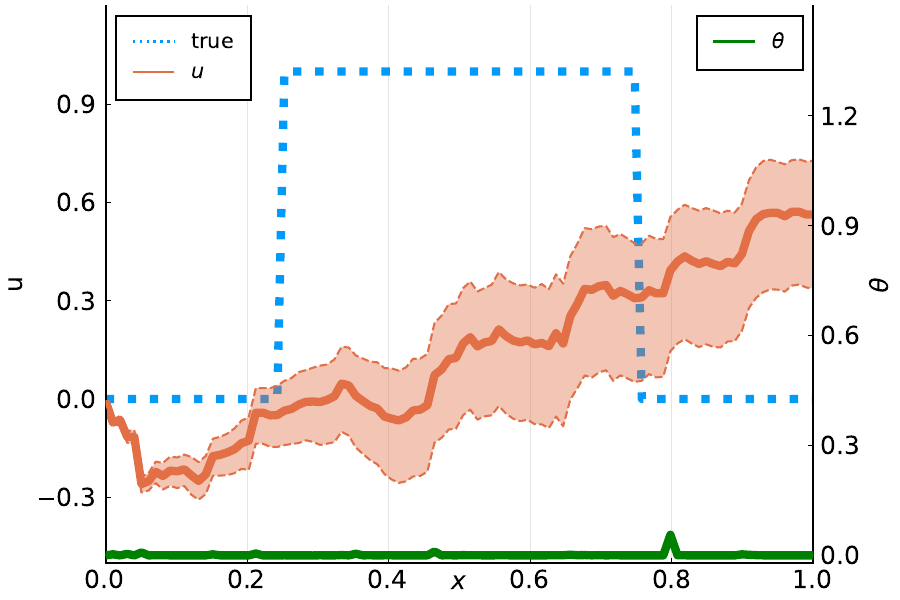} 
    		\caption{Original, rand.\ initialized}
    		\label{fig:Burgers_model4_UQquantile_original_AM_prior}
  	\end{subfigure}%
	\\ 
	\begin{subfigure}[b]{0.45\textwidth}
		\includegraphics[width=\textwidth]{%
      		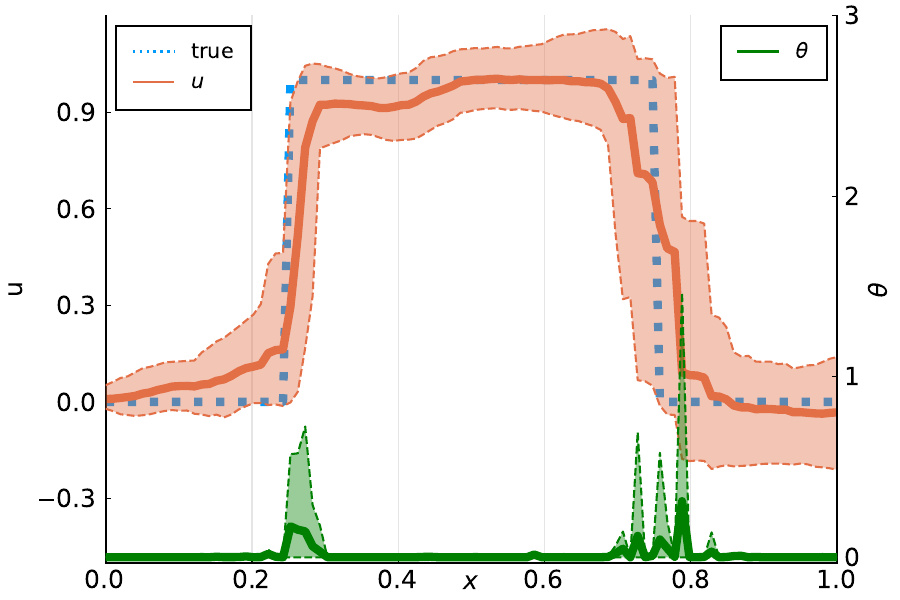} 
    		\caption{Prior-normalized, MAP initialized}
    		\label{fig:Burgers_model4_UQquantile_priorNormalized_AM_MAP}
  	\end{subfigure}%
	~
	\begin{subfigure}[b]{0.45\textwidth}
		\includegraphics[width=\textwidth]{%
      		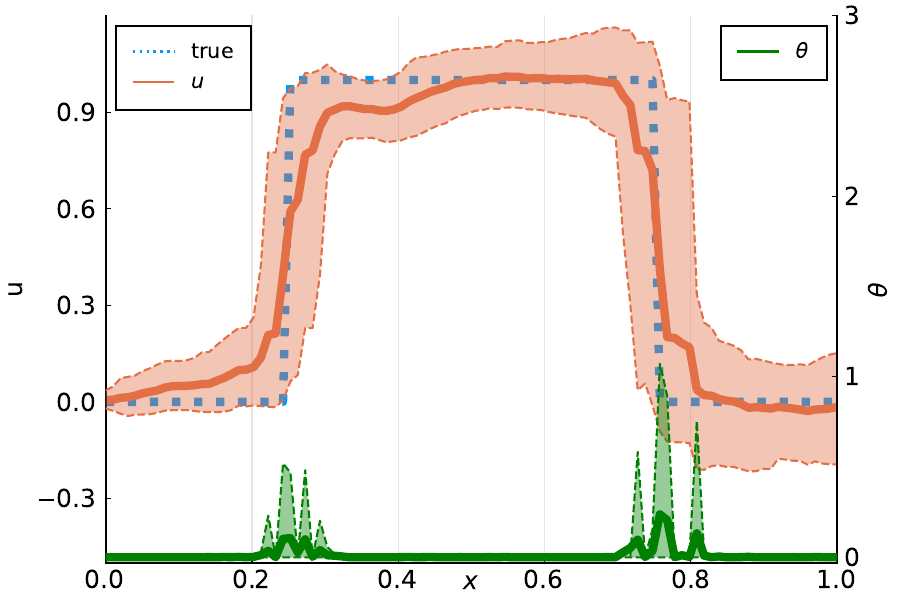} 
    		\caption{Prior-normalized, rand.\ initialized}
    		\label{fig:Burgers_model4_UQquantile_priorNormalized_AM_prior}
  	\end{subfigure}%
  	\caption{ 
		Mean and $90\%$ quantile ranges of the original and prior-normalized posterior with $r=-1$.
		We initialized the AM algorithm with the respective MAP estimate and random prior draws.     
  	}
  	\label{fig:Burgers_UQ_AM}
\end{figure} 

We assess the performance of the AM algorithm of \cite{haario2001adaptive,andrieu2008tutorial,atchade2010limit} with a target mean acceptance rate of 23.4\% on the original and prior-normalized posteriors for $r = -1$ and $\beta, \vartheta$ as in \cref{tab:parameters}.
We generate a single chain with $10^6$ samples, saving every $10^2$th sample. 
\Cref{fig:Burgers_UQ_AM} presents the mean and $90\%$ quantile ranges of the 
$u$- and $\theta$-samples. 
We initialize the chains with the respective MAP estimate (in \cref{fig:Burgers_model4_UQquantile_original_AM_MAP,fig:Burgers_model4_UQquantile_priorNormalized_AM_MAP}) and random prior draws (\cref{fig:Burgers_model4_UQquantile_original_AM_prior,fig:Burgers_model4_UQquantile_priorNormalized_AM_prior}).
We observe for the original posterior in \cref{fig:Burgers_model4_UQquantile_original_AM_MAP,fig:Burgers_model4_UQquantile_original_AM_prior} that the sampler's performance again depends heavily on how it is initialized.
Specifically, the sample means deviate only marginally from the MAP estimate if the latter is used to initialize the chain. 
In contrast, the $\mathbf{x}$-mean deviates significantly from the true initial data when the chain is initialized randomly. 
This observation suggests that the AM sampler struggles to adequately explore the original posterior. 
In contrast, for the prior-normalized posterior in \cref{fig:Burgers_model4_UQquantile_priorNormalized_AM_MAP,fig:Burgers_model4_UQquantile_priorNormalized_AM_prior}, the mean and quantile ranges agree more closely with the underlying initial data, independently of how they are initialized, indicating that the AM sampler explores the prior-normalized posterior more effectively.

\begin{figure}[tb]
	\centering 
	\begin{subfigure}[b]{0.45\textwidth}
		\includegraphics[width=\textwidth]{%
      		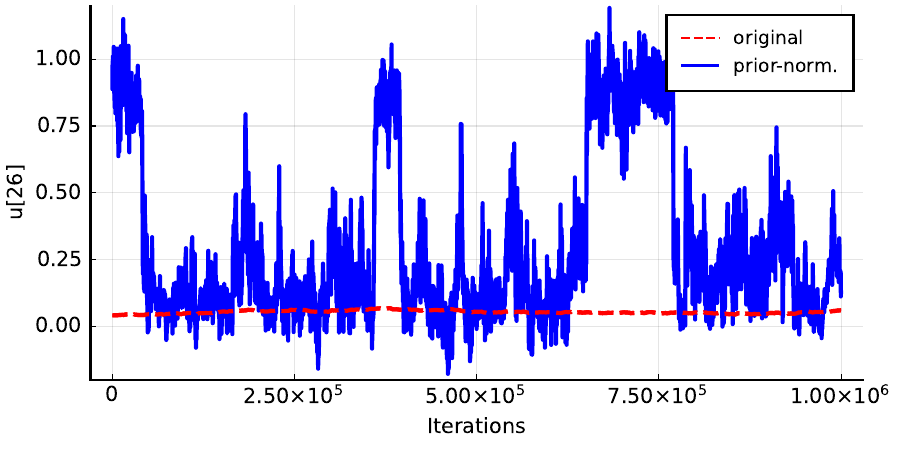} 
    		\caption{MAP initialized}
    		\label{fig:Burgers_model4_traces_u26_AM_MAP}
  	\end{subfigure}%
	~
	\begin{subfigure}[b]{0.45\textwidth}
		\includegraphics[width=\textwidth]{%
      		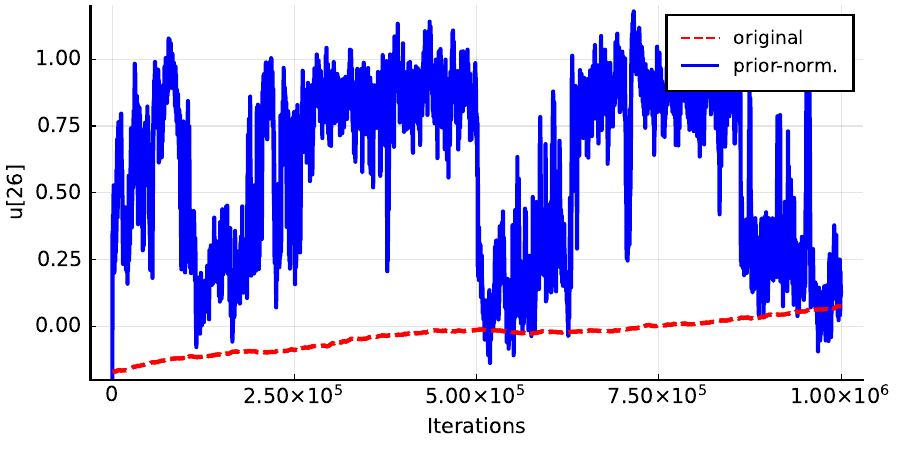} 
    		\caption{Randomly initialized}
    		\label{fig:Burgers_model4_traces_u26_AM_prior}
  	\end{subfigure}%
  	\caption{ 
		Traces for the $u_{26}$-samples of the original and prior-normalized posterior for $r=-1$. 
		Note that $x_{26} \approx 0.252$ is close to the first jump discontinuity and yields $u_0(x_{26}) = 1$. 
		We used the AM algorithm initialized with the respective MAP estimate and random draws from the prior. 
  	}
  	\label{fig:Burgers_traces_AM}
\end{figure} 

Moreover, \cref{fig:Burgers_traces_AM} displays the trace plots for the $u_{26}$-samples, using $r = -1$ and prior samples to initialize the AM sampler. 
Note that $x_{26} \approx 0.252$ is close to the first jump discontinuity and yields $u_0(x_{26}) = 1$.
\cref{fig:Burgers_traces_AM} provides further evidence that the AM sampler explores the original posterior only locally. 
In contrast, the same sampler explores two distinct high-density regions for the prior-normalized posterior. 
As before, this outcome is intuitively desirable, as $x_{26} \approx 0.252$ is slightly to the right of the first jump discontinuity. 
While the true initial data at this grid point is $u_0(t_{22}) = 1$, the Bayesian model should reflect the uncertainty that the jump might occur at the previous grid point, meaning both $u_{26} = 1$ and $u_{26} = 0$ should lie in high-density posterior regions.

\subsection{Impulse image} 
\label{sub:tests_impulse} 

\begin{figure}[tb]
	\centering 
	\begin{subfigure}[b]{0.32\textwidth}
		\includegraphics[width=\textwidth]{%
      		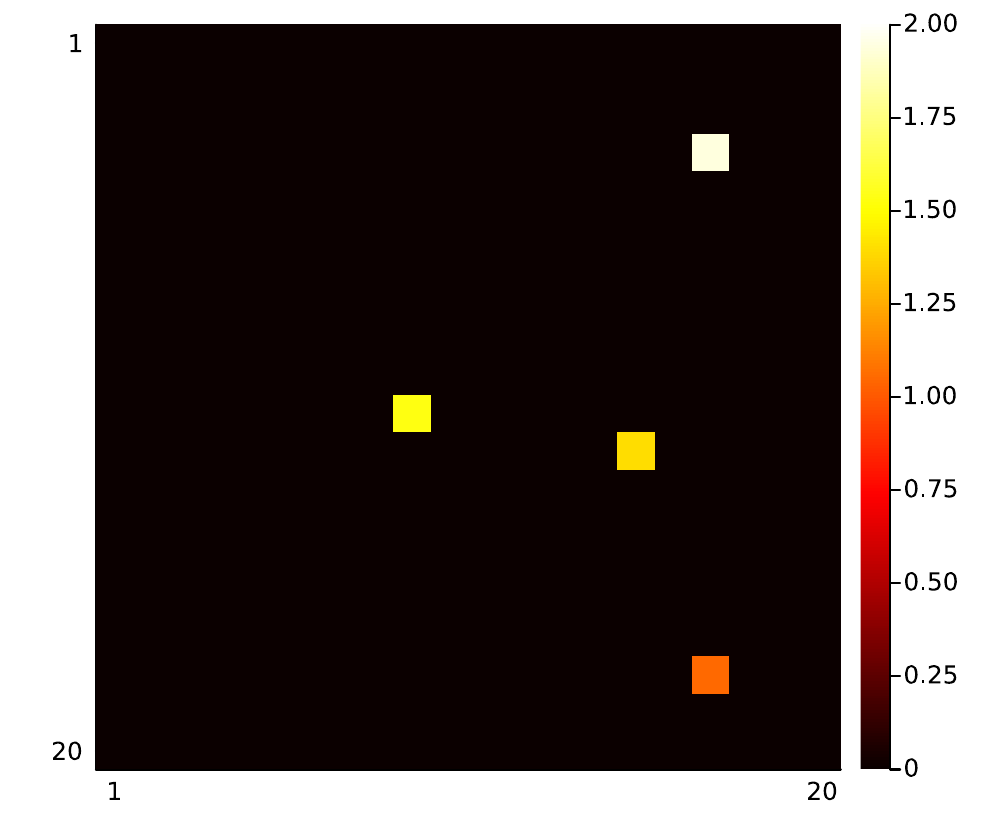} 
    		\caption{True signal}
    		\label{fig:impulse_n20_imageRef}
  	\end{subfigure}%
	\begin{subfigure}[b]{0.32\textwidth}
		\includegraphics[width=\textwidth]{%
      		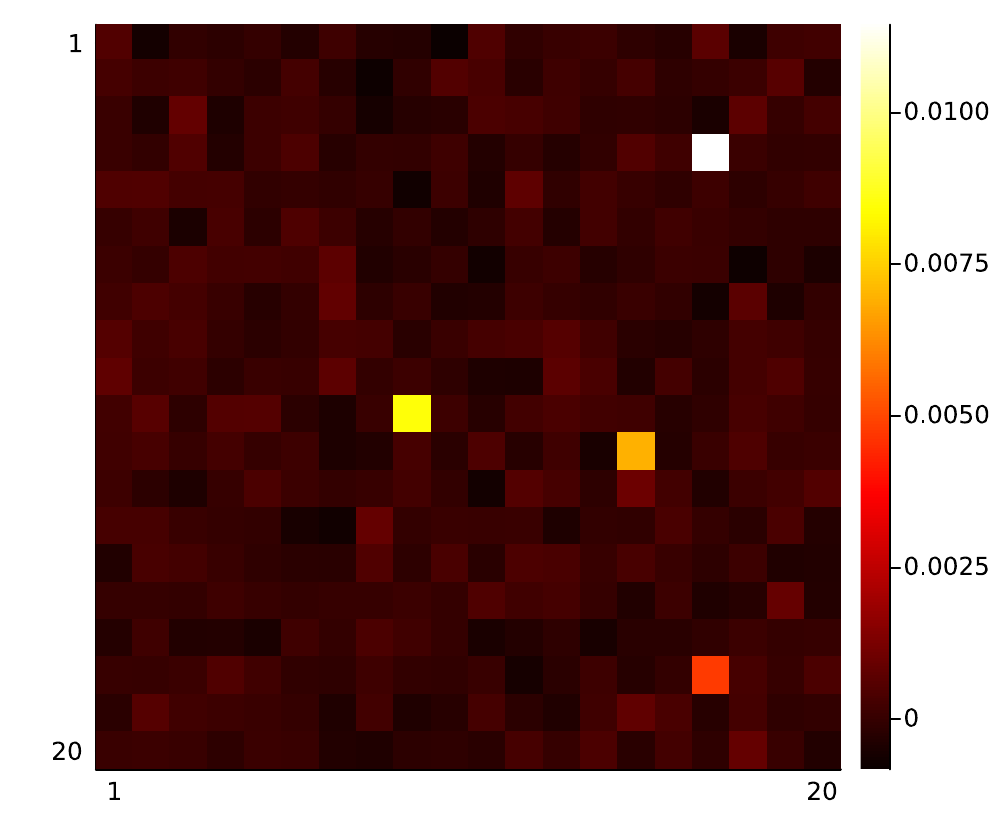} 
    		\caption{Original, $\mathbf{x}$-mean}
    		\label{fig:impulse_n20_model4_x_mean_Gibbs_original}
  	\end{subfigure}%
	\begin{subfigure}[b]{0.32\textwidth}
		\includegraphics[width=\textwidth]{%
      		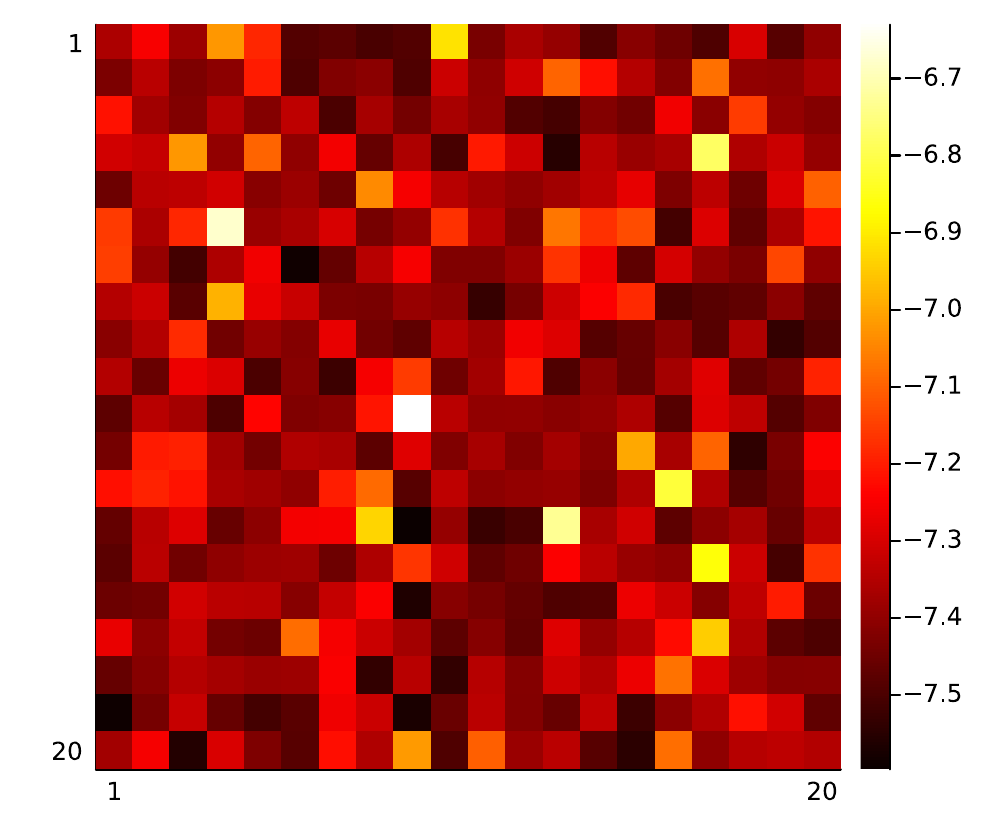} 
    		\caption{Original, log-$\boldsymbol{\theta}$-mean}
    		\label{fig:impulse_n20_model4_theta_mean_Gibbs_original}
  	\end{subfigure}%
	\\ 
	\begin{subfigure}[b]{0.32\textwidth}
		\includegraphics[width=\textwidth]{%
      		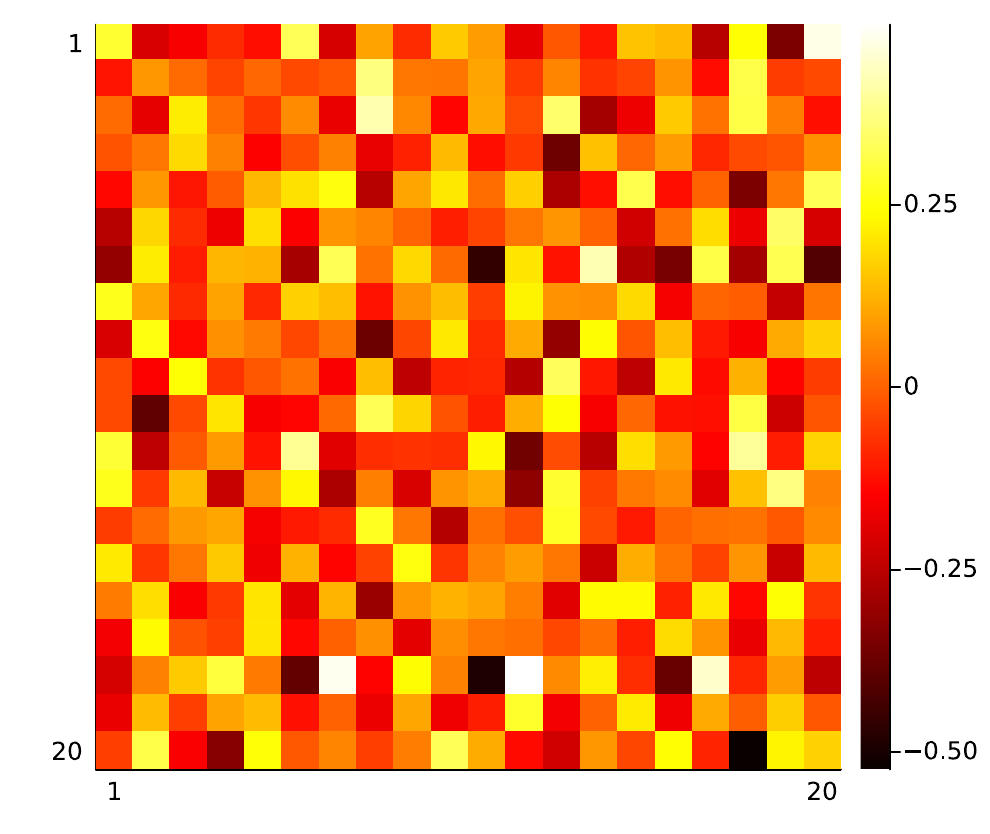} 
    		\caption{Noisy DCT data}
    		\label{fig:impulse_n20_obsData}
  	\end{subfigure}%
	\begin{subfigure}[b]{0.32\textwidth}
		\includegraphics[width=\textwidth]{%
      		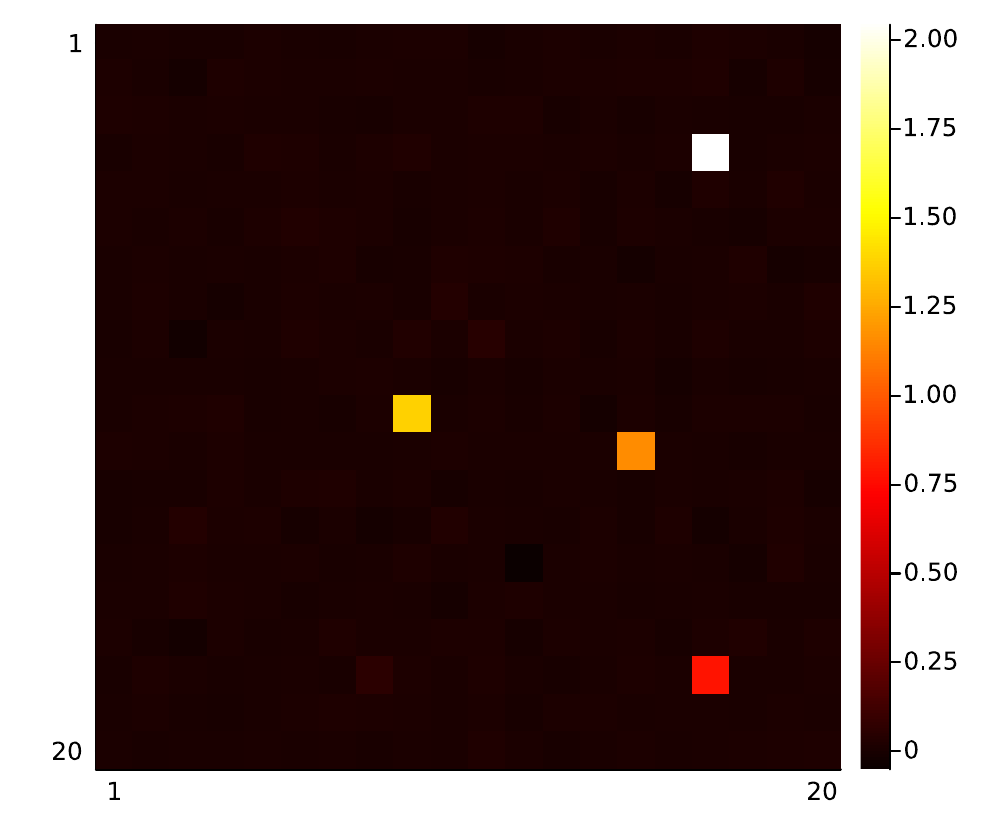} 
    		\caption{Prior-normalized, $\mathbf{x}$-mean}
    		\label{fig:impulse_n20_model4_x_mean_ESS_priorNormalized}
  	\end{subfigure}%
	\begin{subfigure}[b]{0.32\textwidth}
		\includegraphics[width=\textwidth]{%
      		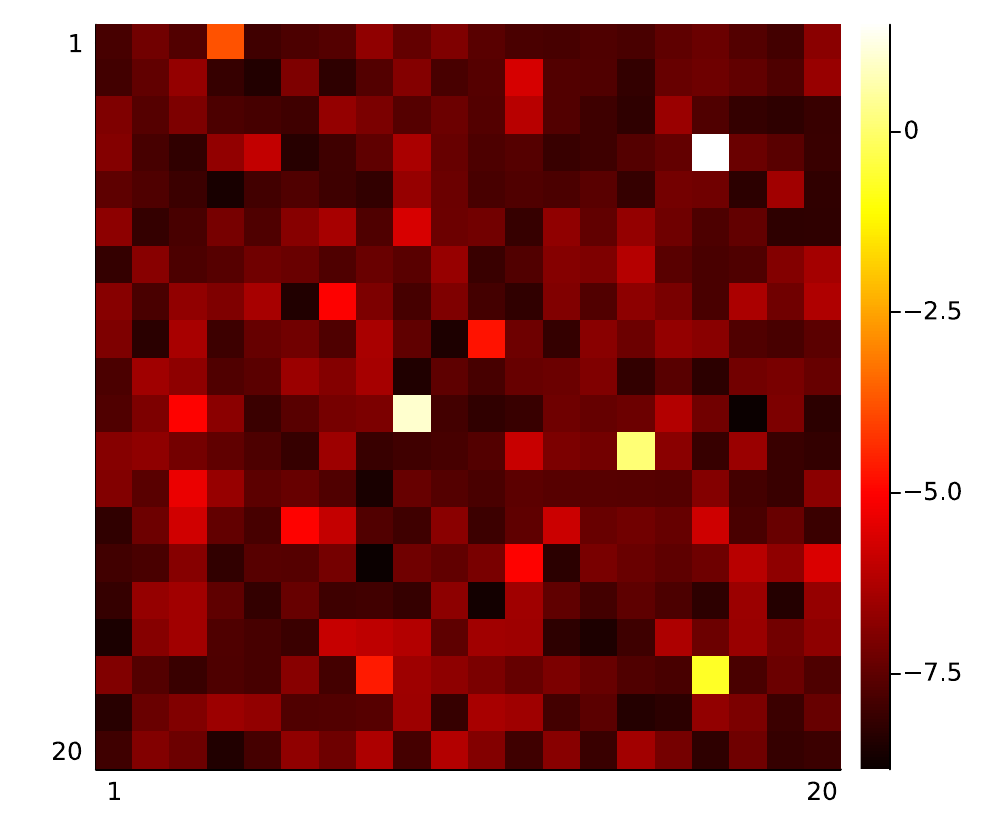} 
    		\caption{Prior-normalized, log-$\boldsymbol{\theta}$-mean}
    		\label{fig:impulse_n20_model4_theta_mean_ESS_priorNormalized}
  	\end{subfigure}%
  	\caption{ 
		(Log-)Means of the $x$- and $\theta$-samples from the original and prior-normalized posterior corresponding to the impulse image.
		The Gibbs sampler generated four chains with $10^5$ samples each from the original posterior ($\approx 20$ minutes runtime).
		The ES sampler generated four chains with $10^4$ samples each from the original posterior ($\approx 3$ minutes runtime).
		Both samplers were initialized using random prior draws.     
  	}
  	\label{fig:impulse}
\end{figure}

We consider the problem of inferring the $20 \times 20$ impulse image in \cref{fig:impulse_n20_imageRef} from its noisy DCT data in \cref{fig:impulse_n20_obsData}. 
The resulting data model can be formulated as  
\begin{equation}\label{eq:impulse_model}
	\mathbf{y} = F \mathbf{x} + \mathbf{e}, 
\end{equation} 
where $\mathbf{x} \in \R^{n}$ with $n = n_1 n_2$ are the pixel intensities of the vectorized image $X \in \R^{n_1 \times n_2}$, $F \in \R^{n \times n}$ is the matrix representation of the DCT, and $\mathbf{y} \in \R^N$ is the noisy  DCT data. 
Here, $n_1 = n_2 = 20$, $n = 400$, and $\mathbf{e} \sim \mathcal{N}(0, \sigma^2 I)$ with $\sigma = 10^{-2}$. 
We model the impulse image having only a few non-zero entries by promoting $\mathbf{x}$ to be sparse and assume that $\boldsymbol{\theta}$ follows an inverse gamma distribution, i.e., $\theta_i \sim \mathcal{GG}(r,\beta,\vartheta)$ with $r = -1$ and $\beta, \vartheta$ as in \cref{tab:parameters}. 
These specific choices yield the following conditional distributions for the original posterior: 
\begin{subequations}\label{eq:impulse_Gibbs}
\begin{align}
	\mathbf{x} | \mathbf{y}, \boldsymbol{\theta} 
		& \sim \mathcal{N}(\boldsymbol{\mu}, \Sigma), \quad 
		\Sigma^{-1} = \frac{1}{\sigma^2} F^T F + \diag( \boldsymbol{\theta} )^{-1}, \quad 
		\boldsymbol{\mu} = \frac{1}{\sigma^2} \Sigma F^T \mathbf{y}, 
		\label{eq:impulse_Gibbs1} \\ 
	\theta_i | \mathbf{y}, \mathbf{x} 
		& \sim \mathcal{GG}(r, \tilde{\beta}, \tilde{\vartheta}_i), \quad 
		\tilde{\beta} = \beta + \frac{1}{2}, \quad
		\tilde{\vartheta}_i = \vartheta + \frac{x_i^2}{2}, \quad
		i=1,\dots,n 
		\label{eq:impulse_Gibbs2}
\end{align}
\end{subequations}   
The conditional distributions \cref{eq:impulse_Gibbs} being standard ones allows us to use a Gibbs sampler \cite{agapiou2014analysis,ascolani2024scalability}, which is popular for hierarchical Bayesian models as it comes with relatively small per-sample costs \cite{tan2010efficient,markkanen2019cauchy,churchill2022sampling,uribe2022hybrid,uribe2023horseshoe}.
Notably, using a Gibbs sampler is no longer possible if the forward operator is non-linear, the noise is not additive and Gaussian, or any other generalized gamma hyper-prior is used.
%
At the same time, we leverage the geometry of the prior-normalized posterior, which has a standard normal prior, by using the ES sampler \cite{murray2010elliptical}. 
The ES sampler is an MCMC method designed for sampling from posteriors with a Gaussian prior by making proposals along elliptical contours defined by the prior. 
It efficiently explores the posterior by drawing an auxiliary direction from the prior and rotating along the ellipse until a point within the slice (i.e., above a likelihood threshold) is found, avoiding the need for step size tuning. 

\Cref{fig:impulse} displays the means of the $x$- and $\theta$-samples for the original and prior-normalized posteriors.
We run four MCMC chains in parallel, generating $10^5$ samples per chain using the Gibbs sampler for the original posterior and $10^4$ samples per chain using the ES sampler for the prior-normalized posterior.
The total runtime for sampling was approximately 20 minutes for the Gibbs sampler and 3 minutes for the ES sampler.\footnote{We used a 2019 MacBook Pro with a 2.6 GHz 6-core Intel Core i7 processor and 16 GB of RAM.} 
Although the ES sampler incurs higher per-sample computational cost than the Gibbs sampler, it explores the posterior more efficiently.
This is evident in \cref{fig:impulse_n20_model4_x_mean_Gibbs_original,fig:impulse_n20_model4_x_mean_ESS_priorNormalized}, where the sample means obtained via the ES sampler are visibly more accurate, despite requiring only about 15\% of the Gibbs sampler's runtime.
In particular, while the support of the impulse image is identifiable in \cref{fig:impulse_n20_model4_x_mean_Gibbs_original}, the reconstructed values are substantially underestimated---for example, the maximum sample mean is on the order of $10^{-2}$, whereas the corresponding true value is approximately 2. 
That is, although the Gibbs sampler is frequently used for hierarchical Bayesian models \cite{tan2010efficient,markkanen2019cauchy,churchill2022sampling,uribe2022hybrid, uribe2023horseshoe}, it may fail to explore the original posterior adequately.
In contrast, existing MCMC samplers---in this case, the ES sampler---demonstrate improved exploration efficiency when applied to the proposed prior-normalized posterior. 
\section{Summary} 
\label{sec:summary} 

We introduced hierarchical prior normalization to transform challenging high-dimensional SBL posteriors, arising from sparsity-promoting hierarchical SBL priors, into prior-normalized posteriors with standard normal priors.
We derived the desired prior-normalizing TMs analytically using the product-like structure of SBL priors and KR rearrangements.
Our numerical experiments, including signal deblurring, inferring the initial condition of the nonlinear inviscid Burgers equation, and recovering an impulse image from noisy DCT data, show that standard MCMC techniques sample the prior-normalized posterior more efficiently than the original one. 
Specifically, even in settings where the Gibbs sampler---widely used in hierarchical Bayesian models due to its relatively low per-sample cost \cite{tan2010efficient,markkanen2019cauchy,churchill2022sampling,uribe2022hybrid,uribe2023horseshoe}---is applicable, it may fail to explore the original posterior adequately. 
In contrast, we found existing MCMC samplers, such as the ES sampler, to exhibit substantially improved exploration efficiency when applied to the proposed prior-normalized posterior.

Although we focused on SBL priors, the same approach applies to a broader class of hierarchical priors based on scale mixtures of normals \cite{beale1959scale,andrews1974scale}, including Laplace \cite{west1987scale,park2008bayesian,flock2024continuous} and horseshoe priors \cite{carvalho2009handling,uribe2023horseshoe,dong2023inducing}, which we will investigate in future work.
We also plan to combine prior normalization with samplers tailored to standard normal priors, such as dimension-robust methods \cite{cotter2013mcmc,pinski2015algorithms,beskos2017geometric}, different variants of the ES sampler \cite{murray2010elliptical,nishihara2014parallel,cabezas2023transport}, and dimension-independent likelihood-informed (DILI) approaches \cite{zahm2022certified,cui2022prior}.
Additionally, we aim to develop a theoretical analysis of the geometry of the prior-normalized posterior to inform and improve sampling strategies.

\section*{Acknowledgements} 
We thank Mathieu Le Provost, Jonathan Lindbloom, and Daniel Sharp for helpful feedback.


\appendix 
\section{Stable implementation of $t^{\tau}$} 
\label{app:implementation}
We employ Julia \cite{bezanson2017julia}, which predominantly uses double precision for floating-point computations.
A naive implementation of \cref{eq:t1}---first computing the CDF $\Phi^0$ and then the inverse CDF (quantile function) $( \mathcal{P}^{\theta} )^{-1}$---results in accumulation of numerical errors. 
This is because $\Phi^0(\tau) = 1$ for $\tau \geq \tau^*$, where $\tau^* \approx 8$. 
To reduce the numerical errors, we use the following stable implementation of $t^{\tau}$. 

Let $F(\theta; r, \beta, \vartheta)$ and $F^{-1}(z; r, \beta, \vartheta)$ respectively be the CDF and the quantile function of the generalized gamma distribution $\mathcal{GG}(r, \beta, \vartheta)$. 
To evaluate $t^{\tau}(\tau)$ in \cref{eq:t1}, we have to compute 
\begin{equation}\label{eq:rem_t_aux1} 
	t^{\tau}(\tau) = F^{-1}( \Phi^{0}(\tau); r, \beta, \vartheta),
\end{equation}
where $\Phi^{0}$ is the CDF of the standard normal distribution. 
Recall that (\cite{gomes2008parameter}) 
\begin{equation}
	F( \theta; r, \beta, \vartheta) = 
	\begin{cases}
		G\left( [\theta/\vartheta]^r; \beta, 1 \right) & \text{if } r > 0, \\ 
		1 - G\left( [\theta/\vartheta]^r; \beta, 1 \right) & \text{if } r < 0,
	\end{cases}
\end{equation}
where $G( \theta; \beta, 1 )$ is the CDF of the gamma distribution with density $g( \theta; \beta, 1 ) \propto \theta^{\beta - 1} e^{- \theta}$. 
Consequently, for the quantile function, we have 
\begin{equation}
	F^{-1}( z; r, \beta, \vartheta) = 
	\begin{cases}
		\vartheta \left[ G^{-1}( z; \beta, 1 ) \right]^{1/r} & \text{if } r > 0, \\ 
		\vartheta \left[ G^{-1}( 1-z; \beta, 1 ) \right]^{1/r} & \text{if } r < 0,
	\end{cases}
\end{equation}
where $G^{-1}( z; \beta, 1 )$ is the quantile function of the gamma distribution, which implies 
\begin{equation}\label{eq:rem_t_aux2}
	t^{\tau}(\tau) = 
	\begin{cases}
		\vartheta \left[ G^{-1}( \Phi^{0}(\tau); \beta, 1 ) \right]^{1/r} & \text{if } r > 0, \\ 
		\vartheta \left[ G^{-1}( 1-\Phi^{0}(\tau); \beta, 1 ) \right]^{1/r} & \text{if } r < 0.
	\end{cases}
\end{equation} 
Now recall that $G( z; \beta, 1 ) = \frac{1}{\Gamma(\beta)} \int_0^z u^{\beta - 1} e^{-u} \intd u = P(z; \beta)$, where $P(z; \beta)$ is the incomplete gamma function ratio; see \cite[\href{https://dlmf.nist.gov/8.2}{Section 8.2}]{NIST:DLMF}.
We can invert $P(z; \beta)$ using SpecialFunctions.jl (or StatsFuns.jl),
which transforms \cref{eq:rem_t_aux2} into 
 \begin{equation}\label{eq:rem_t_aux3}
	t^{\tau}(\tau) = \vartheta \left[ \text{gamma\_inc\_inv}( \beta, p, q ) \right]^{1/r}
\end{equation}
with $p = 1 - \Phi^{0}(\tau)$ if $r>0$, $p = \Phi^{0}(\tau)$ if $r<0$, and $q=1-p$. 
Notably, since $\Phi^{0}(\tau) = \frac{1}{2} [ 1 + \rm{erf}( \tau / \sqrt{2} ) ]$, we can stably compute $p = 1 - \Phi^{0}(\tau)$ as $p = \frac{1}{2} \rm{erfc}( \tau / \sqrt{2} )$. 
Here, $\rm{erf}$ and $\rm{erfc} = 1 - \,\rm{erf}$ are the error and complementary error function, respectively. 


\bibliographystyle{siamplain}
\bibliography{references}

@article{kaipio2007statistical,
  title={Statistical inverse problems: discretization, model reduction and inverse crimes},
  author={Kaipio, Jari and Somersalo, Erkki},
  journal={Journal of Computational and Applied Mathematics},
  volume={198},
  number={2},
  pages={493--504},
  year={2007},
  publisher={Elsevier}
}

@article{stuart2010inverse,
  title={Inverse problems: a {B}ayesian perspective},
  author={Stuart, Andrew M},
  journal={Acta Numerica},
  volume={19},
  pages={451--559},
  year={2010},
  publisher={Cambridge University Press}
}

@book{calvetti2023bayesian,
  title={Bayesian Scientific Computing},
  author={Calvetti, Daniela and Somersalo, Erkki},
  volume={215},
  year={2023},
  publisher={Springer Nature}
}

@article{tipping2001sparse,
  title={Sparse {B}ayesian learning and the relevance vector machine},
  author={Tipping, Michael E},
  journal={Journal of Machine Learning Research},
  volume={1},
  number={Jun},
  pages={211--244},
  year={2001}
}

@article{calvetti2007gaussian,
  title={A {G}aussian hypermodel to recover blocky objects},
  author={Calvetti, Daniela and Somersalo, Erkki},
  journal={Inverse Problems},
  volume={23},
  number={2},
  pages={733},
  year={2007},
  publisher={IOP Publishing}
}

@article{calvetti2019hierachical,
  title={Hierachical {B}ayesian models and sparsity: $\ell_2$-magic},
  author={Calvetti, Daniela and Somersalo, Erkki and Strang, A},
  journal={Inverse Problems},
  volume={35},
  number={3},
  pages={035003},
  year={2019},
  publisher={IOP Publishing}
}

@article{calvetti2020sparse,
  title={Sparse reconstructions from few noisy data: analysis of hierarchical {B}ayesian models with generalized gamma hyperpriors},
  author={Calvetti, Daniela and Pragliola, Monica and Somersalo, Erkki and Strang, Alexander},
  journal={Inverse Problems},
  volume={36},
  number={2},
  pages={025010},
  year={2020},
  publisher={IOP Publishing}
}

@article{calvetti2020sparsity,
  title={Sparsity promoting hybrid solvers for hierarchical {B}ayesian inverse problems},
  author={Calvetti, Daniela and Pragliola, Monica and Somersalo, Erkki},
  journal={SIAM Journal on Scientific Computing},
  volume={42},
  number={6},
  pages={A3761--A3784},
  year={2020},
  publisher={SIAM}
}

@article{glaubitz2023generalized,
  title={Generalized sparse {B}ayesian learning and application to image reconstruction},
  author={Glaubitz, Jan and Gelb, Anne and Song, Guohui},
  journal={SIAM/ASA Journal on Uncertainty Quantification},
  volume={11},
  number={1},
  pages={262--284},
  year={2023},
  publisher={SIAM}
}

@article{xiao2023sequential,
  title={Sequential image recovery using joint hierarchical {B}ayesian learning},
  author={Xiao, Yao and Glaubitz, Jan},
  journal={Journal of Scientific Computing},
  volume={96},
  number={1},
  pages={4},
  year={2023},
  publisher={Springer}
}

@article{sanz2024hierarchical,
  title={Hierarchical {B}ayesian inverse problems: A high-dimensional statistics viewpoint},
  author={Sanz-Alonso, Daniel and Waniorek, Nathan},
  journal={SIAM Review},
  volume={67},
  number={3},
  pages={543--575},
  year={2025},
  publisher={SIAM}
}

@article{lindbloom2024generalized,
  title={Efficient sparsity-promoting {MAP} estimation for {B}ayesian linear inverse problems},
  author={Lindbloom, Jonathan and Glaubitz, Jan and Gelb, Anne},
  journal={Inverse Problems},
  volume={41},
  number={2},
  pages={025001},
  year={2025},
  publisher={IOP Publishing}
}

@article{glaubitz2024leveraging,
  title={Leveraging joint sparsity in hierarchical {B}ayesian learning},
  author={Glaubitz, Jan and Gelb, Anne},
  journal={SIAM/ASA Journal on Uncertainty Quantification},
  volume={12},
  number={2},
  pages={442--472},
  year={2024},
  publisher={SIAM}
}

@article{si2024path,
  title={Path-following methods for Maximum a Posteriori estimators in {B}ayesian hierarchical models: How estimates depend on hyperparameters},
  author={Si, Zilai and Liu, Yucong and Strang, Alexander},
  journal={SIAM Journal on Optimization},
  volume={34},
  number={3},
  pages={2201--2230},
  year={2024},
  publisher={SIAM}
}

@article{lindbloom2025priorconditioned,
  title={Priorconditioned Sparsity-Promoting Projection Methods for Deterministic and {B}ayesian Linear Inverse Problems},
  author={Lindbloom, Jonathan and Pasha, Mirjeta and Glaubitz, Jan and Marzouk, Youssef},
  journal={arXiv preprint arXiv:2505.01827},
  year={2025}
}

@article{beale1959scale,
  title={Scale mixing of symmetric distributions with zero means},
  author={Beale, Evelyn Martin L and Mallows, Colin L},
  journal={The Annals of Mathematical Statistics},
  pages={1145--1151},
  year={1959},
  publisher={JSTOR}
}

@article{andrews1974scale,
  title={Scale mixtures of normal distributions},
  author={Andrews, David F and Mallows, Colin L},
  journal={Journal of the Royal Statistical Society: Series B (Methodological)},
  volume={36},
  number={1},
  pages={99--102},
  year={1974},
  publisher={Wiley Online Library}
}

@article{figueiredo2007majorization,
  title={Majorization--minimization algorithms for wavelet-based image restoration},
  author={Figueiredo, M{\'a}rio AT and Bioucas-Dias, Jos{\'e} M and Nowak, Robert D},
  journal={IEEE Transactions on Image Processing},
  volume={16},
  number={12},
  pages={2980--2991},
  year={2007},
  publisher={IEEE}
}

@article{babacan2009bayesian,
  title={Bayesian compressive sensing using {L}aplace priors},
  author={Babacan, S Derin and Molina, Rafael and Katsaggelos, Aggelos K},
  journal={IEEE Transactions on Image Processing},
  volume={19},
  number={1},
  pages={53--63},
  year={2009},
  publisher={IEEE}
}

@inproceedings{carvalho2009handling,
  title={Handling sparsity via the horseshoe},
  author={Carvalho, Carlos M and Polson, Nicholas G and Scott, James G},
  booktitle={Artificial Intelligence and Statistics},
  pages={73--80},
  year={2009},
  organization={PMLR}
}

@article{markkanen2019cauchy,
  title={Cauchy difference priors for edge-preserving {B}ayesian inversion},
  author={Markkanen, Markku and Roininen, Lassi and Huttunen, Janne MJ and Lasanen, Sari},
  journal={Journal of Inverse and Ill-Posed Problems},
  volume={27},
  number={2},
  pages={225--240},
  year={2019}
}

@article{suuronen2022cauchy,
  title={Cauchy {M}arkov random field priors for {B}ayesian inversion},
  author={Suuronen, Jarkko and Chada, Neil K and Roininen, Lassi},
  journal={Statistics and Computing},
  volume={32},
  number={2},
  pages={33},
  year={2022},
  publisher={Springer}
}

@article{dong2023inducing,
  title={Inducing sparsity via the horseshoe prior in imaging problems},
  author={Dong, Yiqiu and Pragliola, Monica},
  journal={Inverse Problems},
  volume={39},
  number={7},
  pages={074001},
  year={2023},
  publisher={IOP Publishing}
}

@article{papaspiliopoulos2003non,
  title={Non-centered parameterisations for hierarchical models and data augmentation},
  author={Papaspiliopoulos, Omiros and Roberts, Gareth O and Sk{\"o}ld, Martin},
  journal={Bayesian Statistics},
  volume={7},
  pages={307--326},
  year={2003},
  publisher={Oxford University Press New York}
}

@article{papaspiliopoulos2007general,
  title={A general framework for the parametrization of hierarchical models},
  author={Papaspiliopoulos, Omiros and Roberts, Gareth O and Sk{\"o}ld, Martin},
  journal={Statistical Science},
  pages={59--73},
  year={2007},
  publisher={JSTOR}
}

@article{chada2018parameterizations,
  title={Parameterizations for ensemble {K}alman inversion},
  author={Chada, Neil K and Iglesias, Marco A and Roininen, Lassi and Stuart, Andrew M},
  journal={Inverse Problems},
  volume={34},
  number={5},
  pages={055009},
  year={2018},
  publisher={IOP Publishing}
}

@article{dunlop2020hyperparameter,
  title={Hyperparameter estimation in {B}ayesian {MAP} estimation: parameterizations and consistency},
  author={Dunlop, Matthew M and Helin, Tapio and Stuart, Andrew M},
  journal={The SMAI Journal of Computational Mathematics},
  volume={6},
  pages={69--100},
  year={2020}
}

@article{rosenblatt1952remarks,
  title={Remarks on a multivariate transformation},
  author={Rosenblatt, Murray},
  journal={The Annals of Mathematical Statistics},
  volume={23},
  number={3},
  pages={470--472},
  year={1952},
  publisher={JSTOR}
}

@article{knothe1957contributions,
  title={Contributions to the theory of convex bodies.},
  author={Knothe, Herbert},
  journal={Michigan Mathematical Journal},
  volume={4},
  number={1},
  pages={39--52},
  year={1957},
  publisher={University of Michigan, Department of Mathematics}
}

@article{bogachev2005triangular,
  title={Triangular transformations of measures},
  author={Bogachev, Vladimir Igorevich and Kolesnikov, Aleksandr Viktorovich and Medvedev, Kirill Vladimirovich},
  journal={Sbornik: Mathematics},
  volume={196},
  number={3},
  pages={309},
  year={2005},
  publisher={IOP Publishing}
}

@book{villani2009optimal,
  title={Optimal Transport: Old and New},
  author={Villani, C{\'e}dric},
  volume={338},
  year={2009},
  publisher={Springer}
}

@article{ambrosio2013user,
  title={A User’s Guide to Optimal Transport},
  author={Ambrosio, Luigi and Gigli, Nicola},
  journal={Modelling and Optimisation of Flows on Networks},
  pages={1},
  year={2013}
}

@book{santambrogio2015optimal,
  title={Optimal Transport for Applied Mathematicians},
  author={Santambrogio, Filippo},
  year={2015},
  publisher={Springer}
}

@inproceedings{rezende2015variational,
  title={Variational inference with normalizing flows},
  author={Rezende, Danilo and Mohamed, Shakir},
  booktitle={International Conference on Machine Learning},
  pages={1530--1538},
  year={2015},
  organization={PMLR}
}

@article{marzouk2016sampling,
  title={Sampling via measure transport: An introduction},
  author={Marzouk, Youssef and Moselhy, Tarek and Parno, Matthew and Spantini, Alessio},
  journal={Handbook of Uncertainty Quantification},
  volume={1},
  pages={2},
  year={2016},
  publisher={Springer Cham}
}

@article{parno2018transport,
  title={Transport map accelerated {M}arkov chain {M}onte {C}arlo},
  author={Parno, Matthew D and Marzouk, Youssef M},
  journal={SIAM/ASA Journal on Uncertainty Quantification},
  volume={6},
  number={2},
  pages={645--682},
  year={2018},
  publisher={SIAM}
}

@article{baptista2023representation,
  title={On the representation and learning of monotone triangular transport maps},
  author={Baptista, Ricardo and Marzouk, Youssef and Zahm, Olivier},
  journal={Foundations of Computational Mathematics},
  pages={1--46},
  year={2023},
  publisher={Springer}
}

@inproceedings{fleischer2007transformations,
  title={Transformations for accelerating {MCMC} simulations with broken ergodicity},
  author={Fleischer, Mark},
  booktitle={2007 Winter Simulation Conference},
  pages={658--666},
  year={2007},
  organization={IEEE}
}

@article{wang2017bayesian,
  title={Bayesian inverse problems with $\ell_1$ priors: a randomize-then-optimize approach},
  author={Wang, Zheng and Bardsley, Johnathan M and Solonen, Antti and Cui, Tiangang and Marzouk, Youssef M},
  journal={SIAM Journal on Scientific Computing},
  volume={39},
  number={5},
  pages={S140--S166},
  year={2017},
  publisher={SIAM}
}

@article{chen2018robust,
  title={Robust {MCMC} sampling with non-{G}aussian and hierarchical priors in high dimensions},
  author={Chen, Victor and Dunlop, Matthew M and Papaspiliopoulos, Omiros and Stuart, Andrew M},
  journal={arXiv preprint arXiv:1803.03344},
  volume={3},
  year={2018}
}

@article{cui2022prior,
  title={Prior normalization for certified likelihood-informed subspace detection of {B}ayesian inverse problems},
  author={Cui, Tiangang and Tong, Xin T and Zahm, Olivier},
  journal={Inverse Problems},
  volume={38},
  number={12},
  pages={124002},
  year={2022},
  publisher={IOP Publishing}
}

@article{gelman1992inference,
  title={Inference from iterative simulation using multiple sequences},
  author={Gelman, Andrew and Rubin, Donald B},
  journal={Statistical Science},
  volume={7},
  number={4},
  pages={457--472},
  year={1992},
  publisher={Institute of Mathematical Statistics}
}

@article{brooks1998general,
  title={General methods for monitoring convergence of iterative simulations},
  author={Brooks, Stephen P and Gelman, Andrew},
  journal={Journal of Computational and Graphical Statistics},
  volume={7},
  number={4},
  pages={434--455},
  year={1998},
  publisher={Taylor \& Francis}
}

@book{gelman2003bayesian,
  title={Bayesian Data Analysis},
  author={Gelman, Andrew and Carlin, John B and Stern, Hal S and Rubin, Donald B},
  year={2003},
  publisher={Chapman and Hall/CRC}
}

@article{wolff2004monte,
  title={Monte {C}arlo errors with less errors},
  author={Wolff, Ulli and Alpha Collaboration},
  journal={Computer Physics Communications},
  volume={156},
  number={2},
  pages={143--153},
  year={2004},
  publisher={Elsevier}
}

@book{brooks2011handbook,
  title={Handbook of Markov Chain Monte Carlo},
  author={Brooks, Steve and Gelman, Andrew and Jones, Galin and Meng, Xiao-Li},
  year={2011},
  publisher={CRC press}
}

@book{liu2013monte,
  title={Monte Carlo Strategies in Scientific Computing},
  author={Liu, Jun S},
  year={2013},
  publisher={Springer Science \& Business Media}
}

@book{robert2013monte,
  title={Monte Carlo Statistical Methods},
  author={Robert, Christian and Casella, George},
  year={2013},
  publisher={Springer Science \& Business Media}
}

@article{cotter2013mcmc,
  title={{MCMC} Methods for Functions: Modifying Old Algorithms to Make Them Faster},
  author={Cotter, SL and Roberts, GO and Stuart, AM and White, D},
  journal={Statistical Science},
  volume={28},
  number={3},
  pages={424--446},
  year={2013}
}

@article{pinski2015algorithms,
  title={Algorithms for {K}ullback--{L}eibler approximation of probability measures in infinite dimensions},
  author={Pinski, Frank J and Simpson, Gideon and Stuart, Andrew M and Weber, Hendrik},
  journal={SIAM Journal on Scientific Computing},
  volume={37},
  number={6},
  pages={A2733--A2757},
  year={2015},
  publisher={SIAM}
}

@article{fox2016fast,
  title={Fast sampling in a linear-{G}aussian inverse problem},
  author={Fox, Colin and Norton, Richard A},
  journal={SIAM/ASA Journal on Uncertainty Quantification},
  volume={4},
  number={1},
  pages={1191--1218},
  year={2016},
  publisher={SIAM}
}

@article{beskos2017geometric,
  title={Geometric {MCMC} for infinite-dimensional inverse problems},
  author={Beskos, Alexandros and Girolami, Mark and Lan, Shiwei and Farrell, Patrick E and Stuart, Andrew M},
  journal={Journal of Computational Physics},
  volume={335},
  pages={327--351},
  year={2017},
  publisher={Elsevier}
}

@article{vehtari2021rank,
  title={Rank-normalization, folding, and localization: An improved $\hat{R}$ for assessing convergence of {MCMC} (with discussion)},
  author={Vehtari, Aki and Gelman, Andrew and Simpson, Daniel and Carpenter, Bob and B{\"u}rkner, Paul-Christian},
  journal={Bayesian Analysis},
  volume={16},
  number={2},
  pages={667--718},
  year={2021},
  publisher={International Society for Bayesian Analysis}
}

@article{sanz2024first,
  title={A First Course in {M}onte {C}arlo Methods},
  author={Sanz-Alonso, Daniel and Al-Ghattas, Omar},
  journal={arXiv preprint arXiv:2405.16359},
  year={2024}
}

@article{agapiou2014analysis,
  title={Analysis of the {G}ibbs sampler for hierarchical inverse problems},
  author={Agapiou, Sergios and Bardsley, Johnathan M and Papaspiliopoulos, Omiros and Stuart, Andrew M},
  journal={SIAM/ASA Journal on Uncertainty Quantification},
  volume={2},
  number={1},
  pages={511--544},
  year={2014},
  publisher={SIAM}
}

@article{churchill2022sampling,
  title={Sampling-based spotlight {SAR} image reconstruction from phase history data for speckle reduction and uncertainty quantification},
  author={Churchill, Victor and Gelb, Anne},
  journal={SIAM/ASA Journal on Uncertainty Quantification},
  volume={10},
  number={3},
  pages={1225--1249},
  year={2022},
  publisher={SIAM}
}

@article{uribe2022hybrid,
  title={A hybrid {G}ibbs sampler for edge-preserving tomographic reconstruction with uncertain view angles},
  author={Uribe, Felipe and Bardsley, Johnathan M and Dong, Yiqiu and Hansen, Per Christian and Riis, Nicolai AB},
  journal={SIAM/ASA Journal on Uncertainty Quantification},
  volume={10},
  number={3},
  pages={1293--1320},
  year={2022},
  publisher={SIAM}
}

@article{uribe2023horseshoe,
  title={Horseshoe priors for edge-preserving linear {B}ayesian inversion},
  author={Uribe, Felipe and Dong, Yiqiu and Hansen, Per Christian},
  journal={SIAM Journal on Scientific Computing},
  volume={45},
  number={3},
  pages={B337--B365},
  year={2023},
  publisher={SIAM}
}

@article{ascolani2024scalability,
  title={Scalability of {M}etropolis-within-{G}ibbs schemes for high-dimensional {B}ayesian models},
  author={Ascolani, Filippo and Roberts, Gareth O and Zanella, Giacomo},
  journal={arXiv preprint arXiv:2403.09416},
  year={2024}
}

@inproceedings{tan2010efficient,
  title={Efficient sparse {B}ayesian learning via {G}ibbs sampling},
  author={Tan, Xing and Li, Jian and Stoica, Peter},
  booktitle={2010 IEEE International Conference on Acoustics, Speech and Signal Processing},
  pages={3634--3637},
  year={2010},
  organization={IEEE}
}

@article{calvetti2024computationally,
  title={Computationally efficient sampling methods for sparsity promoting hierarchical {B}ayesian models},
  author={Calvetti, Daniela and Somersalo, Erkki},
  journal={SIAM/ASA Journal on Uncertainty Quantification},
  volume={12},
  number={2},
  pages={524--548},
  year={2024},
  publisher={SIAM}
}

@article{west1987scale,
  title={On scale mixtures of normal distributions},
  author={West, Mike},
  journal={Biometrika},
  volume={74},
  number={3},
  pages={646--648},
  year={1987},
  publisher={Oxford University Press}
}

@article{park2008bayesian,
  title={The {B}ayesian {L}asso},
  author={Park, Trevor and Casella, George},
  journal={Journal of the American Statistical Association},
  volume={103},
  number={482},
  pages={681--686},
  year={2008},
  publisher={Taylor \& Francis}
}

@article{flock2024continuous,
  title={Continuous {G}aussian mixture solution for linear {B}ayesian inversion with application to {L}aplace priors},
  author={Flock, Rafael and Dong, Yiqiu and Uribe, Felipe and Zahm, Olivier},
  journal={Inverse Problems},
  year={2024}
}

@article{roberts1998optimal,
  title={Optimal scaling of discrete approximations to {L}angevin diffusions},
  author={Roberts, Gareth O and Rosenthal, Jeffrey S},
  journal={Journal of the Royal Statistical Society: Series B (Statistical Methodology)},
  volume={60},
  number={1},
  pages={255--268},
  year={1998},
  publisher={Wiley Online Library}
}

@article{haario2001adaptive,
  title={An adaptive {M}etropolis algorithm},
  author={Haario, Heikki and Saksman, Eero and Tamminen, Johanna},
  journal={Bernoulli},
  volume={7},
  number={6},
  pages={223--242},
  year={2001}
}

@article{andrieu2008tutorial,
  title={A tutorial on adaptive {MCMC}},
  author={Andrieu, Christophe and Thoms, Johannes},
  journal={Statistics and Computing},
  volume={18},
  pages={343--373},
  year={2008},
  publisher={Springer}
}

@article{atchade2010limit,
  title={Limit theorems for some adaptive {MCMC} algorithms with subgeometric kernels},
  author={Atchad{\'e}, Yves and Fort, Gersende},
  journal={Bernoulli},
  volume={16},
  number={1},
  year={2010},
  publisher={Bernoulli Society for Mathematical Statistics and Probability}
}

@article{yang2020optimal,
  title={Optimal scaling of random-walk {M}etropolis algorithms on general target distributions},
  author={Yang, Jun and Roberts, Gareth O and Rosenthal, Jeffrey S},
  journal={Stochastic Processes and their Applications},
  volume={130},
  number={10},
  pages={6094--6132},
  year={2020},
  publisher={Elsevier}
}

@article{atchade2006adaptive,
  title={An adaptive version for the {M}etropolis adjusted {L}angevin algorithm with a truncated drift},
  author={Atchad{\'e}, Yves F},
  journal={Methodology and Computing in Applied Probability},
  volume={8},
  pages={235--254},
  year={2006},
  publisher={Springer}
}

@article{marshall2012adaptive,
  title={An adaptive approach to {L}angevin {MCMC}},
  author={Marshall, Tristan and Roberts, Gareth},
  journal={Statistics and Computing},
  volume={22},
  pages={1041--1057},
  year={2012},
  publisher={Springer}
}

@article{beskos2013optimal,
  title={Optimal tuning of the hybrid {M}onte {C}arlo algorithm},
  author={Beskos, Alexandros and Pillai, Natesh and Roberts, Gareth and Sanz-Serna, Jesus-Maria and Stuart, Andrew}, 
  journal={Bernoulli},
  volume={19},
  number={5A},
  pages={1501--1534},
  year={2013}
}

@article{zahm2022certified,
  title={Certified dimension reduction in nonlinear {B}ayesian inverse problems},
  author={Zahm, Olivier and Cui, Tiangang and Law, Kody and Spantini, Alessio and Marzouk, Youssef},
  journal={Mathematics of Computation},
  volume={91},
  number={336},
  pages={1789--1835},
  year={2022}
}

@article{bezanson2017julia,
  title={Julia: A fresh approach to numerical computing},
  author={Bezanson, Jeff and Edelman, Alan and Karpinski, Stefan and Shah, Viral B},
  journal={SIAM Review},
  volume={59},
  number={1},
  pages={65--98},
  year={2017},
  publisher={SIAM}
}

@article{rackauckas2017differentialequations,
  title={{DifferentialEquations.jl} {--} {A} Performant and Feature-Rich
         Ecosystem for Solving Differential Equations in {J}ulia},
  author={Rackauckas, Christopher and Nie, Qing},
  journal={Journal of Open Research Software},
  volume={5},
  number={1},
  pages={15},
  year={2017},
  publisher={Ubiquity Press}
}

@article{ranocha2022adaptive,
  title={Adaptive numerical simulations with {T}rixi.jl:
         {A} case study of {J}ulia for scientific computing},
  author={Ranocha, Hendrik and Schlottke-Lakemper, Michael and Winters, Andrew Ross
          and Faulhaber, Erik and Chan, Jesse and Gassner, Gregor J},
  journal={Proceedings of the JuliaCon Conferences},
  volume={1},
  number={1},
  pages={77},
  year={2022},
  month={01},
  publisher={The Open Journal}
}

@article{ranocha2025robustness,
  title={On the robustness of high-order upwind summation-by-parts methods for nonlinear conservation laws},
  author={Ranocha, Hendrik and Winters, Andrew R and Schlottke-Lakemper, Michael and {\"O}ffner, Philipp and Glaubitz, Jan and Gassner, Gregor J},
  journal={Journal of Computational Physics},
  volume={520},
  pages={113471},
  year={2025},
  publisher={Elsevier}
}

@inproceedings{murray2010elliptical,
  title={Elliptical slice sampling},
  author={Murray, Iain and Adams, Ryan and MacKay, David},
  booktitle={Proceedings of the 13th International Conference on Artificial Intelligence and Statistics},
  pages={541--548},
  year={2010}
}

@article{nishihara2014parallel,
  title={Parallel {MCMC} with generalized elliptical slice sampling},
  author={Nishihara, Robert and Murray, Iain and Adams, Ryan P},
  journal={The Journal of Machine Learning Research},
  volume={15},
  number={1},
  pages={2087--2112},
  year={2014},
  publisher={JMLR. org}
}

@inproceedings{cabezas2023transport,
  title={Transport elliptical slice sampling},
  author={Cabezas, Alberto and Nemeth, Christopher},
  booktitle={International Conference on Artificial Intelligence and Statistics},
  pages={3664--3676},
  year={2023},
  organization={PMLR}
}

@misc{NIST:DLMF,
    key = "{\relax DLMF}",
    title = "{\it NIST Digital Library of Mathematical Functions}",
howpublished = "\url{https://dlmf.nist.gov/}, Release 1.2.0 of 2024-03-15", 
    note = "F.~W.~J. Olver, A.~B. {Olde Daalhuis}, D.~W. Lozier, B.~I. Schneider, R.~F. Boisvert, C.~W. Clark, B.~R. Miller, B.~V. Saunders, H.~S. Cohl, and M.~A. McClain, eds."
}

@book{athreya2006measure,
  title={Measure Theory and Probability Theory},
  author={Athreya, Krishna B and Lahiri, Soumendra N},
  volume={19},
  year={2006},
  publisher={Springer}
}

@article{gomes2008parameter,
  title={Parameter estimation of the generalized gamma distribution},
  author={Gom{\`e}s, Oph{\'e}lie and Combes, Catherine and Dussauchoy, Alain},
  journal={Mathematics and Computers in Simulation},
  volume={79},
  number={4},
  pages={955--963},
  year={2008},
  publisher={Elsevier}
}

@book{nocedal2006numerical,
  title={Numerical Optimization},
  author={Nocedal, Jorge and Wright, Stephen J},
  year={2006},
  publisher={Springer}
}

\end{document}